# PROBLEM OF FINDING AN OPTIMAL PIECEWISE LINEAR PATH CONNECTING TWO GIVEN POINTS WITH THE POSSIBILITY OF MAKING $n$ TURNS

**Nefedov V.N.**

## Abstract

We consider the problem of finding an optimal piecewise linear path (polygonal line) connecting two given points $A, B \in \mathbb{R}^2$ with the possibility of making $n \geq 1$ turns at some points $B^{(1)}, \dots, B^{(n)} \in \mathbb{R}^2$ (the absolute value of each turn angle does not exceed a prescribed bound $\varphi \in (0, \pi/2)$). Under the condition $n\varphi < \pi$, we characterize the region to which all interior vertices of such a path must belong (Theorem 1). It is shown that for any point $B^{(1)}$ from this region, there exists a polygonal line satisfying the given constraints (Lemma 1). Based on these findings, an explicit expression is derived (Theorem 2) that describes the collection $\boldsymbol{S}^{(n)}(A, B, \varphi)$ of all admissible sequences $\left(B^{(1)}, \dots, B^{(n)}\right)$ of corner points. This expression is then used to construct a finite family of sequences $\left(B^{(1)}, \dots, B^{(n)}\right) \in \boldsymbol{S}^{(n)}(A, B, \varphi)$ that approximates the aforementioned collection. The resulting finite approximating family serves as the basis for developing algorithms that provide approximate solutions to an optimization problem over $\boldsymbol{S}^{(n)}(A, B, \varphi)$, where the objective function accounts for both the cost of traversing the segments and the cost associated with the turns.

**Keywords:** piecewise linear path, polygonal line, turn angle, circular segment, chord, convex cone, angular sector, approximation, optimization problem.

## Introduction.

There exist practically important problems in which it is required to connect two given points of the two-dimensional plane by a piecewise linear polygonal line. There are constraints on the number of segments of the polygonal line (for example, exactly $n$ segments) and on the turn angles at each point connecting adjacent segments. It is required to describe the set of admissible polygonal lines in order to be able to enumerate them algorithmically in the discrete case, or to approximate them by a finite set of polygonal lines in the continuous case (with subsequent enumeration). Having algorithms for such enumeration, one can then solve optimization problems for some objective function that takes into account the cost of traversing each segment of the polygonal line and the cost of making each turn. It is assumed that the values of the objective function are fairly easy to compute based on available information about the domain on which the optimal polygonal line is sought. In this paper we prove two main statements: Theorem 1 and Lemma 1, which are used to describe the set of admissible polygonal lines of interest.

Let $A, B \in \mathbb{R}^2$, $A \neq B$, $\varphi \in (0, \pi)$. We introduce the set

$$\boldsymbol{S}(A, B, \varphi) = \{\, C \in \mathbb{R}^2 \setminus \{A, B\} \mid \angle(B - C, C - A) \in (-\varphi, \varphi)\,\} \quad (0.1)$$

(we use the assumptions concerning the angle $\angle(C, E)$ for arbitrary $C, E \in \mathbb{R}^2$ according to [1]; more detailed information about angles is given below in Section 1.1). In [1], it is shown that for $\varphi \in (0, \pi)$ the set $\boldsymbol{S}(A, B, \varphi)$ is the union of two circular segments of circles with radius $r_\varphi = \frac{1}{\sin\varphi}$ (left and right parts), sharing the common chord $(A, B) = \{\alpha A + (1 - \alpha)B \mid \alpha \in (0,1)\}$. If $\varphi > \pi$, we also denote $\boldsymbol{S}(A, B, \varphi) = \mathbb{R}^2 \setminus \{A, B\}$ (see also Remark 0.1).

Furthermore, for any $E \in \mathbb{R}^2$ with $E \neq (0,0)$, in the case $\varphi \in (0, \pi/2)$ we denote by $\boldsymbol{C}(E, \varphi) = \{C \in \mathbb{R}^2 \mid \angle(C, E) \in [-\varphi, \varphi]\}$ a convex cone. If (as in [1]) we assume that $\forall C \in \mathbb{R}^2$ $\angle(C,(0,0)) = \angle((0,0),C) = 0$, then $(0,0) \in \boldsymbol{C}(E,\varphi)$. One can also consider $\boldsymbol{C}(E,\varphi)$ for $\varphi \geq \pi/2$. Note that

$$\boldsymbol{C}(E,\pi/2) = \{C \in \mathbb{R}^2 \mid \angle(C,E) \in [-\pi/2,\pi/2]\} = \{C \in \mathbb{R}^2 \mid \langle C,E\rangle \geq 0\}$$

is a half-plane, and for $\varphi \in (\pi/2,\pi)$ the set $\boldsymbol{C}(E,\varphi)$ ceases to be convex. For $\varphi \geq \pi$ we have $\boldsymbol{C}(E,\varphi) = \boldsymbol{C}(E,\pi) = \mathbb{R}^2$. Along with the open interval $(A,B)$ we also consider the closed segment $[A,B] = \{\alpha A + (1-\alpha)B \mid \alpha \in [0,1]\}$, where $A,B \in \mathbb{R}^2$, $A \neq B$.

**Remark 0.1.** One can consider the set $\boldsymbol{S}(A,B,\varphi)$ satisfying (0.1) also for $\varphi \geq \pi$. If $\varphi > \pi$, then $\boldsymbol{S}(A,B,\varphi) = \mathbb{R}^2 \setminus \{A,B\}$ because (see [1], as well as Section 1.1)

$$\forall C \in \mathbb{R}^2 \quad \angle(B-C,C-A) \in (-\pi,\pi] \subset (-\varphi,\varphi),$$

which agrees with the earlier choice for the notation $\boldsymbol{S}(A,B,\varphi)$ when $\varphi > \pi$. The case $\boldsymbol{S}(A,B,\pi)$ is more delicate. To describe it we need some notation. Let

$$A(t) = A + t(B-A), t \in \mathbb{R},$$
$$(A(-\infty),A) = \{A(t) \mid t \in (-\infty,0)\}, (A(-\infty),A] = \{A(t) \mid t \in (-\infty,0]\},$$
$$(B,A(+\infty)) = \{A(t) \mid t \in (1,+\infty)\}, [B,A(+\infty)) = \{A(t) \mid t \in [1,+\infty)\}.$$

Here $A(0) = A, A(1) = B$,

$$\forall C \in [A,B] \quad \angle(B-C,C-A) = 0,$$
$$\forall C \in (A(-\infty),A) \cup \big(B,A(+\infty)\big) \quad \angle(B-C,C-A) = \pi,$$
$$\forall C \notin (A(-\infty),A) \cup [A,B] \cup \big(B,A(+\infty)\big) \quad \angle(B-C,C-A) \in (-\pi,\pi),$$

hence

$$\boldsymbol{S}(A,B,\pi) = \mathbb{R}^2 \setminus ((A(-\infty),A] \cup [B,A(+\infty))). \qquad (0.2)$$

Thus, for $\varphi = \pi$ we define $\boldsymbol{S}(A,B,\varphi)$ according to (0.2).

We will also need the closure of the set $\boldsymbol{S}(A,B,\varphi)$:

$$\text{cl } \boldsymbol{S}(A,B,\varphi) = \{C \in \mathbb{R}^2 \mid \angle(B-C,C-A) \in [-\varphi,\varphi]\}.$$

We note some properties of the set $\boldsymbol{S}(A,B,\varphi)$.

**Property 1.** Let $A, B \in \mathbb{R}^2, A \neq B, \varphi \in (0, \pi)$. Then

$$\boldsymbol{S}(A, B, \varphi) = \boldsymbol{S}(B, A, \varphi).$$

Indeed, using the facts (see [1], as well as Section 1.1) that $\forall D, E \in \mathbb{R}^2$ $\angle(D, E) = -\angle(E, D), \angle(D, E) = \angle(-D, -E)$, we have

$$\angle(B - C, C - A) \in (-\varphi, \varphi) \iff -\angle(B - C, C - A) \in (-\varphi, \varphi)$$
$$\iff \angle(C - A, B - C) \in (-\varphi, \varphi) \iff \angle(A - C, C - B) \in (-\varphi, \varphi),$$

hence

$$\boldsymbol{S}(A, B, \varphi) = \{C \in \mathbb{R}^2 \setminus \{A, B\} \mid \angle(A - C, C - B) \in (-\varphi, \varphi)\} = \boldsymbol{S}(B, A, \varphi).$$

**Property 2.** Let $A, B \in \mathbb{R}^2, A \neq B, \bar{A}, \bar{B} \in (A, B), \varphi \in (0, \pi)$. Then

$$\boldsymbol{S}(\bar{A}, B, \varphi) \subseteq \boldsymbol{S}(A, B, \varphi), \boldsymbol{S}(A, \bar{B}, \varphi) \subseteq \boldsymbol{S}(A, B, \varphi),$$

and if $\bar{A} \neq \bar{B}$, then

$$\boldsymbol{S}(\bar{A}, \bar{B}, \varphi) \subseteq \boldsymbol{S}(A, B, \varphi).$$

It suffices to justify the case $\boldsymbol{S}(\bar{A}, B, \varphi) \subseteq \boldsymbol{S}(A, B, \varphi)$ (the case $\boldsymbol{S}(A, \bar{B}, \varphi) \subseteq \boldsymbol{S}(A, B, \varphi)$ is proved similarly, and $\boldsymbol{S}(\bar{A}, \bar{B}, \varphi) \subseteq \boldsymbol{S}(A, B, \varphi)$ follows easily from these two cases together with Property 1). Indeed, if $A, B, C$ are pairwise distinct points, $\bar{A} \in (A, B)$, then one can easily show (using simple geometric considerations) that in the case $\alpha = \angle(B - C, C - A) > 0$ we have $\bar{\alpha} = \angle(B - C, C - \bar{A}) > \alpha$, and in the case $\alpha = \angle(B - C, C - A) < 0$ we have $\bar{\alpha} = \angle(B - C, C - \bar{A}) < \alpha$, which implies that $C \in \boldsymbol{S}(\bar{A}, B, \varphi) \Rightarrow C \in \boldsymbol{S}(A, B, \varphi)$.

We will use the following main statements (announced in [2]).

**Theorem 1.** *Let $\varphi \in (0, \pi), n \geq 1, n\varphi < \pi, A = B^{(0)}, B = B^{(n+1)}, A \neq B$, $B^{(0)}, B^{(1)}, \ldots, B^{(n)}, B^{(n+1)}$ a sequence of points in $\mathbb{R}^2, B^{(1)} \neq A, B^{(n)} \neq B$,*

$$\boldsymbol{L} = \left[B^{(0)}, B^{(1)}\right] \cup \left[B^{(1)}, B^{(2)}\right] \cup \cdots \cup \left[B^{(n)}, B^{(n+1)}\right] \qquad (0.3)$$

*a polygonal line, and*

$$\angle\left(B^{(i+1)} - B^{(i)}, B^{(i)} - B^{(i-1)}\right) \in [-\varphi, \varphi], i = 1, \dots, n, \qquad (0.4)$$
$$B^{(i)} \neq B^{(i+1)}, i = 0,1, \dots, n.$$

*Then*

1. *for* $n \geq 2$ *we have* $B^{(1)}, \dots, B^{(n)} \in \boldsymbol{S}(A, B, n\varphi)$*, and in particular* $B^{(i)} \notin \{A, B\}, i = 1, \dots, n$*. Moreover,* $B^{(i)} \neq B^{(j)}$ *for* $i, j \in \{0,1, \dots, n+1\}, i \neq j$.
2. *for* $n = 1$ *we have* $B^{(1)} \in \mathrm{cl}\ \boldsymbol{S}(A, B, \varphi) \setminus \{A, B\}$.

The case $n = 1$ is quite obvious. The proof of Theorem 1 for $n \geq 2$ will be given in Section 2.

**Remark 0.2.** Condition (0.4) allows the cases $\angle\left(B^{(i+1)} - B^{(i)}, B^{(i)} - B^{(i-1)}\right) = 0$; when this occurs, the point $B^{(i)}$ is not a turn point, i.e., each such case reduces the number of turn points $n$ by 1, which simplifies the problem. In this connection, the case with the additional condition

$$\angle\left(B^{(i+1)} - B^{(i)}, B^{(i)} - B^{(i-1)}\right) \neq 0 \qquad (0.5)$$

is of greatest interest. However, it should be noted that in some applied problems an additional constraint of the form

$$| B^{(i)} - B^{(i-1)} | \leq s, i = 1, \dots, n+1,$$

is often imposed, where $s > 0$ is a given number (a constraint on the “span” length). In such problems, a zero turn at some intermediate point $B^{(i)}$ is quite appropriate, and it can be used to ensure that the specified condition holds. However, in the present work, this additional constraint is not taken into account.

**Lemma 1.** *Under the assumptions of Theorem 1, a polygonal line* $\boldsymbol{L}$ *satisfying conditions* (0.3), (0.4), (0.5) *exists for every* $B^{(1)} \in \boldsymbol{S}(A, B, n\varphi)$*, and for* $n = 1$ *it exists for every* $B^{(1)} \in \mathrm{cl}\ \boldsymbol{S}(A, B, \varphi) \setminus \{A, B\}$.

The proof of Lemma 1 will be given in Section 3.

**Remark 0.3.** Continuing Remark 0.1, we make clarifications for Theorem 1 and Lemma 1 in the case $n\varphi \geq \pi$. By analogy with the proof of Lemma 1, one can show that for $n\varphi > \pi$ a polygonal line $L$ satisfying (0.3), (0.4), (0.5) exists for every $B^{(1)} \in \mathbb{R}^2 \setminus \{A, B\}$, which corresponds to the condition $B^{(1)} \in \boldsymbol{S}(A, B, n\varphi) = \mathbb{R}^2 \setminus \{A, B\}$ in this case. It should be noted, however, that the condition $B^{(i)} \neq B^{(j)}$ for $i, j \in \{0,1, \dots, n+1\}, i \neq j$ (see Theorem 1) may then fail. In the more delicate case $n\varphi = \pi, n \geq 2$, one can also show that a polygonal line $\boldsymbol{L}$ satisfying (0.3), (0.4), (0.5) exists for every $B^{(1)} \in \mathbb{R}^2 \setminus ((A(-\infty), A] \cup [B, A(+\infty)))$, which corresponds to the condition $B^{(1)} \in \boldsymbol{S}(A, B, n\varphi) = \boldsymbol{S}(A, B, \pi)$ (see (0.2)) for this case. For $n = 1$ we have $\operatorname{cl} \boldsymbol{S}(A, B, \pi) \setminus \{A, B\} = \mathbb{R}^2 \setminus \{A, B\}$ and, since

$$\forall C \in \mathbb{R}^2 \quad \angle(B - C, C - A) \in (-\pi, \pi],$$

the indicated polygonal line $\boldsymbol{L}$ exists for every $B^{(1)} \in \operatorname{cl} \ \boldsymbol{S}(A, B, \varphi) \setminus \{A, B\}$ also in this case. Thus, Lemma 1 can be generalized to any value $n\varphi > 0$.

Let $\boldsymbol{S}^{(n)}(A, B, \varphi)$ denote the set of all possible sequences $\left(B^{(1)}, \dots, B^{(n)}\right)$ of interior points of a polygonal line $\boldsymbol{L}$ satisfying conditions (0.3), (0.4), where $A = B^{(0)}, B = B^{(n+1)}$. A consequence of Theorem 1 and Lemma 1 is

**Theorem 2.** *Under the assumptions of Theorem 1,*

$$\begin{gathered}\left(B^{(1)}, \dots, B^{(n)}\right) \in \mathbf{S}^{(n)}(A, B, \varphi) \Leftrightarrow \left(B^{(1)}, \dots, B^{(n)}\right) \in \\ \in \mathbf{S}(A, B, n\varphi) \times \left\{\mathbf{S}\left(B^{(1)}, B, (n-1)\varphi\right) \cap \left[B^{(1)} + \mathbf{C}\left(B^{(1)} - A, \varphi\right)\right]\right\} \times \\ \times \left\{\boldsymbol{S}\left(B^{(2)}, B, (n-2)\varphi\right) \cap \left[B^{(2)} + \mathbf{C}\left(B^{(2)} - B^{(1)}, \varphi\right)\right]\right\} \times \cdots \times \\ \times \left\{\mathbf{S}\left(B^{(n-2)}, B, 2\varphi\right) \cap \left[B^{(n-2)} + \mathbf{C}\left(B^{(n-2)} - B^{(n-3)}, \varphi\right)\right]\right\} \times \\ \times \left\{\left[\operatorname{cl} \mathbf{S}\left(B^{(n-1)}, B, \varphi\right) \setminus \{A, B\}\right] \cap \left[B^{(n-1)} + \mathbf{C}\left(B^{(n-1)} - B^{(n-2)}, \varphi\right)\right]\right\}. \quad (0.6)\end{gathered}$$

*and moreover, for any sequential choice of elements* $B^{(1)}, \dots, B^{(n)}$ *in accordance with the right-hand side of* (0.6), *the set from which the next element of this sequence is chosen will be nonempty for any collection of previous members.*

The proof of Theorem 2 will be given in Section 4.

**Remark 0.4.** One could add the natural condition (0.5) to the assumptions of Theorem 2, but this would complicate formula (0.6). In that case, from each $i$-th term of the Cartesian product, $i = 2, \dots, n$, one would have to exclude points lying on the line through $B^{(i-2)}, B^{(i-1)}$.

**Remark 0.5.** In view of Remark 0.3, formula (0.6) can also be applied in the case $n\varphi \geq \pi$, i.e., in all possible cases.

**Remark 0.6.** In formula (0.6), the principle of dynamic programming is essentially used [3-5].

Note that in the above reasoning no additional restrictions (apart from the conditions of Theorem 1 and the natural restriction (0.5)) were imposed on the points $A, B, B^{(1)}, \dots, B^{(n)}$ defining the polygonal line $\boldsymbol{L}$ satisfying (0.3), (0.4). Meanwhile, additional restrictions are possible, and there may be several of them.

One possible restriction is that the points $A, B, B^{(1)}, \dots, B^{(n)}$ belong to some admissible set $\boldsymbol{Q} \subseteq \mathbb{R}^2$. It is desirable that this set have a simple structure, for example a coordinate parallelepiped of the form

$$\boldsymbol{Q} = \{\, C = (c_1, c_2) \in \mathbb{R}^2 \mid q_i \leq c_i \leq p_i, = 1,2 \,\}, \qquad (0.7)$$

where $q_i, p_i \in \mathbb{R}$, $q_i < p_i$, $i = 1,2$.

Let $\tau > 0$. Consider the finite set

$$\boldsymbol{Q}^{(\tau)} = \{C = (c_1, c_2) \in \mathbb{R}^2 | c_i = j_i\tau, q_i \leq c_i \leq p_i, j_i \in \mathbb{Z}, i = 1,2\}. \qquad (0.8)$$

If we use the set $\boldsymbol{Q}^{(\tau)}$ instead of $\boldsymbol{Q}$, we arrive at the discrete case. When a constraint of the form $B^{(1)}, \dots, B^{(n)} \in \boldsymbol{Q}^{(\tau)}$ is present, it must be taken into account in each factor of the direct product (0.6), i.e., we now use $\boldsymbol{S}^{(n)}(\boldsymbol{Q}^{(\tau)}, A, B, \varphi) = \boldsymbol{S}^{(n)}(A, B, \varphi) \cap [\boldsymbol{Q}^{(\tau)}]^n$ instead of $\boldsymbol{S}^{(n)}(A, B, \varphi)$.

We propose to use the finite set $\boldsymbol{S}^{(n)}(\boldsymbol{Q}^{(\tau)}, A, B, \varphi)$ as an approximating set for $\boldsymbol{S}^{(n)}(A, B, \varphi)$ for sufficiently small $\tau > 0$. It is shown in Section 5 that if (0.7), (0.8) hold and, in addition, the condition $\boldsymbol{S}(A, B, n\varphi) \subset Q$ is satisfied, then the finite set $\boldsymbol{S}^{(n)}(\boldsymbol{Q}^{(\tau)}, A, B, \varphi)$ converges in the Hausdorff metric [6] to the set $\boldsymbol{S}^{(n)}(A, B, \varphi)$ as $\tau \to 0^+$.

Using the finite set $\boldsymbol{S}^{(n)}(\boldsymbol{Q}^{(\tau)}, A, B, \varphi)$, we can further consider discrete optimization problems of the form:

$$f(B_1, \dots, B_n) \to \min \,; (B_1, \dots, B_n) \in \boldsymbol{S}^{(n)}(\boldsymbol{Q}^{(\tau)}, A, B, \varphi), \qquad (0.9)$$

where $f(B_1, \dots, B_n)$ is an objective function that accounts for the cost of moving along the polygonal line $\boldsymbol{L}$ (see (0.3), (0.4)) as well as the cost of making turns at the points $B_1, \dots, B_n$ (by analogy with the objective function $f$ in [1]). We present some conditions for the convergence of solutions of the approximate problem (0.9) as $\tau \to 0^+$ to the solutions of the "exact" problem

$$f(B_1, \dots, B_n) \to \inf \;\; (= f^*); (B_1, \dots, B_n) \in \boldsymbol{S}^{(n)}(A, B, \varphi). \qquad (0.10)$$

If for some $W^* = (B_1^*, \dots, B_n^*) \in \boldsymbol{S}^{(n)}(A, B, \varphi)$ we have $f(W^*) = f^*$, then we call $W^*$ a **regular solution** of problem (0.10).

**Statement 0.1.** Suppose we are under the conditions of Theorem 1, $\boldsymbol{Q}, \boldsymbol{Q}^{(\tau)}$ satisfy (0.7), (0.8), and $\boldsymbol{S}(A, B, n\varphi) \subset \boldsymbol{Q}$. Then, if the function $f(B_1, \dots, B_n)$ is continuous on $\boldsymbol{S}^{(n)}(A, B, \varphi)$ and there exists at least one regular solution of problem (0.10), then

$$\lim_{\tau \to 0^+} \min \ f(\boldsymbol{S}^{(n)}(\boldsymbol{Q}^{(\tau)}, A, B, \varphi)) = f^*.$$

## 1. Preliminary information

### 1.1. Angle between vectors. Operations with angles

For arbitrary non-zero vectors $A = (a_1, a_2) \in \mathbb{R}^2$, $B = (b_1, b_2) \in \mathbb{R}^2$, the angle between these vectors (assuming they originate from the same point $\mathbf{0} = (0,0) \in \mathbb{R}^2$) will be denoted by $\angle(A, B)$. We take into account the orientation of the angle: from vector $A$ towards vector $B$, so that $\angle(A, B) = -\angle(B, A)$. The only exception is the case of oppositely directed vectors. In this case,

$$\angle(A, -A) = \angle(-A, A) = \pi.$$

The clockwise direction is considered positive. Angle values lie in the interval $(-\pi, \pi]$ (for the angle between oppositely directed vectors, we choose the value $\pi$). Furthermore, for greater generality of reasoning, in the case of a zero vector $A$ or $B$, we set $\angle(A, B) = 0$. It is easy to show that if $B \neq -tA$ for $t > 0$, then

$$\angle(A, B) > 0 \Leftrightarrow \det \begin{bmatrix} a_1 & a_2 \\ b_1 & b_2 \end{bmatrix} = \det \begin{bmatrix} a_1 & b_1 \\ a_2 & b_2 \end{bmatrix} < 0$$

(as a consequence of the fact that the cross product of vectors $(a_1, a_2, 0), (b_1, b_2, 0) \in \mathbb{R}^3$ forms a right-handed triple [7]). The only exception is the case $B, A \neq (0,0)$, $B = -tA$, $t > 0$. Then $\angle(A, B) = \pi > 0$, and

$$\det \begin{bmatrix} a_1 & a_2 \\ b_1 & b_2 \end{bmatrix} = \det \begin{bmatrix} a_1 & a_2 \\ -ta_1 & -ta_2 \end{bmatrix} = 0.$$

Obviously,

1. $\forall A, B \in \mathbb{R}^2, \angle(A, B) = \angle(-A, -B)$;
2. $\forall t_1, t_2 > 0, \forall A, B \in \mathbb{R}^2, \angle(A, B) = \angle(t_1 A, t_2 B)$.

**Example 1.1.** a) For $A = (0, 1) \in \mathbb{R}^2$, $B = (1, 0) \in \mathbb{R}^2$, we have:

$$\angle(A,B) = \frac{\pi}{2} > 0, \det\begin{bmatrix} a_1 & a_2 \\ b_1 & b_2 \end{bmatrix} = \begin{vmatrix} 0 & 1 \\ 1 & 0 \end{vmatrix} = -1 < 0,$$

$$\angle(B,A) = -\frac{\pi}{2} < 0, \det\begin{bmatrix} b_1 & b_2 \\ a_1 & a_2 \end{bmatrix} = \begin{vmatrix} 1 & 0 \\ 0 & 1 \end{vmatrix} = 1 > 0.$$

b) For $A = (-1,3) \in \mathbb{R}^2, B = (-5,2) \in \mathbb{R}^2$, we have:

$$\det\begin{bmatrix} a_1 & a_2 \\ b_1 & b_2 \end{bmatrix} = \begin{vmatrix} -1 & 3 \\ -5 & 2 \end{vmatrix} = -2 + 15 = 13 > 0 \Rightarrow \angle(A,B) < 0.$$

For points $A, B \in \mathbb{R}^2$, denote by $[A,B] = \{(x,y) = \alpha A + (1-\alpha)B \in \mathbb{R}^2 \mid \alpha \in [0,1]\}$ the line segment joining points $A$ and $B$. Accordingly, $(A,B) = \{(x,y) = \alpha A + (1-\alpha)B \in \mathbb{R}^2 \mid \alpha \in (0,1)\}$.

For arbitrary points $A, B, C \in \mathbb{R}^2$ with $A \neq B$, one can analytically determine where point $C$ lies relative to the oriented line passing through $A$ and $B$ and directed from $A$ to $B$. We briefly denote this line by $\overrightarrow{A,B}$. It is sometimes convenient to consider this same line as unoriented, namely as the set of points belonging to it, i.e., $\{\overrightarrow{A,B}\} = \{A + t(B-A) \mid t \in \mathbb{R}\}$. Then, if point $C$ does not lie on this line, we find the sign of the angle (by computing the corresponding determinant) $\alpha = \angle(B-A, C-A)$ (as will be shown in Statement 1.5, the second point $A$ can be replaced by any point on this line, for example $B$). If $\alpha > 0, \alpha \neq \pi$, then point $C$ lies **strictly to the right** of this line; if $\alpha < 0$, then **strictly to the left.** If point $C$ lies on the line $\overrightarrow{A,B}$ (i.e., $C \in \{\overrightarrow{A,B}\} = \{A + t(B-A) \mid t \in \mathbb{R}\}$), we write briefly $C \in \overrightarrow{A,B}$. It is easy to see (see Statement 1.2) that

$$C \in \overrightarrow{A,B} \iff \begin{vmatrix} b_1 - a_1 & b_2 - a_2 \\ c_1 - a_1 & c_2 - a_2 \end{vmatrix} = \begin{vmatrix} b_1 - a_1 & c_1 - a_1 \\ b_2 - a_2 & c_2 - a_2 \end{vmatrix} = 0 \iff$$

$$\iff \alpha = \angle(B-A, C-A) \in \{0, \pi\}.$$

We shall say that point $C$ is **not strictly to the right** (**left**) of the line $\overrightarrow{A,B}$ if either $\alpha > 0$ $(\alpha < 0)$ or $C \in \overrightarrow{A,B}$.

**Statement 1.1.** *Let* $A, B, C \in \mathbb{R}^2$, $A \neq B$. *The point* $C$ *lies on the line* $\overrightarrow{A, B}$ *if and only if*

$$\angle(B-A, C-A), \angle(C-A, B-A) \in \{0, \pi\},$$

*or if and only if*

$$\angle(B-A, C-B), \angle(C-B, B-A) \in \{0, \pi\}.$$

***Proof.*** Let point $C$ lie on the line $\overrightarrow{A, B}$. Then for some $t \in \mathbb{R}$ we have $C = A + t(B-A)$, and then

$$\angle(B-A, C-A) = \angle(B-A, t(B-A)) = \begin{cases} 0, & \text{if } t \geq 0, \\ \pi, & \text{if } t < 0, \end{cases}$$

($\angle(C-A, B-A) \in \{0, \pi\}$ is proved similarly). Accordingly,

$$\angle(B-A, C-B) = \angle(B-A, (t-1)(B-A)) = \begin{cases} 0, & \text{if } t \geq 1, \\ \pi, & \text{if } t < 1. \end{cases}$$

Note also that $\angle(B-A, C-A) = 0 \Leftrightarrow \angle(C-A, B-A) = 0$, and $\angle(B-A, C-A) = \pi \Leftrightarrow \angle(C-A, B-A) = \pi$.

Now suppose $\angle(B-A, C-A) \in \{0, \pi\}$ (the case $\angle(B-A, C-B) \in \{0, \pi\}$ is treated similarly). We show that point $C$ lies on the line $\overrightarrow{A, B}$. If $\angle(B-A, C-A) = 0$, then $C - A = t(B-A)$ for some $t \geq 0$, whence $C = A + t(B-A) \in \overrightarrow{A, B}$. If $\angle(B-A, C-A) = \pi$, then $C - A = t(B-A)$ for some $t < 0$, whence $C = A + t(B-A) \in \overrightarrow{A, B}$. □

**Statement 1.2.** *Let* $A, B, C \in \mathbb{R}^2$, $A \neq B$. *Then the point* $C$ *lies on the line* $\overrightarrow{A, B}$ *if and only if*

$$\begin{vmatrix} b_1 - a_1 & b_2 - a_2 \\ c_1 - a_1 & c_2 - a_2 \end{vmatrix} = 0.$$

***Proof.*** Let $C \in \overrightarrow{A, B}$. Then for some $t \in \mathbb{R}$ we have $C = A + t(B-A)$, and then

$$\begin{vmatrix} b_1 - a_1 & b_2 - a_2 \\ c_1 - a_1 & c_2 - a_2 \end{vmatrix} = \begin{vmatrix} b_1 - a_1 & b_2 - a_2 \\ t(b_1 - a_1) & t(b_2 - a_2) \end{vmatrix} = 0.$$

Conversely,

$$\begin{vmatrix} b_1 - a_1 & b_2 - a_2 \\ c_1 - a_1 & c_2 - a_2 \end{vmatrix} = 0 \implies \exists \alpha, \beta \neq 0:$$

$$\alpha(B - A) + \beta(C - A) = (0,0) \in \mathbb{R}^2 \implies$$

$$\implies C = A + (-\alpha/\beta)(B - A) \in \overrightarrow{A, B}. \square$$

Thus, the following holds:

**Statement 1.3.** *Let $A, B, C \in \mathbb{R}^2, A \neq B, \alpha = \angle(B - A, C - A)$. Then*

1. *Point $C$ lies **strictly to the right** of the line $\overrightarrow{A, B}$ iff $\alpha \in (0, \pi)$.*

   *Correspondingly, point $C$ lies **strictly to the left** of the line $\overrightarrow{A, B}$ iff $\alpha \in (-\pi, 0)$.*

2. *Point $C$ lies **non-strictly to the right** of the line $\overrightarrow{A, B}$ iff $\alpha \in [0, \pi]$.*

   *Correspondingly, point $C$ lies **non-strictly to the left** of the line $\overrightarrow{A, B}$ iff $\alpha \in (-\pi, 0] \cup \{\pi\}$.*

**Statement 1.4.** *Let $A, B, C \in \mathbb{R}^2, A \neq B, A(t) = (a_1(t), a_2(t)) = A + t(B - A), t \in \mathbb{R}$. Then*

$$\begin{vmatrix} b_1 - a_1 & b_2 - a_2 \\ c_1 - a_1 & c_2 - a_2 \end{vmatrix} = \begin{vmatrix} b_1 - a_1 & b_2 - a_2 \\ c_1 - a_1(t) & c_2 - a_2(t) \end{vmatrix}, \forall t \in \mathbb{R}.$$

***Proof.*** Let $t \in \mathbb{R}$. Then

$$\begin{vmatrix} b_1 - a_1 & b_2 - a_2 \\ c_1 - a_1(t) & c_2 - a_2(t) \end{vmatrix} =$$

$$= \begin{vmatrix} b_1 - a_1 & b_2 - a_2 \\ (c_1 - a_1) - t(b_1 - a_1) & (c_2 - a_2) - t(b_2 - a_2) \end{vmatrix} =$$

$$= \begin{vmatrix} b_1 - a_1 & b_2 - a_2 \\ c_1 - a_1 & c_2 - a_2 \end{vmatrix}. \square$$

**Statement 1.5.** *Let $A, B, C \in \mathbb{R}^2, A \neq B, A(t) = (a_1(t), a_2(t)) = A + t(B - A), t \in \mathbb{R}$. Then*

*(a) if $\angle(B - A, C - A) \in (0, \pi)$, then*

$$\angle(B-A, C-A(t)) \in (0,\pi), \forall t \in \mathbb{R},$$

*and in particular, for* $t=1$: $A(t)=B$, $\angle(B-A, C-B) \in (0,\pi)$.

*(b) if* $\angle(B-A, C-A) \in (-\pi, 0)$*, then*

$$\angle(B-A, C-A(t)) \in (-\pi,0), \forall t \in \mathbb{R},$$

*and in particular, for* $t=1$: $A(t)=B$, $\angle(B-A, C-B) \in (-\pi,0)$.

*(c) if* $\angle(B-A, C-A) \in \{0,\pi\}$*, then*

$$\angle(B-A, C-A(t)) \in \{0,\pi\}, \forall t \in \mathbb{R}.$$

***Proof.*** We prove the first case (the others are analogous). By assumption, $\angle(B-A, C-A) > 0$ and $\angle(B-A, C-A) \neq \pi$, whence

$$\begin{vmatrix} b_1-a_1 & b_2-a_2 \\ c_1-a_1 & c_2-a_2 \end{vmatrix} < 0.$$

For $A(t)=(a_1(t), a_2(t)) = A + t(B-A)$, $t \in \mathbb{R}$, by Statement 1.4 we have

$$\begin{vmatrix} b_1-a_1 & b_2-a_2 \\ c_1-a_1(t) & c_2-a_2(t) \end{vmatrix} = \begin{vmatrix} b_1-a_1 & b_2-a_2 \\ c_1-a_1 & c_2-a_2 \end{vmatrix} < 0,$$

i.e., $\angle(B-A, C-A(t)) \in (0,\pi)$ for all $t \in \mathbb{R}$. In particular, for $t=1$: $A(t) = B$, $\angle(B-A, C-B) \in (0,\pi)$. □

A corollary of Statement 1.5 is:

**Statement 1.6.** *Let* $A, B, C \in \mathbb{R}^2$, $A \neq B$*. Then*

$$\angle(B-A, C-A) \in (0,\pi) \iff \angle(B-A, C-B) \in (0,\pi),$$

$$\angle(B-A, C-A) \in (-\pi,0) \iff \angle(B-A, C-B) \in (-\pi,0),$$

$$\angle(B-A, C-A) \in \{0,\pi\} \iff \angle(B-A, C-B) \in \{0,\pi\}.$$

If during arithmetic operations it turns out that for some vectors $A, B \in \mathbb{R}^2$ the value $\alpha = \angle(A,B) \notin (-\pi,\pi]$, then we determine its value modulo $2\pi$, i.e., we set $\angle(A,B) = \lfloor\alpha\rfloor_{2\pi}$, where

$$\lfloor\alpha\rfloor_{2\pi} = \alpha + 2k\pi \in (-\pi,\pi], k \in \mathbb{Z}. \qquad (1.1)$$

**Statement 1.7.** *The value of* $k \in \mathbb{Z}$ *in condition* (1.1) *exists and is uniquely determined for every* $\alpha \in \mathbb{R}$.

*Proof.* **Existence:** obviously, for every $\alpha \in \mathbb{R}$ there exists $k \in \mathbb{Z}$ such that $2k\pi - \pi < \alpha \leq 2k\pi + \pi$ (since $\mathbb{R} = \bigcup_{k\in\mathbb{Z}}( - \pi + 2k\pi, \pi + 2k\pi])$.

**Uniqueness:** indeed, suppose that

$$\alpha_1 = \lfloor \alpha \rfloor_{2\pi} = \alpha + 2k\pi \in (-\pi, \pi], \alpha_2 = \alpha + 2j\pi \in (-\pi, \pi], k \neq j.$$

Then $k \neq j \Rightarrow \mid \alpha_1 - \alpha_2 \mid \geq 2\pi$, which contradicts the fact that $\alpha_1, \alpha_2 \in (-\pi, \pi] \Rightarrow \mid \alpha_1 - \alpha_2 \mid < 2\pi$. □

Note that for every $\alpha \in \mathbb{R}$ we have

$$\forall j \in \mathbb{Z}, \lfloor \alpha \rfloor_{2\pi} = \lfloor \alpha + 2j\pi \rfloor_{2\pi} \qquad (1.2)$$

(if $\lfloor \alpha \rfloor_{2\pi} = \alpha + 2k\pi \in (-\pi, \pi], k \in \mathbb{Z}$, then $\alpha + 2j\pi + 2(k - j)\pi = \alpha + 2k\pi \in (-\pi, \pi]$, i.e., $\lfloor \alpha + 2j\pi \rfloor_{2\pi} = \alpha + 2k\pi = \lfloor \alpha \rfloor_{2\pi}$).

For any $\alpha, \beta \in \mathbb{R}$ denote

$$\alpha \oplus \beta = \lfloor \alpha + \beta \rfloor_{2\pi}.$$

We show that

$$\forall \alpha, \beta \in \mathbb{R} \quad \lfloor \alpha + \beta \rfloor_{2\pi} = \lfloor \alpha \rfloor_{2\pi} \oplus \lfloor \beta \rfloor_{2\pi}. \qquad (1.3)$$

Indeed, let $\lfloor \alpha \rfloor_{2\pi} = \alpha + 2k\pi$, $\lfloor \beta \rfloor_{2\pi} = \beta + 2j\pi$, $k, j \in \mathbb{Z}$. Then by (1.2)

$$\lfloor \alpha \rfloor_{2\pi} \oplus \lfloor \beta \rfloor_{2\pi} = \lfloor \lfloor \alpha \rfloor_{2\pi} + \lfloor \beta \rfloor_{2\pi} \rfloor_{2\pi} = \lfloor \alpha + \beta + 2(k + j)\pi \rfloor_{2\pi} = \lfloor \alpha + \beta \rfloor_{2\pi},$$

so (1.3) is proved.

From (1.3), using mathematical induction, one easily proves:

**Statement 1.8.**

$$\forall \alpha_1, \dots, \alpha_r \in \mathbb{R} \quad \lfloor \alpha_1 + \cdots + \alpha_r \rfloor_{2\pi} = \lfloor \alpha_1 \rfloor_{2\pi} \oplus \cdots \oplus \lfloor \alpha_r \rfloor_{2\pi}. \qquad (1.4)$$

Moreover, using (1.1), (1.2), (1.3), we obtain

$$\forall \alpha, \beta \in \mathbb{R}, \lfloor \alpha + \lfloor \beta \rfloor_{2\pi} \rfloor_{2\pi} = \lfloor \alpha + \beta \rfloor_{2\pi} = \lfloor \alpha \rfloor_{2\pi} \oplus \lfloor \beta \rfloor_{2\pi}. \qquad (1.5)$$

Using (1.5), we show the associativity of $\oplus$ for angles $\alpha, \beta, \gamma \in \mathbb{R}$:

$$\alpha \oplus (\beta \oplus \gamma) = \alpha \oplus \lfloor \beta + \gamma \rfloor_{2\pi} = \lfloor \alpha + \lfloor \beta + \gamma \rfloor_{2\pi} \rfloor_{2\pi} = \lfloor \alpha + \beta + \gamma \rfloor_{2\pi},$$

$$(\alpha \oplus \beta) \oplus \gamma = \lfloor \alpha + \beta \rfloor_{2\pi} \oplus \gamma = \lfloor \lfloor \alpha + \beta \rfloor_{2\pi} + \gamma \rfloor_{2\pi} = \lfloor \alpha + \beta + \gamma \rfloor_{2\pi},$$

i.e.,

$$\alpha \oplus (\beta \oplus \gamma) = \lfloor \alpha + \beta + \gamma \rfloor_{2\pi} = (\alpha \oplus \beta) \oplus \gamma.$$

Commutativity of $\oplus$ is also obviously true.

**Example 1.2.** a) $3\pi/2 = \pi + \pi/2 \notin (-\pi, \pi] \Rightarrow \lfloor 3\pi/2 \rfloor_{2\pi} = 3\pi/2 - 2\pi = -\pi/2 \in (-\pi, \pi]$. b) $-5\pi \notin (-\pi, \pi] \Rightarrow \lfloor -5\pi \rfloor_{2\pi} = -5\pi + 3 \cdot 2\pi = \pi \in (-\pi, \pi]$. Note that $-5\pi + 2 \cdot 2\pi = -\pi \notin (-\pi, \pi]$, i.e., $\lfloor -5\pi \rfloor_{2\pi} \neq -\pi$.

Let $A, B \in \mathbb{R}^2$, $A, B \neq (0,0)$. Denote by $\overrightarrow{\angle(A,B)}$ the full clockwise angle from $A$ to $B$. Then $\overrightarrow{\angle(A,B)} \in [0,2\pi)$ ( $\forall t > 0, \overrightarrow{\angle(A,tA)} = 0, \overrightarrow{\angle(A,-tA)} = \pi$ ), and $\angle(A,B) = \lfloor \overrightarrow{\angle(A,B)} \rfloor_{2\pi}$. Geometrically, the value of the angle $\overrightarrow{\angle(A,B)}$ equals the length of the arc of the unit circle when traversing that circle from the point $\frac{1}{|A|}A$ to the point $\frac{1}{|B|}B$ clockwise.

Then, if $\forall t > 0 \;\; A \neq tB$ (i.e., $\frac{1}{|A|}A \neq \frac{1}{|B|}B$), using Statement 1.8 we obtain

$$\overrightarrow{\angle(A,B)} + \overrightarrow{\angle(B,A)} = 2\pi, \angle(A,B) = \lfloor \overrightarrow{\angle(A,B)} \rfloor_{2\pi} =$$
$$\lfloor 2\pi - \overrightarrow{\angle(B,A)} \rfloor_{2\pi} = \lfloor -\overrightarrow{\angle(B,A)} \rfloor_{2\pi}. \qquad (1.6)$$

We say that the angle $\overrightarrow{\angle(A,B)}$ **passes through** a vector $C \neq (0,0)$ (write: $C \in \overrightarrow{\angle(A,B)}$) if $\overrightarrow{\angle(A,C)} \leq \overrightarrow{\angle(A,B)}$. In this case (because when traversing the unit circle clockwise from $\frac{1}{|A|}A$ to $\frac{1}{|B|}B$ we pass the intermediate point $\frac{1}{|C|}C$, so that the length of the traversed arc from $\frac{1}{|A|}A$ to $\frac{1}{|B|}B$ equals the sum of the arc lengths from $\frac{1}{|A|}A$ to $\frac{1}{|C|}C$ and from $\frac{1}{|C|}C$ to $\frac{1}{|B|}B$) we have

$$\overrightarrow{\angle(A,B)} = \overrightarrow{\angle(A,C)} + \overrightarrow{\angle(C,B)} \qquad (1.7)$$

(in particular, it may be that $\frac{1}{|A|}A = \frac{1}{|C|}C$ or $\frac{1}{|B|}B = \frac{1}{|C|}C$), whence by (1.3)

$$\angle(A,B) = \lfloor \overrightarrow{\angle(A,B)} \rfloor_{2\pi} = \lfloor \overrightarrow{\angle(A,C)} + \overrightarrow{\angle(C,B)} \rfloor_{2\pi} =$$
$$\lfloor \overrightarrow{\angle(A,C)} \rfloor_{2\pi} \oplus \lfloor \overrightarrow{\angle(C,B)} \rfloor_{2\pi} = \angle(A,C) \oplus \angle(C,B). \qquad (1.8)$$

Now suppose the angle $\overrightarrow{\angle(A,B)}$ **does not pass through** the vector $C$ (write: $C \notin \overrightarrow{\angle(A,B)}$ ). Then $\overrightarrow{\angle(A,C)} > \overrightarrow{\angle(A,B)}$, i.e., analogously to (1.7)

$$\overrightarrow{\angle(A,C)} = \overrightarrow{\angle(A,B)} + \overrightarrow{\angle(B,C)},$$

hence

$$\overrightarrow{\angle(A,B)} = \overrightarrow{\angle(A,C)} - \overrightarrow{\angle(B,C)},$$

and consequently, using (1.3) and (1.6), we obtain

$$\angle(A,B) = \lfloor\overrightarrow{\angle(A,B)}\rfloor_{2\pi} = \lfloor\overrightarrow{\angle(A,C)} - \overrightarrow{\angle(B,C)}\rfloor_{2\pi} =$$

$$\lfloor\overrightarrow{\angle(A,C)}\rfloor_{2\pi} \oplus \lfloor-\overrightarrow{\angle(B,C)}\rfloor_{2\pi} =$$

$$\lfloor\overrightarrow{\angle(A,C)}\rfloor_{2\pi} \oplus \lfloor\overrightarrow{\angle(C,B)}\rfloor_{2\pi} = \angle(A,C) \oplus \angle(C,B),$$

i.e., (1.8) holds again.

**Remark 1.2.** Equality (1.8) also holds in the case $A = tB$, $B \neq (0,0)$, $t > 0$ (with $\angle(A,B) = 0$) for an arbitrary $C \neq (0,0)$. Indeed, if $C = \alpha B$ for some $\alpha < 0$, then $\angle(A,C) = \angle(tB,\alpha B) = \angle(B,-B) = \pi$, $\angle(C,B) = \angle(\alpha B,tB) = \angle(-B,B) = \pi$, and then

$$0 = \angle(A,B) = \angle(A,C) \oplus \angle(C,B) = \pi \oplus \pi = \lfloor 2\pi\rfloor_{2\pi} = 0.$$

If $C \neq \alpha B$ for all $\alpha < 0$, then $\angle(A,C), \angle(C,B) \neq \pi$, $\angle(A,C) = \angle(tB,C) = \angle(B,C) = -\angle(C,B)$, and then

$$0 = \angle(A,B) = \angle(A,C) \oplus \angle(C,B) = \lfloor-\angle(C,B) + \angle(C,B)\rfloor_{2\pi} = \lfloor 0\rfloor_{2\pi} = 0.$$

Using mathematical induction, it is now easy to prove

**Statement 1.9.** *Let* $A, B, C^{(1)}, \dots, C^{(r)} \in \mathbb{R}^2$*,* $A, B, C^{(1)}, \dots, C^{(r)} \neq (0,0)$*. Then*

$$\angle(A,B) = \angle(A,C^{(1)}) \oplus \angle(C^{(1)},C^{(2)}) \oplus \cdots \oplus \angle(C^{(r-1)},C^{(r)}) \oplus \angle(C^{(r)},B). \quad (1.9)$$

Using (1.4) (see Statement 1.8), it is easy to prove

**Statement 1.10.** *Let the conditions of Statement 1.9 hold. Then, if*

$$\angle(A,C^{(1)}) + \angle(C^{(1)},C^{(2)}) + \cdots + \angle(C^{(r-1)},C^{(r)}) + \angle(C^{(r)},B) \in (-\pi,\pi],$$

*we have*

$$\angle(A,B) = \angle(A,C^{(1)}) \oplus \angle(C^{(1)},C^{(2)}) \oplus \cdots \oplus \angle(C^{(r-1)},C^{(r)}) \oplus \angle(C^{(r)},B)$$

$$= \lfloor \angle(A,C^{(1)}) + \angle(C^{(1)},C^{(2)}) + \cdots + \angle(C^{(r-1)},C^{(r)}) + \angle(C^{(r)},B) \rfloor_{2\pi}$$

$$= \angle(A,C^{(1)}) + \angle(C^{(1)},C^{(2)}) + \cdots + \angle(C^{(r-1)},C^{(r)}) + \angle(C^{(r)},B).$$

We now present some additional statements concerning angles.

**Statement 1.11.** *Let* $\alpha \in \mathbb{R}$. *Then*

1. $\pi \oplus 0 = \lfloor \pi + 0 \rfloor_{2\pi} = \lfloor \pi \rfloor_{2\pi} = \pi,$
2. $\pi \oplus \pi = \lfloor \pi + \pi \rfloor_{2\pi} = \lfloor 2\pi \rfloor_{2\pi} = 2\pi - 2\pi = 0,$
3. $\alpha \in (0,\pi) \Rightarrow \pi \oplus \alpha = \lfloor \pi + \alpha \rfloor_{2\pi} = \pi + \alpha - 2\pi = \alpha - \pi \in (-\pi, 0),$
4. $\alpha \in (-\pi,0) \Rightarrow \pi \oplus \alpha = \lfloor \pi + \alpha \rfloor_{2\pi} = \pi + \alpha \in (-\pi, 0).$

**Statement 1.12.** *Let* $A, B \in \mathbb{R}^2$, $A, B \neq (0,0)$. *Then*

1. $\angle(A,-B) = \angle(-A,B) = \pi \oplus \angle(A,B),$
2. $\angle(A,B) \in (0,\pi) \iff \angle(A,-B) \in (-\pi,0) \iff \angle(-A,B) \in (-\pi,0),$
3. $\angle(A,B) \in (-\pi,0) \iff \angle(A,-B) \in (0,\pi) \iff \angle(-A,B) \in (0,\pi),$
4. $\angle(A,B) = 0 \iff \angle(A,-B) = \pi \iff \angle(-A,B) = \pi,$
5. $\angle(A,B) = \pi \iff \angle(A,-B) = 0 \iff \angle(-A,B) = 0.$

***Proof.*** We prove part 1) (parts 2)–5) follow from Statement 1.11 and part 1)). Using Statement 1.9, we have

$$\angle(B,A) \oplus \angle(A,-B) = \angle(B,-B) = \pi,$$

hence (using commutativity of $\oplus$)

$$\pi \oplus \angle(A,B) = \angle(A,B) \oplus \pi = \angle(A,B) \oplus [\angle(B,A) \oplus \angle(A,-B)],$$

and therefore (using associativity of $\oplus$ and Statement 1.9)

$$\pi \oplus \angle(A,B) = [\angle(A,B) \oplus \angle(B,A)] \oplus \angle(A,-B) =$$

$$\angle(A,A) \oplus \angle(A,-B) = \angle(A,-B). \square$$

**Statement 1.13.** *Let* $\alpha \in (0,\pi)$, $\beta \in (-\pi,\pi]$, *and suppose* $\alpha \oplus \beta = \pi$. *Then* $\beta = \pi - \alpha \in (0,\pi)$.

***Proof.*** From $\alpha \oplus \beta = \pi$ we have

$$\alpha + \beta \in \{\pi + 2j\pi \mid j \in \mathbb{Z}\}. \qquad (1.10)$$

Since $\alpha \in (0, \pi), \beta \in (-\pi, \pi]$, we obtain

$$\alpha + \beta \in (-\pi, 2\pi). \qquad (1.11)$$

Condition (1.11) in the case of $\alpha + \beta = \pi + 2j\pi, j \in \mathbb{Z}$ (see condition (1.10)), is satisfied only for $j = 0$, which corresponds to the equality $\beta = \pi - \alpha$ (for $j = 1$ we get $\alpha + \beta = \pi + 2j\pi = 3\pi \notin (-\pi, 2\pi)$; for $j = -1$ we get $\alpha + \beta = \pi + 2j\pi = -\pi \notin (-\pi, 2\pi)$; for all other values of $j$ condition (1.11) is not satisfied either). □

**Statement 1.14.** *Let* $\alpha \in (-\pi, 0)$*,* $\beta \in (-\pi, \pi]$*, and suppose* $\alpha \oplus \beta = \pi$*. Then* $\beta = -\pi - \alpha \in (-\pi, 0)$*.*

***Proof.*** From $\alpha \oplus \beta = \pi$ we have

$$\alpha + \beta \in \{\pi + 2j\pi \mid j \in \mathbb{Z}\}. \qquad (1.12)$$

Since $\alpha \in (-\pi, 0), \beta \in (-\pi, \pi]$, we obtain

$$\alpha + \beta \in (-2\pi, \pi). \qquad (1.13)$$

Condition (1.12) together with (13) forces $j = -1$ only, which gives $\beta = -\pi - \alpha$ (for $j = 0$ we get $\alpha + \beta = \pi \notin (-2\pi, \pi)$; for $j = -2$ we get $\alpha + \beta = -3\pi \notin (-2\pi, \pi)$; other $j$ are excluded). □

**Statement 1.15.** *Let either*

*1.* $\delta \in [0, \pi]$*,* $\gamma \in (-\pi, \pi]$*,* $\alpha \in (0, \pi]$*,*

*or*

*2.* $\delta \in [0, \pi)$*,* $\gamma \in (-\pi, \pi]$*,* $\alpha \in [0, \pi)$*,*

*and suppose* $\delta \oplus \gamma = \alpha$*. Then* $\delta + \gamma = \alpha$*.*

***Proof.*** Consider case 1) (case 2) is analogous). From $\delta \oplus \gamma = \alpha$ we have

$$\delta + \gamma \in \{\alpha + 2j\pi \mid j \in \mathbb{Z}\}. \qquad (1.14)$$

Since $\delta \in [0, \pi]$ (resp. $[0, \pi)$), $\gamma \in (-\pi, \pi]$, we obtain

$$\delta + \gamma \in (-\pi, 2\pi] \text{ (resp. } (-\pi, 2\pi)). \qquad (1.15)$$

Condition (1.14) together with (1.15) forces $j = 0$ only, giving $\delta + \gamma = \alpha$ (for $j = 1$ we have $\delta + \gamma = \alpha + 2\pi \notin (-\pi, 2\pi]$ (resp. $(-\pi, 2\pi)$) because $\alpha + 2\pi \in (2\pi, 3\pi]$ (resp. $[2\pi, 3\pi)$); for $j = -1$ we have $\delta + \gamma = \alpha - 2\pi \notin (-\pi, 2\pi]$ (resp. $(-\pi, 2\pi)$) because $\alpha - 2\pi \in (-2\pi, -\pi]$ (resp. $[-2\pi, -\pi)$); other $j$ are excluded). □

**Statement 1.16.** *Let* $\delta \in (-\pi, 0]$, $\gamma \in (-\pi, \pi]$, $\alpha \in (-\pi, 0]$, *and suppose* $\delta \oplus \gamma = \alpha$. *Then* $\delta + \gamma = \alpha$.

***Proof.*** From $\delta \oplus \gamma = \alpha$ we have

$$\delta + \gamma \in \{\alpha + 2j\pi \mid j \in \mathbb{Z}\}. \qquad (1.16)$$

Since $\delta \in (-\pi, 0]$, $\gamma \in (-\pi, \pi]$, we obtain

$$\delta + \gamma \in (-2\pi, \pi]. \qquad (1.17)$$

Condition (1.16) together with (1.17) forces $j = 0$ only, giving $\delta + \gamma = \alpha$ (for $j = 1$ we get $\delta + \gamma = \alpha + 2\pi \notin (-2\pi, \pi]$ because $\alpha + 2\pi \in (\pi, 2\pi]$; for $j = -1$ we get $\delta + \gamma = \alpha - 2\pi \notin (-2\pi, \pi]$ because $\alpha - 2\pi \in (-3\pi, -2\pi]$; other $j$ excluded). □

**Statement 1.17.** *Let (see fig. 1.1)* $A, B, C, D, E \in \mathbb{R}^2$, $A \neq B$, $D \neq A$, *points* $C, E$ *lie strictly to the right of the line* $\overrightarrow{A, B}$, *i.e.* $\angle(B - A, C - A), \angle(B - A, E - A) \in (0, \pi)$ *(see Statement 1.3), point* $D$ *lies not strictly to the right of the line* $\overrightarrow{A, B}$, *i.e.* $\angle(B - A, D - A) \in [0, \pi]$ *(see Statement 1.3). Let*

$$\alpha = \angle(E - A, C - A) \in (-\pi, 0], \beta = \angle(D - A, C - A) \geq \alpha,$$
$$\gamma = \angle(D - A, E - A).$$

*Then*

1. $\gamma = \beta - \alpha \in [0, \pi)$, $\angle(E - A, D - A) = -\gamma \in (-\pi, 0]$, *and consequently point* $D$ *lies not strictly to the left of the line* $\overrightarrow{A, E}$ *(see Statement 1.3).*

2. *If* $\beta > \alpha$*, then* $\gamma = \beta - \alpha \in (0, \pi)$*, hence* $\angle(E - A, D - A) = -\gamma \in (-\pi, 0)$*, i.e. point* $D$ *lies strictly to the left of the line* $\overrightarrow{A, E}$ *(see Statement 1.3), and by Statement 1.6 we have* $\angle(E - A, D - E) \in (-\pi, 0)$*,* $\angle(D - E, E - A) \in (0, \pi)$*.*
3. *If* $\gamma = 0$*, then* $\angle(D - E, E - A) \in \{0, \pi\}$*.*

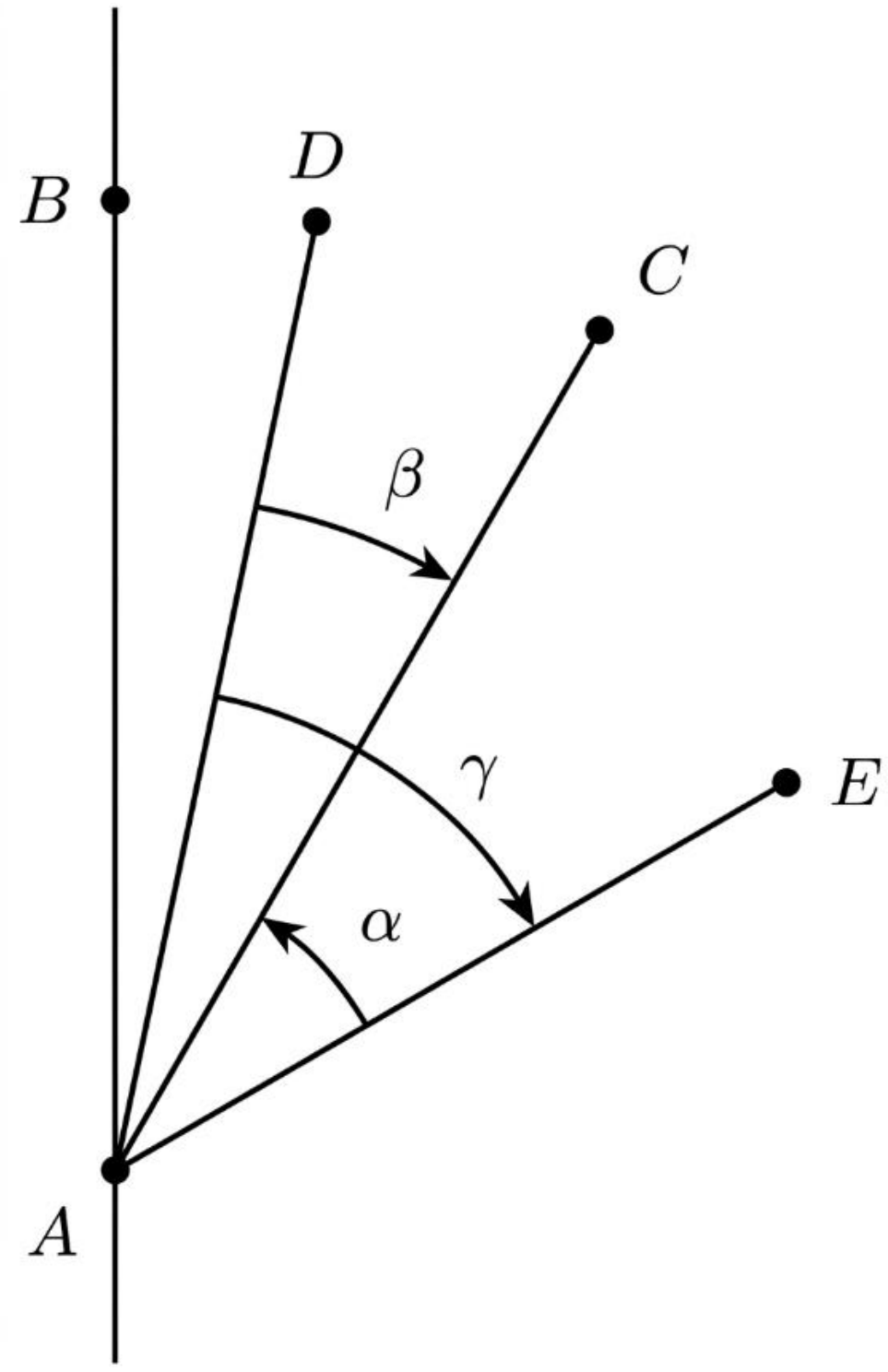

**Fig. 1.1.** The layout of points and angles following statement 1.17

***Proof.*** We prove 1). By Statement 1.9 (since $B \neq A, D \neq A, E \neq A$):

$$\angle(B - A, E - A) = \angle(B - A, D - A) \oplus \angle(D - A, E - A),$$

and by Statement 1.15 (using $\angle(B - A, E - A) \in (0, \pi), \angle(B - A, D - A) \in [0, \pi]$) we obtain

$$\angle(B - A, E - A) = \angle(B - A, D - A) + \angle(D - A, E - A),$$

and therefore

$$\angle(B - A, E - A) = \angle(B - A, D - A) + \gamma. \qquad (1.18)$$

Further, by Statement 1.9 (using $B \neq A, D \neq A, C \neq A$):

$$\angle(B-A,C-A)=\angle(B-A,D-A)\oplus\angle(D-A,C-A),$$

and by Statement 1.15 (using $\angle(B-A,C-A)\in(0,\pi)$, $\angle(B-A,D-A)\in[0,\pi]$) we get

$$\angle(B-A,C-A)=\angle(B-A,D-A)+\angle(D-A,C-A),$$

and therefore

$$\angle(B-A,C-A)=\angle(B-A,D-A)+\beta. \qquad (1.19)$$

Again, by Statement 1.9:

$$\angle(B-A,E-A)=\angle(B-A,C-A)\oplus\angle(C-A,E-A),$$

and by Statement 1.15 (since $\angle(B-A,E-A),\angle(B-A,C-A)\in(0,\pi)$) we have

$$\angle(B-A,E-A)=\angle(B-A,C-A)+\angle(C-A,E-A)=\angle(B-A,C-A)-\alpha,$$

or, using (1.19),

$$\angle(B-A,E-A)=\angle(B-A,D-A)+\beta-\alpha. \qquad (1.20)$$

From (1.18) and (1.20), using $\beta\geq\alpha$, we get $\gamma=\angle(B-A,E-A)-\angle(B-A,D-A)=\beta-\alpha\geq 0$, hence $\gamma\in[0,\pi)$ ($\gamma\neq\pi$; if $\gamma=\pi$ then $E-A=t(A-D)$ with $t>0$, and then $\angle(B-A,A-D)=\angle(B-A,E-A)\in(0,\pi)$, implying by Statement 1.12(2) that $\angle(B-A,D-A)\in(-\pi,0)$, contradicting the assumption $\angle(B-A,D-A)\in[0,\pi]$).

2) follows from the equality $\gamma=\beta-\alpha$.

3) If $\gamma=0$ then $D-A=t(E-A)$ with $t>0$ (since $D\neq A$), so $D-E=(t-1)(E-A)$ and $\angle(D-E,E-A)=\angle((t-1)(E-A),E-A)\in\{0,\pi\}$. □

**Statement 1.18.** *Let (see fig. 1.2)* $A,B,C,D\in\mathbb{R}^2$, $A\neq B$, $B\neq C$. *Point* $D$ *lies not strictly to the right of the line* $\overrightarrow{B,C}$ *and not strictly to the left of the line* $\overrightarrow{A,B}$. *Let* $\alpha=\angle(C-B,B-A)\in(0,\pi)$. *Then*

$$\beta=\angle(C-B,D-B)\in[0,\alpha],\gamma=\angle(D-B,B-A)\in[0,\alpha]. \qquad (1.21)$$

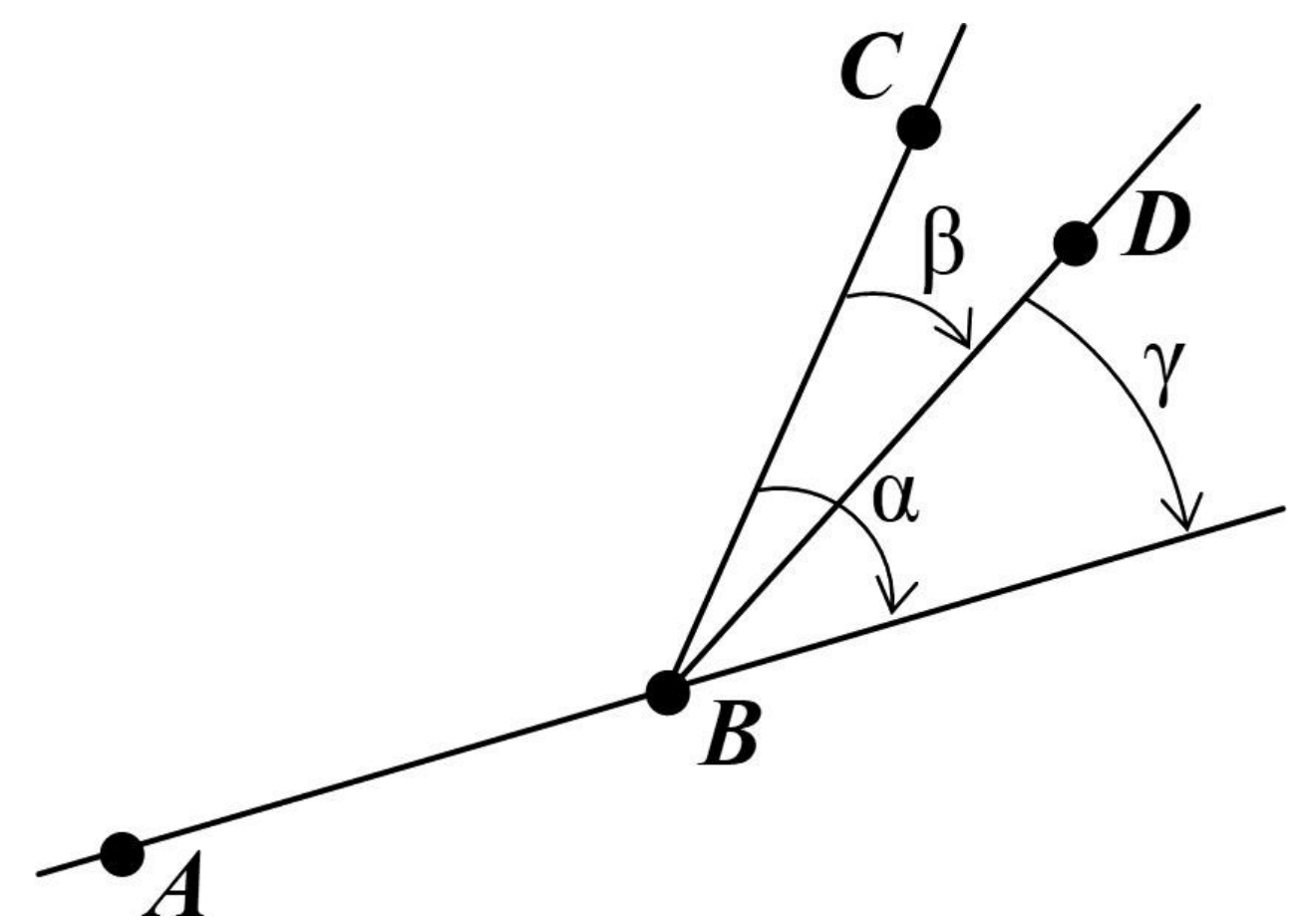

**Fig. 1.2.** The layout of points and angles following statement 1.18

***Proof.*** Consider the nontrivial case $D \neq B$. By Statement 1.9, $\alpha = \beta \oplus \gamma$. Since $D$ is not strictly to the right of $\overrightarrow{B,C}$, we have $\beta = \angle(C - B, D - B) \in [0, \pi]$ (Statement 1.3(2)). Since $D$ is not strictly to the left of $\overrightarrow{A,B}$, we have $\angle(B - A, D - B) \in (-\pi, 0] \cup \{\pi\}$ (Statement 1.3(2)), hence $\gamma = \angle(D - B, B - A) \in [0, \pi]$. By Statement 1.15 we obtain $\alpha = \beta + \gamma$, and because $\beta, \gamma \geq 0$ we get (1.21). □

**Statement 1.19.** *Let* $A, B, C, D \in \mathbb{R}^2$*,* $A \neq B$*. Points* $C, D$ *lie strictly to the right of the line* $\overrightarrow{A,B}$*, i.e.* $\alpha = \angle(B - A, D - A) \in (0, \pi)$*,* $\delta = \angle(B - A, C - A) \in (0, \pi)$ *(see Statement 1.3(1)). Suppose additionally* $\delta \leq \alpha$*. Then* $\gamma = \angle(C - A, D - A) = \alpha - \delta \in [0, \pi)$*,* $\angle(D - A, C - A) = -\gamma \in (-\pi, 0]$*, i.e. point* $C$ *lies not strictly to the left of the line* $\overrightarrow{A,D}$ *(see Statement 1.3(2)).*

***Proof.*** From the conditions, $C \neq A$. By Statement 1.9:

$$\delta \oplus \gamma = \angle(B - A, C - A) \oplus \angle(C - A, D - A) = \angle(B - A, D - A) = \alpha.$$

By Statement 1.8, $\delta + \gamma = \alpha$. Since $\delta \leq \alpha$ and $\alpha \in (0, \pi)$, we get $\gamma = \alpha - \delta \in [0, \pi)$. □

**Statement 1.20.** *Let* $A, B, C \in \mathbb{R}^2$*,* $A \neq B$*, and suppose point* $C$ *lies strictly to the left (resp. right) of the line* $\overrightarrow{A,B}$*. Then point* $B$ *lies strictly to the right*

*(resp. left) of the line* $\overrightarrow{A,C}$*, and point* $A$ *lies strictly to the left (resp. right) of the line* $\overrightarrow{B,C}$*.*

***Proof.*** Assume $C$ is strictly left of $\overrightarrow{A,B}$ (the other case is analogous). By Statement 1.3(1), $\angle(B-A,C-A) \in (-\pi,0)$, so $\angle(C-A,B-A) \in (0,\pi)$, hence by Statement 1.3(1), $B$ lies strictly to the right of $\overrightarrow{A,C}$. Next, by Statement 6, $\angle(C-A,B-C) \in (0,\pi)$, hence $\angle(B-C,C-A) = \angle(C-B,A-C) \in (-\pi,0)$; therefore by Statement 1.3(1), $A$ lies strictly to the left of $\overrightarrow{B,C}$. □

**Statement 1.21***. Let (see fig. 1.3)* $A,B,C$ *be pairwise distinct points in* $\mathbb{R}^2$*,* $\alpha = \angle(C-B,B-A) \in (0,\pi)$*. Then*

$$\gamma = \angle(C-A,B-A) \in (0,\alpha), \delta = \angle(C-B,C-A) \in (0,\alpha).$$

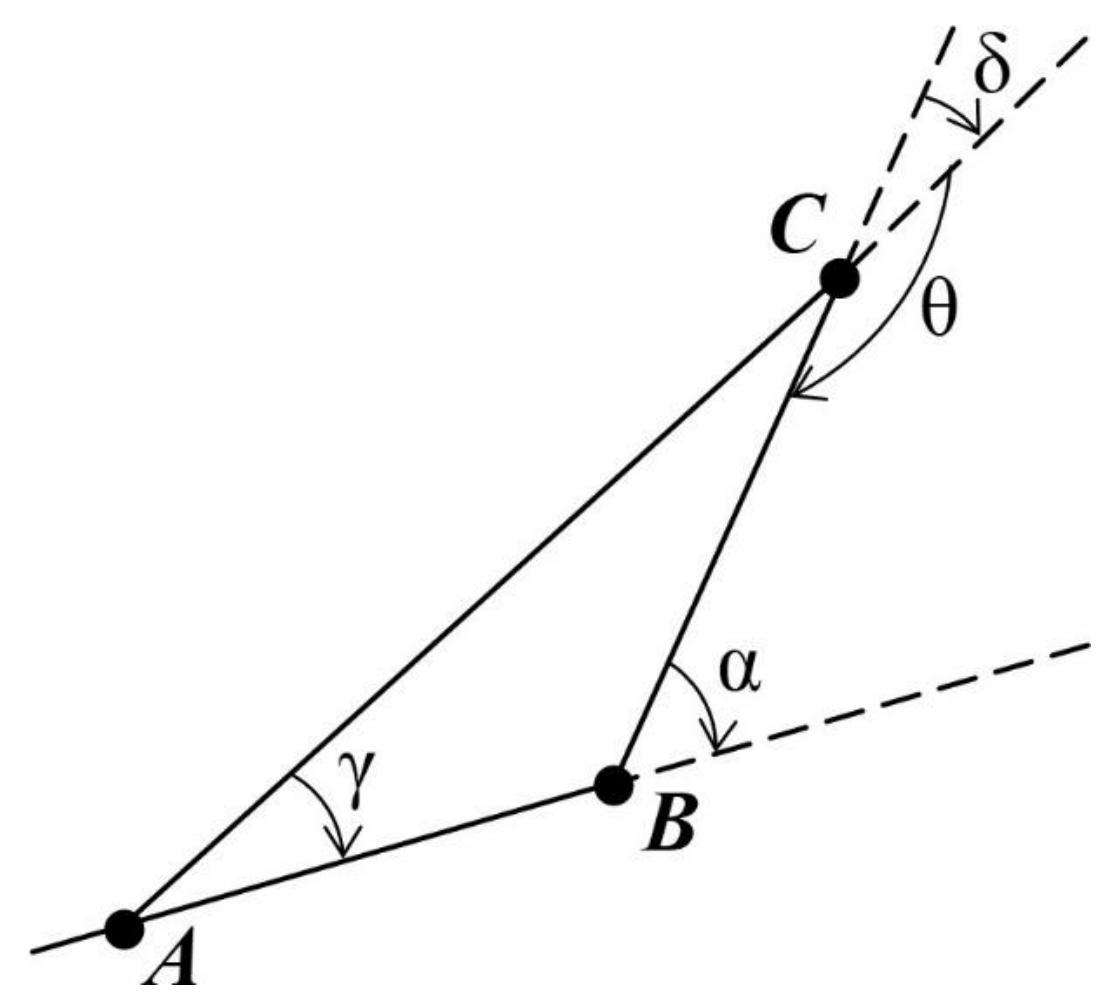

**Fig. 3.** The layout of points and angles following statement 1.21

***Proof.*** From $\alpha \in (0,\pi)$ we have $\angle(B-A,C-B) = -\alpha \in (-\pi,0)$. By Statement 1.6, $\angle(B-A,C-A) \in (-\pi,0)$, so $\gamma = \angle(C-A,B-A) \in (0,\pi)$ and again by Statement 1.6, $\theta = \angle(C-A,B-C) \in (0,\pi)$.

By Statement 1.9:

$$\delta \oplus \theta = \angle(C-B,C-A) \oplus \angle(C-A,B-C) = \angle(C-B,B-C) = \pi,$$

and by Statement 1.13 (using $\theta \in (0,\pi)$) we get $\delta = \pi - \theta \in (0,\pi)$. Moreover,

$$\delta \oplus \gamma = \angle(C-B,C-A) \oplus \angle(C-A,B-A) = \angle(C-B,B-A) = \alpha,$$

and by Statement 1.15 (since $\alpha, \delta, \gamma \in (0, \pi)$) we have $\delta + \gamma = \alpha$, hence $\gamma = \alpha - \delta \in (0, \alpha)$ and $\delta = \alpha - \gamma \in (0, \alpha)$. □

**Remark 1.3.** *Under the conditions of Statement 1.21, let* $\theta = \angle(C - A, B - C)$*. Then* $\gamma = \angle(C - A, B - A) \in (0, \theta)$.

***Proof.*** Set $\beta = \angle(B - A, B - C)$. Since $\alpha = \angle(C - B, B - A) \in (0, \pi)$, by Statement 12(3) we have $\angle(C - B, A - B) \in (-\pi, 0)$, so

$$\beta = \angle(B - A, B - C) = \angle(A - B, C - B) = -\angle(C - B, A - B) \in (0, \pi).$$

Note that

$$\gamma \oplus \beta = \angle(C - A, B - A) \oplus \angle(B - A, B - C) = \angle(C - A, B - C) = \theta,$$

with $\gamma, \beta, \theta \in (0, \pi)$ (see proof of Statement 21). By Statement 15, $\theta = \gamma \oplus \beta = \gamma + \beta$, so $\gamma = \theta - \beta \in (0, \theta)$. □

**Statement 1.22.** *Let (see fig. 1.4)* $A, B, C$ *be pairwise distinct points in* $\mathbb{R}^2$*,* $\alpha = \angle(C - B, B - A) \in (-\pi, 0)$*. Then*
$\gamma = \angle(C - A, B - A) \in (\alpha, 0), \delta = \angle(C - B, C - A) = \angle(B - C, A - C) \in (\alpha, 0)$.

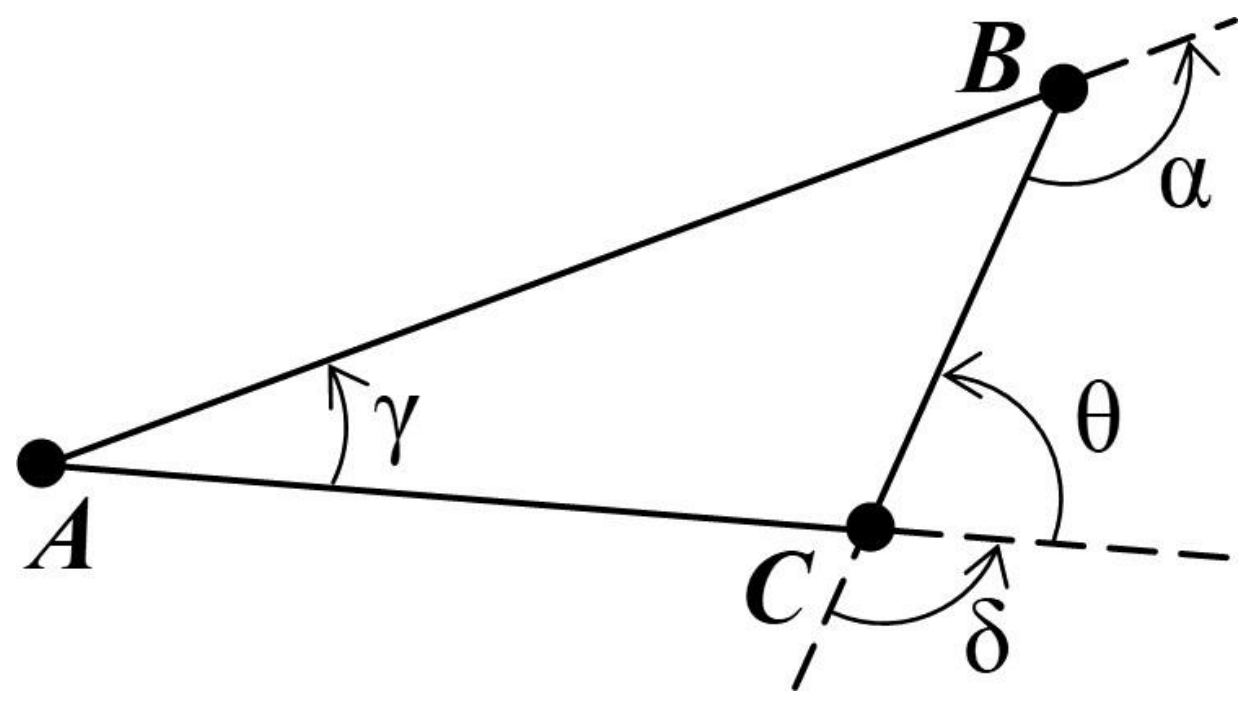

**Figure 1.4.** The layout of points and angles following statement 1.22

***Proof.*** From $\alpha \in (-\pi, 0)$ we get $\angle(B - A, C - B) = -\alpha \in (0, \pi)$. By Statement 1.6, $\angle(B - A, C - A) \in (0, \pi)$, so $\gamma = \angle(C - A, B - A) \in (-\pi, 0)$ and by Statement 1.6, $\theta = \angle(C - A, B - C) \in (-\pi, 0)$. By Statement 1.9:

$$\delta \oplus \theta = \angle(C - B, C - A) \oplus \angle(C - A, B - C) = \angle(C - B, B - C) = \pi,$$

and by Statement 1.14 (since $\theta \in (-\pi, 0)$) we obtain $\delta = -\pi - \theta \in (-\pi, 0)$. Also,

$$\delta \oplus \gamma = \angle(C-B, C-A) \oplus \angle(C-A, B-A) = \angle(C-B, B-A) = \alpha,$$

and by Statement 1.16 (since $\alpha, \delta, \gamma \in (-\pi, 0)$) we have $\delta + \gamma = \alpha$, hence $\gamma = \alpha - \delta \in (\alpha, 0)$ and $\delta = \alpha - \gamma \in (\alpha, 0)$. □

**Statement 1.23.** *Let* $A, B, E, D \in \mathbb{R}^2$, $A \neq B$, $E \neq D$. *Suppose points* $E, D$ *lie strictly to the right of the line* $\overrightarrow{A, B}$ *and*

$$\angle(D-B, E-A) < 0. \qquad (1.22)$$

*Then the lines* $\overrightarrow{A, E}$ *and* $\overrightarrow{B, D}$ *intersect at a unique point* $C \in \mathbb{R}^2$ *that lies strictly to the right of the line* $\overrightarrow{A, B}$ *.*

***Proof.*** The sets of points belonging to the lines $\overrightarrow{A, E}$ and $\overrightarrow{B, D}$, can be described parametrically:

$$\{\overrightarrow{A, E}\} = \{A(t_1) = A + t_1(E - A) \mid t_1 \in \mathbb{R}\},$$

$$\{\overrightarrow{B, D}\} = \{B(t_2) = B + t_2(D - B) \mid t_2 \in \mathbb{R}\}.$$

Intersection condition $A(t_1) = B(t_2)$ leads to the linear system

$$\begin{cases} t_1(e_1 - a_1) + t_2(b_1 - d_1) = b_1 - a_1, \\ t_1(e_2 - a_2) + t_2(b_2 - d_2) = b_2 - a_2. \end{cases} \qquad (1.23)$$

A unique solution exists iff

$$\Delta = \begin{vmatrix} e_1 - a_1 & b_1 - d_1 \\ e_2 - a_2 & b_2 - d_2 \end{vmatrix} \neq 0.$$

From (1.22) we have

$$0 < \begin{vmatrix} d_1 - b_1 & d_2 - b_2 \\ e_1 - a_1 & e_2 - a_2 \end{vmatrix} = \begin{vmatrix} d_1 - b_1 & e_1 - a_1 \\ d_2 - b_2 & e_2 - a_2 \end{vmatrix} = \begin{vmatrix} e_1 - a_1 & b_1 - d_1 \\ e_2 - a_2 & b_2 - d_2 \end{vmatrix},$$

so $\Delta > 0$ and the intersection point $C$ exists and is unique. From (1.23) by Cramer's rule we have

$$C = A + \frac{\Delta_1}{\Delta}(E - A), \Delta_1 = \begin{vmatrix} b_1 - a_1 & b_1 - d_1 \\ b_2 - a_2 & b_2 - d_2 \end{vmatrix}.$$

Since $D$ lies strictly to the right of $\overrightarrow{A, B}$, by Statement 3(1) we have $\angle(B - A, D - A) \in (0, \pi)$ and by Statement 1.6 $\angle(B - A, D - B) \in (0, \pi)$, hence

$$\begin{vmatrix} b_1 - a_1 & b_2 - a_2 \\ d_1 - b_1 & d_2 - b_2 \end{vmatrix} = \begin{vmatrix} b_1 - a_1 & d_1 - b_1 \\ b_2 - a_2 & d_2 - b_2 \end{vmatrix} < 0,$$

so $\Delta_1 > 0$ and $\Delta_1/\Delta > 0$. Because $E$ lies strictly to the right of $\overrightarrow{A,B}$, by Statement 1.3(1) we obtain

$$\angle(B-A, C-A) = \angle(B-A, \frac{\Delta_1}{\Delta}(E-A)) = \angle(B-A, E-A) \in (0,\pi),$$

i.e. point $C$ also lies strictly to the right of the line $\overrightarrow{A,B}$. □

## 1.2. Angular sector and some of its properties

Let $\varphi \in (0,\pi)$, $E \neq (0,0)$. Denote

$$\boldsymbol{C}(E,\varphi) = \{C \in \mathbb{R}^2 \mid \angle(C,E) \in [-\varphi,\varphi]\}$$

– a two-dimensional angular sector. Note that $\mathbf{0} = (0,0) \in \boldsymbol{C}(E,\varphi)$, since $\angle((0,0),E) = 0$.

**Statement 1.24.** *Let (see fig. 1.5) $A, B$ be distinct points in $\mathbb{R}^2$, $\varphi \in (0,\pi)$. Then $B + \boldsymbol{C}(B-A,\varphi) \subseteq A + \boldsymbol{C}(B-A,\varphi)$.*

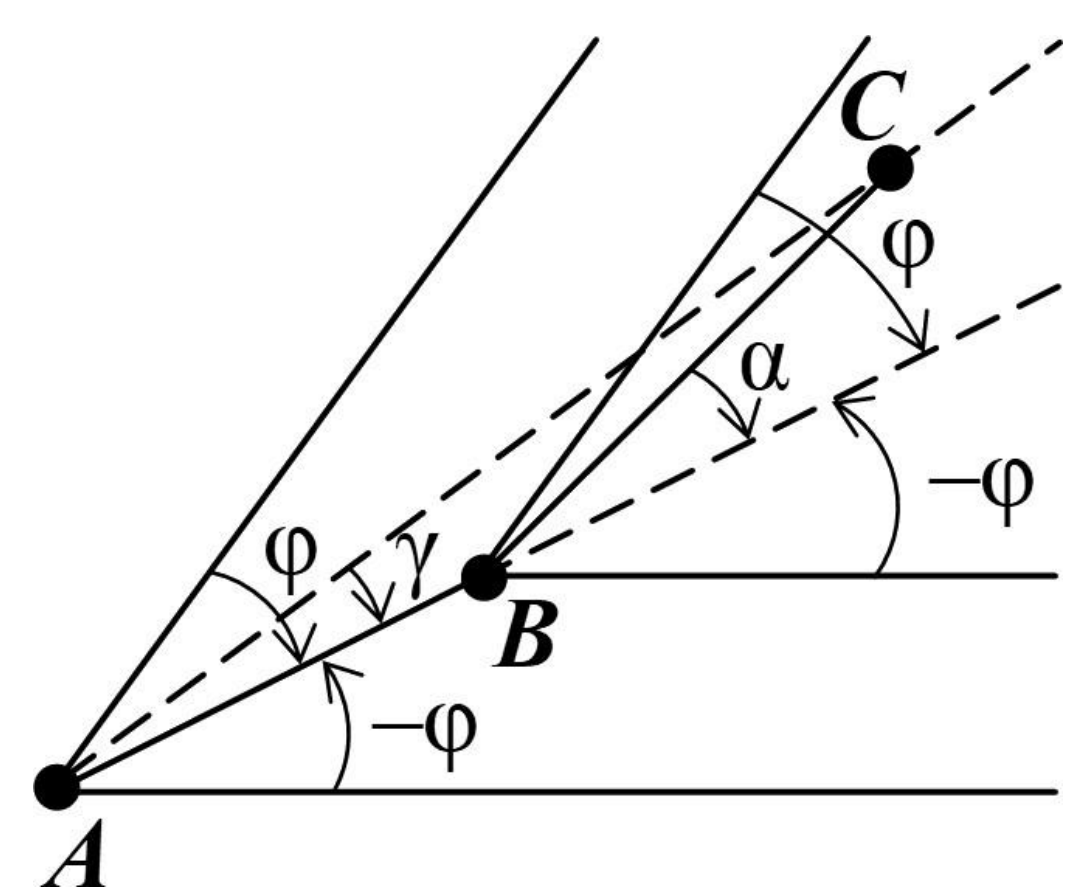

**Fig. 1.5.** The layout of points and angles following statement 1.24

***Proof.*** Let $C \in B + \boldsymbol{C}(B-A,\varphi)$ (see Fig. 1.5). Then

$$C - B \in \boldsymbol{C}(B-A,\varphi), \alpha = \angle(C-B,\, B-A) \in [-\varphi,\varphi].$$

In the case $\alpha = 0$ we have: $C - B = t(B-A)$ with $t \geq 0$, whence $C - A = (t+1)(B-A)$; consequently,

$$\angle(C-A,\, B-A) = \angle((t+1)(B-A),\, B-A) = 0 \in [-\varphi,\varphi].$$

Then $C - A \in \boldsymbol{C}(B - A, \varphi)$, whence $C \in A + \boldsymbol{C}(B - A, \varphi)$.

Now let $\alpha = \angle(C - B,\, B - A) \neq 0$ (hence $C \neq B$, also $C \neq A$, since $\alpha \in [-\varphi, \varphi]$, $\varphi \in (0, \pi)$). Consider the case $\alpha \in (0, \varphi]$ (the case $\alpha \in [-\varphi, 0)$ is treated similarly, using Statement 1.22). Then by Statement 1.21,

$$\gamma = \angle(C - A,\, B - A) \in (0, \alpha) \subseteq (0, \varphi),$$

whence $C - A \in \boldsymbol{C}(B - A, \varphi)$ or $C \in A + \boldsymbol{C}(B - A, \varphi)$. □

**Statement 1.25.** *Let* $A, B$ *be distinct points in* $\mathbb{R}^2$, $\varphi_1, \varphi_2 > 0$, $\varphi_1 + \varphi_2 \in (0, \pi)$, $C \in B + \boldsymbol{C}(B - A, \varphi_1)$, $C \neq B$, $D \in C + \boldsymbol{C}(C - B, \varphi_2)$. *Then*
$D \in B + \boldsymbol{C}(B - A, \varphi_1 + \varphi_2) \subseteq A + \boldsymbol{C}(B - A, \varphi_1 + \varphi_2)$.

***Proof.*** From the conditions of the statement and using Statement 1.24, we obtain:

$$\angle(C - B, B - A) \in [-\varphi_1, \varphi_1], \qquad (1.24)$$

$$D \in C + \boldsymbol{C}(C - B, \varphi_2) \subseteq B + \boldsymbol{C}(C - B, \varphi_2),$$

hence

$$D - B \in \boldsymbol{C}(C - B, \varphi_2) \Rightarrow \angle(D - B, C - B) \in [-\varphi_2, \varphi_2]. \qquad (1.25)$$

Then, by Statement 1.9 (since $D \neq B$ because $D \in C + \boldsymbol{C}(C - B, \varphi_2)$ and $\angle(B - C, C - B) = \pi < \varphi_2$, i.e. $B - C \notin \boldsymbol{C}(C - B, \varphi_2)$), we have

$$\angle(D - B, B - A) = \angle(D - B, C - B) \oplus \angle(C - B, B - A).$$

Using (1.24), (1.25), the fact that $\varphi_1 + \varphi_2 \in (0, \pi)$, and Statement 1.10, we obtain:

$$\angle(D - B, B - A) = \angle(D - B, C - B) + \angle(C - B, B - A) \in [-\varphi_1 - \varphi_2, \varphi_1 + \varphi_2],$$

hence (by Statement 1.24),

$$D - B \in \boldsymbol{C}(B - A, \varphi_1 + \varphi_2) \;\Rightarrow\; D \in B + \boldsymbol{C}(B - A, \varphi_1 + \varphi_2)$$
$$\subseteq A + \boldsymbol{C}(B - A, \varphi_1 + \varphi_2). \;□$$

**Statement 1.26.** *Let* $A, B$ *be distinct points in* $\mathbb{R}^2$, $\varphi_1, \varphi_2, \varphi_3 > 0$, $\varphi_1 + \varphi_2 + \varphi_3 \in (0,\pi)$, $C \in B + \boldsymbol{C}(B - A, \varphi_1)$, $C \neq B$, $D \in C + \boldsymbol{C}(C - B, \varphi_2)$, $D \neq C$, $E \in D + \boldsymbol{C}(D - C, \varphi_3)$. *Then*

$$E \in B + \boldsymbol{C}(B - A, \varphi_1 + \varphi_2 + \varphi_3) \subseteq A + \boldsymbol{C}(B - A, \varphi_1 + \varphi_2 + \varphi_3).$$

***Proof.*** From Statement 1.25, using the conditions

$$C \neq B, D \in C + \boldsymbol{C}(C - B, \varphi_2), D \neq C,$$
$$E \in D + \boldsymbol{C}(D - C, \varphi_3), \varphi_2 + \varphi_3 \in (0,\pi),$$

we obtain

$$E \in C + \boldsymbol{C}(C - B, \varphi_2 + \varphi_3). \qquad (1.26)$$

Then, again applying Statement 1.25, from the conditions $C \neq B, \varphi_1 + \varphi_2 + \varphi_3 \in (0,\pi)$, $C \in B + \boldsymbol{C}(B - A, \varphi_1)$, and (1.26), we obtain the desired statement. □

Using mathematical induction similarly to Statement 1.26, one can easily prove the general

**Statement 1.27.** *Let* $n \geq 1$, $A = B^{(0)}, B^{(1)}, \dots, B^{(n+1)} \in \mathbb{R}^2$, $B^{(i)} \neq B^{(i-1)}$, $i = 1, \dots, n+1$, $\varphi_1, \dots, \varphi_n > 0$, *and suppose*

$$B^{(i+1)} \in B^{(i)} + \boldsymbol{C}(B^{(i)} - B^{(i-1)}, \varphi_i), i = 1, \dots, n.$$

*Then, if* $\varphi_1 + \dots + \varphi_n \in (0,\pi)$, *we have*

$$B^{(n+1)} \in B^{(1)} + \boldsymbol{C}(B^{(1)} - A, \varphi_1 + \dots + \varphi_n) \subseteq A + \boldsymbol{C}(B^{(1)} - A, \varphi_1 + \dots + \varphi_n).$$

*Moreover, for any* $0 \leq j < i \leq n$ *with* $\varphi_{j+1} + \dots + \varphi_i \in (0,\pi)$ *the following is true:*

$$B^{(i+1)} \in B^{(j+1)} + \boldsymbol{C}(B^{(j+1)} - B^{(j)}, \varphi_{j+1} + \dots + \varphi_i)$$
$$\subseteq B^{(j)} + \boldsymbol{C}(B^{(j+1)} - B^{(j)}, \varphi_{j+1} + \dots + \varphi_i).$$

A corollary of Statement 1.27 is

**Statement 1.28.** *Suppose we are under conditions of Statement 1.27. Then, if* $\varphi_1 + \dots + \varphi_n \in (0,\pi)$, *we have* $B^{(n+1)} \neq B^{(0)}$. *Moreover, in the* *case* $0 \le j < i \le n$ *and* $\varphi_{j+1} + \cdots + \varphi_i \in (0, \pi)$ *the following holds:* $i + 1 \ge j + 2 \ge 2$ *,* $B^{(i+1)} \ne B^{(j)}$*.*

***Proof.*** We prove $B^{(n+1)} \ne B^{(0)}$ (the proof that $B^{(i+1)} \ne B^{(j)}$ is analogous). Suppose $B^{(n+1)} = B^{(0)}$. Then by Statement 1.27,

$$B^{(n+1)} = B^{(0)} \in B^{(1)} + \boldsymbol{C}(B^{(1)} - B^{(0)}, \varphi_1 + \cdots + \varphi_n),$$

hence

$$B^{(0)} - B^{(1)} \in \boldsymbol{C}(B^{(1)} - B^{(0)}, \varphi_1 + \cdots + \varphi_n),$$

which contradicts the fact that

$$\angle(B^{(0)} - B^{(1)}, B^{(1)} - B^{(0)}) = \pi > \varphi_1 + \cdots + \varphi_n. \square$$

**Remark 1.4.** Statement 1.28 can be strengthened by replacing the condition $B^{(n+1)} \ne B^{(0)}$ with the weaker condition

$$B^{(n+1)} \ne B^{(1)} + t(B^{(0)} - B^{(1)}) \text{ for all } t > 0.$$

Correspondingly, for $0 \le j < i \le n$ with $\varphi_{j+1} + \cdots + \varphi_i \in (0, \pi)$, the condition $B^{(i+1)} \ne B^{(j)}$ (where $i + 1 \ge j + 2 \ge 2$) can be replaced by

$$B^{(i+1)} \ne B^{(j+1)} + t(B^{(j)} - B^{(j+1)}) \text{ for all } t > 0.$$

### 1.3. Piecewise linear broken lines and some of their properties

First, we give an auxiliary statement.

**Statement 1.29.** *Suppose that for some numbers* $x, y, \alpha_1, \ldots, \alpha_n$ *we have*

$$\min\{0, \alpha_1, \ldots, \alpha_n\} \le x \le \max\{0, \alpha_1, \ldots, \alpha_n\},$$

*with* $y = 0$ *when* $x = 0$*, and when* $x \ne 0$ *we have* $|y| \in (0, |x|)$ *and* $y$ *has the same sign as* $x$*. Then*

$$\min\{0, \alpha_1, \ldots, \alpha_n\} \le y \le \max\{0, \alpha_1, \ldots, \alpha_n\}.$$

**Remark 1.5.** The presence of the number 0 in the bounds for $x, y$ is essential. For example, take $n = 2, \alpha_1 = -3, \alpha_2 = -4, x = -3, y = -2$. Then

$$\min\{-3, -4\} = -4 \le x = -3 = \max\{-3, -4\},$$

$|y| \in (0, |x|)$, $y$ has the same sign as $x$, but the inequalities

$$\min\{-3,-4\} = -4 \le y = -2 \le \max\{-3,-4\} = -3$$

do not hold. However, for the same $x = -3, y = -2$ we have

$$\min\{0,-3,-4\} = -4 \le x = -3 \le 0 = \max\{0,-3,-4\},$$

$$\min\{0,-3,-4\} = -4 \le y = -2 \le 0 = \max\{0,-3,-4\}.$$

**Statement 1.30.** *Let $n \ge 2$, $A = B^{(0)}, B^{(1)}, B^{(2)}, \dots, B^{(n)}, B^{(n+1)} = B$ be points in $\mathbb{R}^2$,*

$\alpha_i = \angle(B^{(i+1)} - B^{(i)}, B^{(i)} - B^{(i-1)}) \notin \{0, \pi\}$,

$\gamma_i = \angle(B^{(i+1)} - A, B^{(1)} - A), \bar{\gamma}_i = \angle(B^{(i+1)} - B^{(1)}, B^{(1)} - A), i = 1, \dots, n,$

$\beta_{i,j} = \angle(B^{(i+1)} - B^{(j)}, B^{(j+1)} - B^{(j)}), \bar{\beta}_{i,j} = \angle(B^{(i+1)} - B^{(j+1)}, B^{(j+1)} - B^{(j)}),$

$i = 1, \dots, n,\ j = 0, \dots, i-1$

*(see Figs. 1.6, 1.7). Then*

1. *$B^{(i)} \ne B^{(i-1)}$, $i = 1, \dots, n+1$.*
2. *The angle $\gamma_1 = \angle(B^{(2)} - A, B^{(1)} - A)$ has the same sign as $\alpha_1 = \bar{\gamma}_1 = \angle(B^{(2)} - B^{(1)}, B^{(1)} - B^{(0)}) \ne 0$ and $|\gamma_1| \in (0, |\alpha_1|)$.*
3. *For each $i \in \{2, \dots, n\}$ we have:*

   *(a) If $\bar{\gamma}_i = 0$, then $B^{(i+1)} \in \overrightarrow{A, B^{(1)}}$ and $\gamma_i = \angle(B^{(i+1)} - A, B^{(1)} - A) = \angle(B^{(i+1)} - B^{(1)}, B^{(1)} - A) = \bar{\gamma}_i = 0$.*

   *(b) If $\bar{\gamma}_i = \pi$, then $B^{(i+1)} \in \overrightarrow{A, B^{(1)}}$ and for some $t \in \mathbb{R}$ the following is true*

$$B^{(i+1)} = B^{(1)} + t(A - B^{(1)}) = B^{(1)} + t(B^{(0)} - B^{(1)}), t > 0. \qquad (1.27)$$

   *(c) If $\bar{\gamma}_i \notin \{0, \pi\}$ (i.e. $\bar{\gamma}_i \in (-\pi, 0) \cup (0, \pi)$), then $\gamma_i = \angle(B^{(i+1)} - A, B^{(1)} - A)$ has the same sign as $\bar{\gamma}_i = \angle(B^{(i+1)} - B^{(1)}, B^{(1)} - A)$ and $|\gamma_i| \in (0, |\bar{\gamma}_i|)$.*

4. *For each $i \in \{1, \dots, n\}$ we have*

$$\beta_{i,0} = \gamma_i, \bar{\beta}_{i,0} = \angle(B^{(i+1)} - B^{(1)}, B^{(1)} - B^{(0)}) = \bar{\gamma}_i,$$

$$\bar{\beta}_{i,i-1} = \angle\left(B^{(i+1)} - B^{(i)}, B^{(i)} - B^{(i-1)}\right) = \alpha_i. \qquad (1.28)$$

*5. For each* $i \in \{1, \dots, n\}, j \in \{0, \dots, i-1\}$ *the following hold:*

*(a) If* $\bar{\beta}_{i,j} = 0$*, then* $B^{(i+1)} \in \overrightarrow{B^{(j)}, B^{(j+1)}}$ *and* $\beta_{i,j} = \bar{\beta}_{i,j} = 0$*.*

*(b) If* $\bar{\beta}_{i,j} = \pi$*, then* $i \geq j + 2 \geq 2$ *and for some* $t \in \mathbb{R}$ *the following is true*

$$B^{(i+1)} = B^{(j+1)} + t(B^{(j)} - B^{(j+1)}), t > 0. \qquad (1.29)$$

*(c) If* $\bar{\beta}_{i,j} \notin \{0, \pi\}$ *(i.e.* $\bar{\beta}_{i,j} \in (-\pi, 0) \cup (0, \pi)$*), then* $\beta_{i,j}$ *has the same sign as* $\bar{\beta}_{i,j}$ *and* $|\beta_{i,j}| \in (0, |\bar{\beta}_{i,j}|)$*.*

*6. For any* $i \in \{2, \dots, n\}, j \in \{0, \dots, i-2\}$*, if* $B^{(i+1)} \neq B^{(j+1)}$*, the recurrence formula holds:*

$$\bar{\beta}_{i,j} = \angle\left(B^{(i+1)} - B^{(j+1)}, B^{(j+1)} - B^{(j)}\right) =$$

$$\angle\left(B^{(i+1)} - B^{(j+1)}, B^{(j+2)} - B^{(j+1)}\right) \oplus \angle\left(B^{(j+2)} - B^{(j+1)}, B^{(j+1)} - B^{(j)}\right) =$$

$$\beta_{i,j+1} \oplus \alpha_{j+1}. \qquad (1.30)$$

*7. For any* $i \in \{2, \dots, n\}$*, if* $B^{(i+1)} \neq B^{(1)}$*, then*

$$\bar{\gamma}_i = \angle(B^{(i+1)} - B^{(1)}, B^{(1)} - A) = \angle(B^{(i+1)} - B^{(1)}, B^{(2)} - B^{(1)}) \oplus$$

$$\angle(B^{(2)} - B^{(1)}, B^{(1)} - A) = \beta_{i,1} \oplus \alpha_1. \qquad (1.31)$$

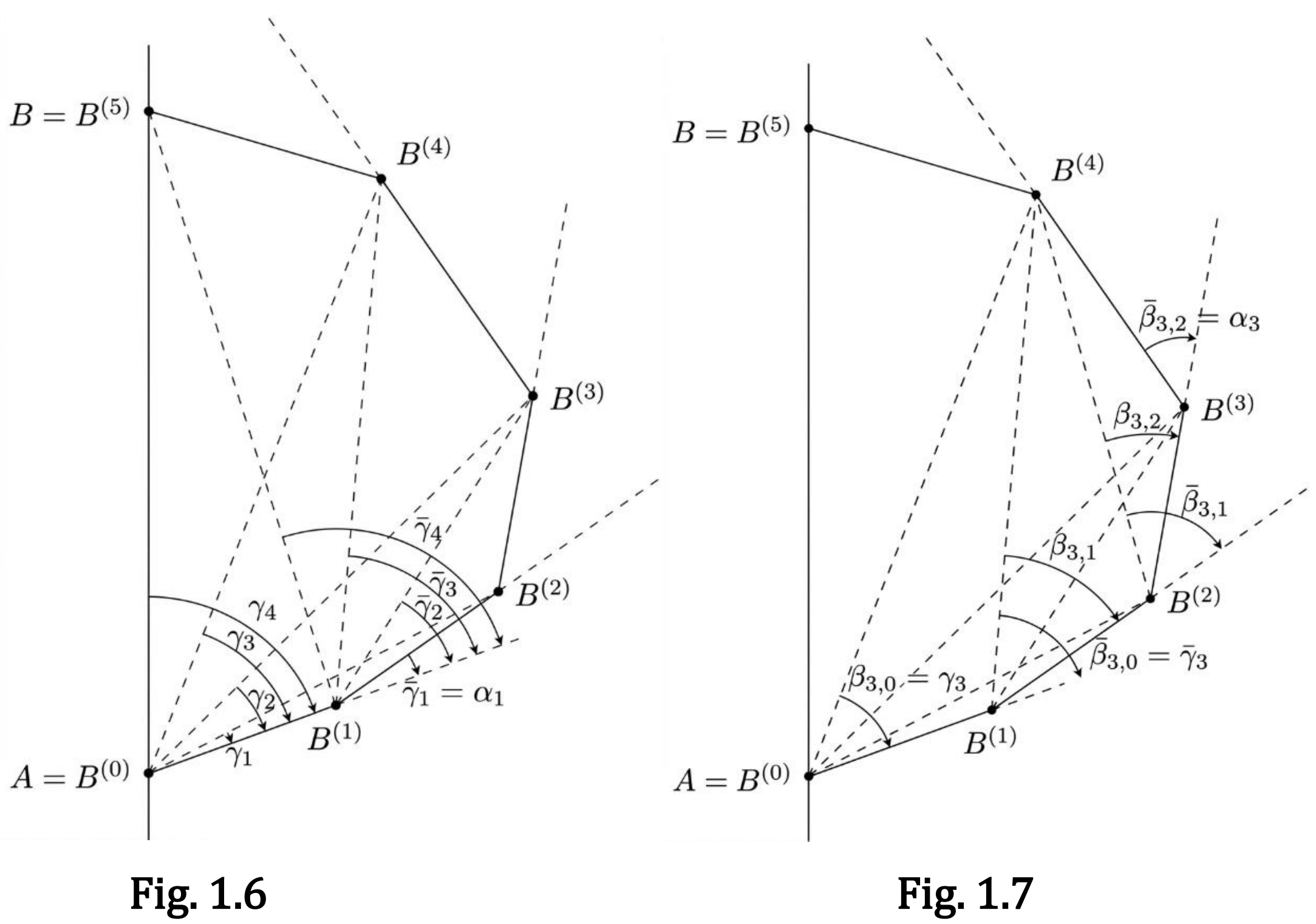

**Fig. 1.6** **Fig. 1.7**

***Proof.*** 1) Follows from $\alpha_i = \angle(B^{(i+1)} - B^{(i)}, B^{(i)} - B^{(i-1)}) \neq 0, i = 1, \dots, n$. 2) Follows from Statements 1.21 and 1.22 and the condition $\alpha_1 \notin \{0, \pi\}$. We prove 3). (a) If $\bar{\gamma}_i = \angle(B^{(i+1)} - B^{(1)}, B^{(1)} - A) = 0$, then by Statement 1.1, $B^{(i+1)} \in \overrightarrow{A, B^{(1)}}$ and for some $t \geq 0$ we have $B^{(i+1)} - B^{(1)} = t(B^{(1)} - A)$, hence $B^{(i+1)} - A = (t+1)(B^{(1)} - A)$, and therefore

$$\gamma_i = \angle(B^{(i+1)} - A, B^{(1)} - A) = \angle((t+1)(B^{(1)} - A), B^{(1)} - A) = 0.$$

(b) If $\bar{\gamma}_i = \angle(B^{(i+1)} - B^{(1)}, B^{(1)} - A) = \pi$, then by Statement 1.1, $B^{(i+1)} \in \overrightarrow{A, B^{(1)}}$ and for some $t \in \mathbb{R}$ we have

$$B^{(i+1)} - B^{(1)} = -t(B^{(1)} - A), t > 0,$$

which yields (1.27).

(c) Follows from Statements 1.21 and 1.22.

Statement 4) is obvious. Let us prove Statement 5).

(a) If $\bar{\beta}_{i,j} = \angle\left(B^{(i+1)} - B^{(j+1)}, B^{(j+1)} - B^{(j)}\right) = 0$, then by virtue of

Statement 1.1 we have $B^{(i+1)} \in \overrightarrow{B^{(j)}, B^{(j+1)}}$ and for some $t \geq 0$ the equality $B^{(i+1)} - B^{(j+1)} = t(B^{(j+1)} - B^{(j)})$ holds, whence $B^{(i+1)} - B^{(j)} = (t + 1)(B^{(j+1)} - B^{(j)})$, and consequently

$$\beta_{i,j} = \angle(B^{(i+1)} - B^{(j)}, B^{(j+1)} - B^{(j)}) =$$
$$\angle((t+1)(B^{(j+1)} - B^{(j)}), B^{(j+1)} - B^{(j)}) = 0.$$

(b) If $\bar{\beta}_{i,j} = \angle\left(B^{(i+1)} - B^{(j+1)}, B^{(j+1)} - B^{(j)}\right) = \pi$, then by Statement 1.1 we have $B^{(i+1)} \in \overrightarrow{B^{(j)}, B^{(j+1)}}$ and for some $t \in \mathbb{R}$ the following holds:

$$B^{(i+1)} - B^{(j+1)} = -t(B^{(j+1)} - B^{(j)}),\ t > 0,$$

from which (1.29) follows. Moreover, since $\bar{\beta}_{i,i-1} = \alpha_i \notin \{0, \pi\}$, we have $j \neq i - 1$, i.e., $j \leq i - 2$, whence $i \geq j + 2 \geq 2$.

(c) If $\bar{\beta}_{i,j} \notin \{0, \pi\}$ (i.e., $\bar{\beta}_{i,j} \in (-\pi, 0) \cup (0, \pi)$), then by Statements 1.21, 1.22 the angle $\beta_{i,j} = \angle(B^{(i+1)} - B^{(j)}, B^{(j+1)} - B^{(j)})$ has the same sign as the angle $\bar{\beta}_{i,j} = \angle\left(B^{(i+1)} - B^{(j+1)}, B^{(j+1)} - B^{(j)}\right)$ and moreover $|\beta_{i,j}| \in (0, |\bar{\beta}_{i,j}|)$.

Statements 6) and 7) follow from Statement 1.9 (addition of angles). □

**Statement 1.31.** *Assume we are under the conditions of Statement 1.30 and*

$$\theta_i = \angle(B^{(i+1)} - B^{(i)}, B^{(i+1)} - A),$$
$$\bar{\theta}_i = \angle(B^{(i+1)} - B^{(i)}, B^{(i)} - A),\ i = 1, \dots, n$$

*(see Fig. 1.8). Then*

1. *Angle* $\theta_1 = \angle(B^{(2)} - B^{(1)}, B^{(2)} - A)$ *has the same sign as* $\alpha_1 = \bar{\theta}_1 = \angle(B^{(2)} - B^{(1)}, B^{(1)} - B^{(0)}) \neq 0$ *and moreover* $|\theta_1| \in (0, |\alpha_1|)$.
2. *For each* $i \in \{2, \dots, n\}$ *we have:*

   *(a) If* $\bar{\theta}_i = 0$*, then* $B^{(i+1)} \in \overrightarrow{A, B^{(i)}}$ *and*

   $\theta_i = \angle(B^{(i+1)} - B^{(i)}, B^{(i+1)} - A) = \angle(B^{(i+1)} - B^{(i)}, B^{(i)} - A) = \bar{\theta}_i = 0$.

*(b) If* $\bar{\theta}_i = \pi$*, then* $B^{(i+1)} \in \overrightarrow{A, B^{(i)}}$ *and moreover* $\theta_i \in \{0, \pi\}$*.*

*(c) If* $\bar{\theta}_i \notin \{0, \pi\}$ *(i.e.,* $\bar{\theta}_i \in (-\pi, 0) \cup (0, \pi)$*), then the angle* $\theta_i = \angle(B^{(i+1)} - B^{(i)}, B^{(i+1)} - A)$ *has the same sign as* $\bar{\theta}_i = \angle(B^{(i+1)} - B^{(i)}, B^{(i)} - A)$*, and moreover* $|\theta_i| \in (0, |\bar{\theta}_i|)$*.*

3. *For any* $i \in \{2, \dots, n\}$ *in the case* $B^{(i)} \neq A$ *we have*

$$\bar{\theta}_i = \angle(B^{(i+1)} - B^{(i)}, B^{(i)} - A) = \angle(B^{(i+1)} - B^{(i)}, B^{(i)} - B^{(i-1)}) \oplus \angle(B^{(i)} - B^{(i-1)}, B^{(i)} - A) = \theta_{i-1} \oplus \alpha_i. \quad (1.32)$$

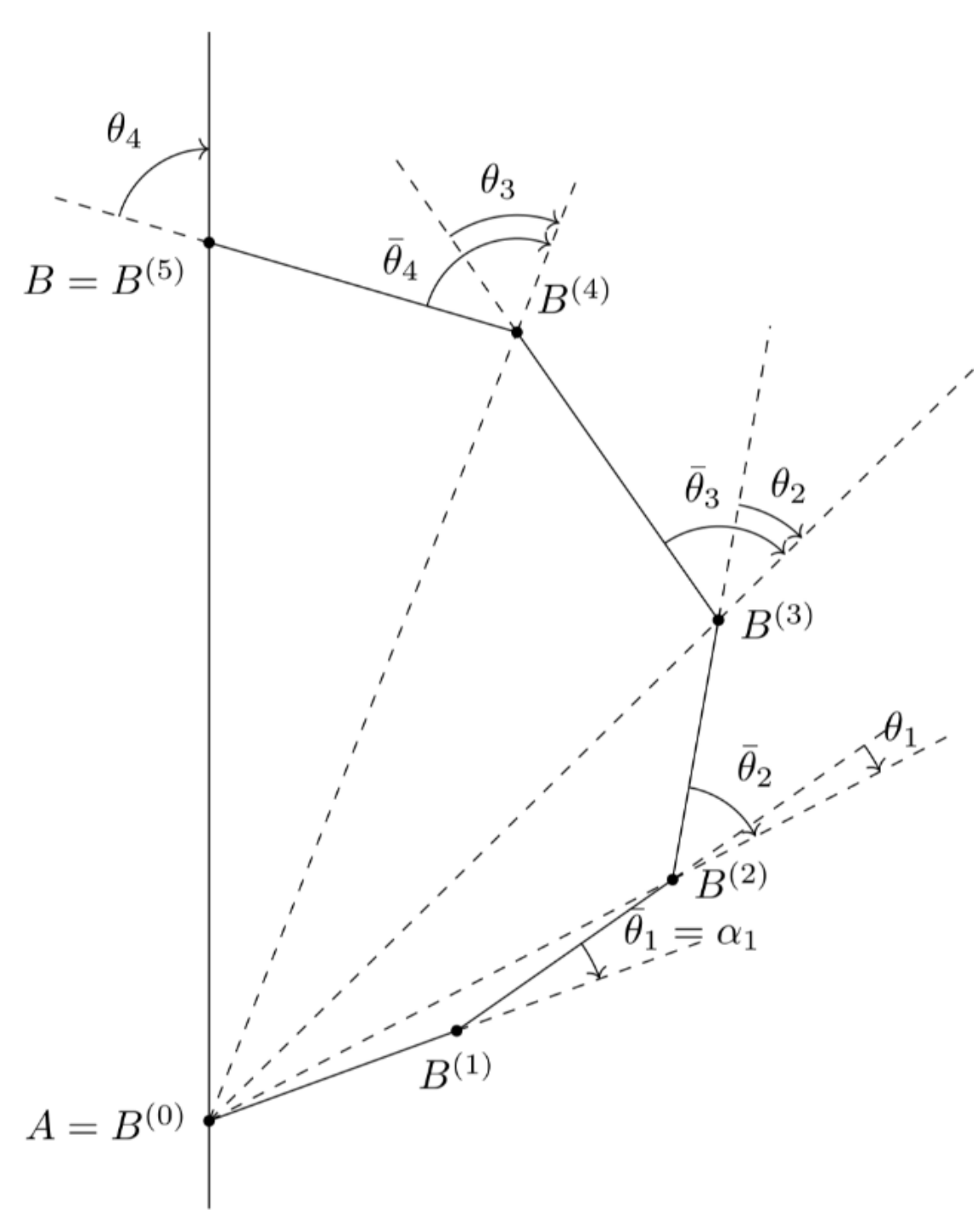

**Fig. 1.8.** The layout of points and angles following statement 1.31

***Proof.*** Statement 1) follows from Statements 1.21, 1.22, and also from the condition $\alpha_1 \notin \{0, \pi\}$. Let us prove Statement 1.2).

(a) If $\bar{\theta}_i = 0$, then by Statement 1.1 we have $B^{(i+1)} \in \overrightarrow{A, B^{(i)}}$ and for some $t \geq 0$ the equality $B^{(i+1)} - B^{(i)} = t(B^{(i)} - A)$ holds, whence $B^{(i+1)} - A = (t+1)(B^{(i)} - A)$, and consequently

$$\theta_i = \angle(B^{(i+1)} - B^{(i)}, B^{(i+1)} - A) = \angle(t(B^{(i)} - A), (t+1)(B^{(i)} - A)) = 0.$$

(b) If $\bar{\theta}_i = \pi$, then by Statement 1.1 we have $B^{(i+1)} \in \overrightarrow{A, B^{(i)}}$ and for some $t \in \mathbb{R}$ the following holds: $B^{(i+1)} - B^{(i)} = t(A - B^{(i)})$, $t > 0$; consequently, $B^{(i+1)} - A = (t-1)(A - B^{(i)})$, $t > 0$, and moreover $\theta_i = \angle(B^{(i+1)} - B^{(i)}, B^{(i+1)} - A) = \angle(t(A - B^{(i)}), (t-1)(A - B^{(i)}))$, whence $\theta_i = \pi$ for $t \in (0,1)$ and $\theta_i = 0$ for $t \geq 1$.

(c) Follows from Statements 1.21, 1.22.

Statement 3) follows from Statement 1.9 taking into account that $B^{(i+1)} \neq B^{(i)}, i = 0, \dots, n$. □

We will need some notation. Let $\alpha_1, \dots, \alpha_n \in \mathbb{R}$ be arbitrary numbers. For any $i \in \{1, \dots, n\}, j \in \{0, \dots, i-1\}$, denote:

$$\mu_{i,j} = \mu_{i,j}(\alpha_1, \dots, \alpha_n) = \max\{0, \alpha_{j+1}, \alpha_{j+1} + \alpha_{j+2}, \dots, \alpha_{j+1} + \dots + \alpha_i\},$$
$$\bar{\mu}_{i,j} = \bar{\mu}_{i,j}(\alpha_1, \dots, \alpha_n) = \min\{0, \alpha_{j+1}, \alpha_{j+1} + \alpha_{j+2}, \dots, \alpha_{j+1} + \dots + \alpha_i\},$$
$$\sigma_{i,j} = \sigma_{i,j}(\alpha_1, \dots, \alpha_n) = |\alpha_{j+1}| + \dots + |\alpha_i|,$$
$$\mu_i = \mu_i(\alpha_1, \dots, \alpha_n) = \mu_{i,0} = \max\{0, \alpha_1, \alpha_1 + \alpha_2, \dots, \alpha_1 + \dots + \alpha_i\},$$
$$\bar{\mu}_i = \bar{\mu}_i(\alpha_1, \dots, \alpha_n) = \bar{\mu}_{i,0} = \min\{0, \alpha_1, \alpha_1 + \alpha_2, \dots, \alpha_1 + \dots + \alpha_i\},$$
$$\sigma_i = \sigma_i(\alpha_1, \dots, \alpha_n) = \sigma_{i,0} = |\alpha_1| + \dots + |\alpha_i|. \quad (1.33)$$

**Example 1.3.** 1) Let $n = 2$. Then

$\sigma_1 = |\alpha_1|, \sigma_2 = |\alpha_1| + |\alpha_2|,$

$\sigma_{1,0} = \sigma_1 = |\alpha_1|, \sigma_{2,0} = \sigma_2 = |\alpha_1| + |\alpha_2|, \sigma_{2,1} = |\alpha_2|,$

$\mu_1 = \max\{0, \alpha_1\}, \mu_2 = \max\{0, \alpha_1, \alpha_1 + \alpha_2\},$

$\mu_{1,0} = \mu_1 = \max\{0, \alpha_1\}, \mu_{2,0} = \mu_2 = \max\{0, \alpha_1, \alpha_1 + \alpha_2\}, \mu_{2,1} = \max\{0, \alpha_2\},$

$\bar{\mu}_1 = \min\{0, \alpha_1\}, \bar{\mu}_2 = \min\{0, \alpha_1, \alpha_1 + \alpha_2\},$

$\bar{\mu}_{1,0} = \bar{\mu}_1 = \min\{0, \alpha_1\}, \bar{\mu}_{2,0} = \bar{\mu}_2 = \min\{0, \alpha_1, \alpha_1 + \alpha_2\}, \bar{\mu}_{2,1} = \min\{0, \alpha_2\}.$

2. Let $n = 3, \alpha_1 = -1, \alpha_2 = 1, \alpha_3 = 1$. Then

$\mu_{3,1} = \max\{0, \alpha_2, \alpha_2 + \alpha_3\} = \max\{0,1,1+1\} = 2,$

$$\mu_1 = \max\{0, \alpha_1\} = 0,$$

$$\mu_{3,1} + \mu_1 = 2,$$

$$\mu_3 = \max\{0, \alpha_1, \alpha_1 + \alpha_2, \alpha_1 + \alpha_2 + \alpha_3\}$$
$$= \max\{0, -1, -1 + 1, -1 + 1 + 1\} = 1,$$

$$\max\{0, \alpha_2, \alpha_2 + \alpha_3\} + \max\{0, \alpha_1\} = \mu_{3,1} + \mu_1 = 2 > 1 = \mu_3$$
$$= \max\{0, \alpha_1, \alpha_1 + \alpha_2, \alpha_1 + \alpha_2 + \alpha_3\},$$

$$\mu_{3,1} + \alpha_1 = 1,$$

$$\mu_{3,1} + \alpha_1 = \max\{0, \alpha_2, \alpha_2 + \alpha_3\} + \alpha_1 = \max\{\alpha_1, \alpha_1 + \alpha_2, \alpha_1 + \alpha_2 + \alpha_3\}$$
$$\leq \max\{0, \alpha_1, \alpha_1 + \alpha_2, \alpha_1 + \alpha_2 + \alpha_3\} = \mu_{3,0} = \mu_3.$$

The following statement is obvious.

**Statement 1.32.** *Let* $\alpha_1, \dots, \alpha_n \in \mathbb{R}$ *be arbitrary numbers and let equalities (1.33) hold. Then for any* $i \in \{1, \dots, n\}, j \in \{0, \dots, i-1\}$*:*

$$\sigma_j + \sigma_{i,j} = \sigma_i,$$

$$-\sigma_i \leq \bar{\mu}_i \leq \mu_i \leq \sigma_i, -\sigma_{i,j} \leq \bar{\mu}_{i,j} \leq \mu_{i,j} \leq \sigma_{i,j},$$

*and, moreover, for any* $i \in \{2, \dots, n\}, j \in \{1, \dots, i-1\}$*:*

$$\sigma_{i,j} + \alpha_j = \alpha_j + |\alpha_{j+1}| + \dots + |\alpha_i| \leq \sigma_{i,j-1}, \qquad (1.34)$$

$$-\sigma_{i,j-1} \leq \alpha_j - |\alpha_{j+1}| - \dots - |\alpha_i| = -\sigma_{i,j} + \alpha_j,$$

$$-\sigma_i = -\sigma_j - \sigma_{i,j} \leq \bar{\mu}_j + \bar{\mu}_{i,j} \leq \mu_j + \mu_{i,j} \leq \sigma_j + \sigma_{i,j} = \sigma_i,$$

$$\mu_{i,j} + \alpha_j = \max\{\alpha_j, \alpha_j + \alpha_{j+1}, \dots, \alpha_j + \alpha_{j+1} + \dots + \alpha_i\} \leq \mu_{i,j-1}, \qquad (1.35)$$

$$\bar{\mu}_{i,j-1} \leq \bar{\mu}_{i,j} + \alpha_j = \min\{\alpha_j, \alpha_j + \alpha_{j+1}, \dots, \alpha_j + \alpha_{j+1} + \dots + \alpha_i\}. \qquad (1.36)$$

**Statement 1.33.** *Suppose we are under the conditions of Statement 1.30 and for some* $i \in \{2, \dots, n\}$ *we have*

$$\sigma_i = \sigma_i(\alpha_1, \dots, \alpha_n) = \sum_{k=1}^{i} |\alpha_k| < \pi. \qquad (1.37)$$

*Then*

*1.* $\bar{\beta}_{i,j} \neq \pi$ *for* $j \in \{0, \dots, i-1\}$.

*2. For* $j \in \{0, \dots, i-1\}$ *the following inequalities hold:*

$$-\sigma_{i,j} \leq \bar{\mu}_{i,j} \leq \bar{\beta}_{i,j} \leq \mu_{i,j} \leq \sigma_{i,j}, -\sigma_{i,j} \leq \bar{\mu}_{i,j} \leq \beta_{i,j} \leq \mu_{i,j} \leq \sigma_{i,j}, \quad (1.38)$$

*3. The inequalities hold:*

$$-\sigma_i \leq \bar{\mu}_i = \bar{\mu}_{i,0} \leq \bar{\gamma}_i = \bar{\beta}_{i,0} \leq \mu_{i,0} = \mu_i \leq \sigma_i, \quad (1.39)$$

$$-\sigma_i \leq \bar{\mu}_i \leq \gamma_i \leq \mu_i \leq \sigma_i. \quad (1.40)$$

***Proof.*** From (1.37) and from Statement 1.30(1), it follows that for $\varphi_k = |\alpha_k|, k = 1, \dots, i$, we are under the conditions of Statements 1.27, 1.28, and also Remark 1.4. We prove 1). If for some $j \in \{1, \dots, i-2\}$ ( $\bar{\beta}_{i,i-1} = \alpha_i \notin \{0, \pi\}$ ) we have $\bar{\beta}_{i,j} = \pi$, then due to Statement 1.30(5) holds (1.29), which contradicts the assertion of Remark 1.4. We prove 2). Let $i \in \{2, \dots, n\}$. We prove inequalities (1.38) for $\bar{\beta}_{i,j}$, where $j \in \{0, \dots, i-1\}$. Then by Statement 1.29 and Statement 1.30(5), they will also hold for $\beta_{i,j}$ (for example, in the case of $\bar{\beta}_{i,j} = 0$, we have $\beta_{i,j} = 0, \bar{\mu}_{i,j} \leq 0 \leq \mu_{i,j}$). For $j = i-1$ these inequalities become

$$-|\alpha_i| \leq \min\{0, \alpha_i\} \leq \bar{\beta}_{i,i-1} \leq \max\{0, \alpha_i\} \leq |\alpha_i|,$$

and follow from the equality $\bar{\beta}_{i,i-1} = \alpha_i$ (see Statement 1.30(4)). Assume the inequalities (1.38) hold for some $\bar{\beta}_{i,j}$ with $j \in \{1, \dots, i-1\}$. We show they will also hold for $j-1 \in \{0, \dots, i-2\}$. Since inequalities (1.38) hold for $\bar{\beta}_{i,j}$, as already noted, they will also hold for $\beta_{i,j}$. Note that by Statement 1.28 we cannot have $B^{(i)} = B^{(j)}$, and consequently, by Statement 1.30(6) (see (1.30)) the equality $\bar{\beta}_{i,j-1} = \beta_{i,j} \oplus \alpha_j$ holds . Note that from (1.38), using (1.34) and (1.37), we have

$$-\pi < -\sigma_{i,j-1} \leq -\sigma_{i,j} + \alpha_j \leq \beta_{i,j} + \alpha_j \leq \sigma_{i,j} + \alpha_j \leq \sigma_{i,j-1} < \pi,$$

hence by Statement 1.10, $\bar{\beta}_{i,j-1} = \beta_{i,j} \oplus \alpha_j = \beta_{i,j} + \alpha_j$, and then from (1.38) together with (1.35), (1.36) we obtain

$$\bar{\beta}_{i,j-1} = \beta_{i,j} + \alpha_j \le \mu_{i,j} + \alpha_j \le \mu_{i,j-1} \le \sigma_{i,j-1},$$
$$-\sigma_{i,j-1} \le \bar{\mu}_{i,j-1} \le \bar{\mu}_{i,j} + \alpha_j \le \beta_{i,j} + \alpha_j = \bar{\beta}_{i,j-1},$$

i.e., inequalities (1.38) are proved for all $j \in \{0, \dots, i-1\}$.

We prove 3). Inequality (1.39) follows from (1.38) with $j = 0$. Inequality (1.40) follows from (1.39), Statement 1.29, and Statement 1.30(3). Moreover, from (1.37) and (1.39) we obtain that the case of Statement 1.30(3,b), where $\bar{\gamma}_i = \pi$, is impossible. □

We will need additional notation. Let $\alpha_1, \dots, \alpha_n \in \mathbb{R}$ be arbitrary numbers. For any $i \in \{1, \dots, n\}$ denote

$$\vartheta_i = \vartheta_i(\alpha_1, \dots, \alpha_n) = \max\{0, \alpha_i, \alpha_i + \alpha_{i-1}, \dots, \alpha_i + \alpha_{i-1} + \dots + \alpha_1\},$$
$$\bar{\vartheta}_i = \bar{\vartheta}_i(\alpha_1, \dots, \alpha_n) = \min\{0, \alpha_i, \alpha_i + \alpha_{i-1}, \dots, \alpha_i + \alpha_{i-1} + \dots + \alpha_1\}.$$

Obviously, for any $i \in \{1, \dots, n\}$

$$-|\alpha_1| - \dots - |\alpha_i| = -\sigma_i \le \bar{\vartheta}_i \le \vartheta_i \le \sigma_i = |\alpha_1| + \dots + |\alpha_i|.$$

Furthermore, for any $i \in \{1, \dots, n-1\}$

$$\vartheta_i + \alpha_{i+1} = \max\{\alpha_{i+1}, \alpha_{i+1} + \alpha_i, \alpha_{i+1} + \alpha_i + \alpha_{i-1}, \dots, \alpha_{i+1} + \alpha_i + \dots + \alpha_1\} \le \vartheta_{i+1}, \qquad (1.41)$$
$$\bar{\vartheta}_i + \alpha_{i+1} = \min\{\alpha_{i+1}, \alpha_{i+1} + \alpha_i, \alpha_{i+1} + \alpha_i + \alpha_{i-1}, \dots, \alpha_{i+1} + \alpha_i + \dots + \alpha_1\} \ge \bar{\vartheta}_{i+1}. \qquad (1.42)$$

**Statement 1.34.** *Suppose we are under the conditions of Statement 1.30 and for some* $i \in \{2, \dots, n\}$ *condition (1.37) holds. Then*

1. *For* $j = 1, \dots, i$ *the following inequalities hold:*

$$-\sigma_j \le \bar{\vartheta}_j \le \bar{\theta}_j = \angle(B^{(j+1)} - B^{(j)}, B^{(j)} - A) \le \vartheta_j \le \sigma_j, \qquad (1.43)$$
$$-\sigma_j \le \bar{\vartheta}_j \le \theta_j = \angle(B^{(j+1)} - B^{(j)}, B^{(j+1)} - A) \le \vartheta_j \le \sigma_j. \qquad (1.44)$$

2. *Let* $i \in \{1, \dots, n\}, j \in \{0, \dots, i\}$*, and let (1.37) hold. Then* $B^{(i+1)} \neq B^{(j)}$*, and the inequalities hold*:

$$-\sigma_i \leq \angle(B^{(i+1)} - B^{(j)}, B^{(j)} - A) = \bar{\beta}_{i,j-1} \oplus \theta_{j-1} = \bar{\beta}_{i,j-1} + \theta_{j-1} \leq \sigma_i. \quad (1.45)$$

**Proof.** We prove part 1) by induction on $j = 1, \dots, i$. For $j = 1$, by Statement 1.31(1), and, taking Statement 1.29 into account, we have

$$\min\{0, \alpha_1\} \leq \bar{\theta}_1 = \alpha_1 \leq \max\{0, \alpha_1\}, \min\{0, \alpha_1\} \leq \theta_1 = \alpha_1 \leq \max\{0, \alpha_1\},$$

whence, using the inequalities $\vartheta_1 = \max\{0, \alpha_1\} \leq |\alpha_1| = \sigma_1, -\sigma_1 = -|\alpha_1| \leq \min\{0, \alpha_1\} = \bar{\vartheta}_1$, we obtain the validity of (1.43), (1.44) for $j = 1$. Suppose (1.43), (1.44) hold for some index $j - 1 \in \{1, \dots, i-1\}$. We show the validity of (1.43), (1.44) for index $j \in \{2, \dots, i\}$. From the inductive hypothesis it follows that

$$-\sigma_{j-1} \leq \bar{\vartheta}_{j-1} \leq \bar{\theta}_{j-1} \leq \vartheta_{j-1} \leq \sigma_{j-1}, -\sigma_{j-1} \leq \bar{\vartheta}_{j-1} \leq \theta_j \leq \vartheta_{j-1} \leq \sigma_{j-1}. \quad (1.46)$$

By (1.37) (taking into account Statement 1.28) we have $B^{(j)} \neq A, j = 1, \dots, i$. Therefore, we can use the recurrence formula (1.32), which, given condition (1.37) and inequalities (1.46), takes the form: $\bar{\theta}_j = \theta_{j-1} + \alpha_j, j = 1, \dots, i$. But then, using (1.46), (1.41), (1.42), we obtain

$$-\sigma_j \leq \bar{\vartheta}_j \leq \bar{\vartheta}_{j-1} + \alpha_j \leq \bar{\theta}_j \leq \vartheta_{j-1} + \alpha_j \leq \vartheta_j \leq \sigma_j.$$

Next, using Statement 1.31(2), and Statement 1.29, we have $-\sigma_j \leq \bar{\vartheta}_j \leq \theta_j \leq \vartheta_j \leq \sigma_j$.

To prove 2), note that by (1.32) we are under the conditions of Statement 1.28, by which $B^{(i+1)} \neq B^{(j)}$. But then by Statement 1.9 the formula holds:

$$\angle(B^{(i+1)} - B^{(j)}, B^{(j)} - A) = \angle(B^{(i+1)} - B^{(j)}, B^{(j)} - B^{(j-1)}) \oplus \angle\left(B^{(j)} - B^{(j-1)}, B^{(j)} - A\right) = \bar{\beta}_{i,j-1} \oplus \theta_{j-1}. \quad (1.47)$$

From (1.37), taking into account inequalities (1.44), (1.38), and also Statement 1.32, we obtain

$$-\pi < -\sigma_i = -\sigma_{i,j-1} - \sigma_{j-1} \le \bar{\beta}_{i,j-1} + \theta_{j-1} \le \sigma_{i,j-1} + \sigma_{j-1} = \sigma_i < \pi, \quad (1.48)$$

hence by Statement 1.10 we have $\bar{\beta}_{i,j-1} \oplus \theta_{j-1} = \bar{\beta}_{i,j-1} + \theta_{j-1}$, and also (1.45). □

**Statement 1.35.** *Suppose we are under the conditions of Statement 1.30 and* $1 \le r \le n, \alpha_i \in [-\varphi, \varphi], i = 1,2, \dots, r, r\varphi < \pi$. *Suppose the point* $B^{(r-1)}$ *lies not strictly to the left of the line* $\overrightarrow{A, B^{(r)}}$ *and the point* $B^{(r+1)}$ *lies strictly to the left of the line* $\overrightarrow{A, B^{(r)}}$ *(i.e., due to the Statements 1.3. 1.6,* $\angle(B^{(r)} - A,\ B^{(r+1)} - B^{(r)}) \in (-\pi, 0)$*). Then*

1. $\angle(B^{(r)} - A,\ B^{(r-1)} - B^{(r)}) \in (-\pi, 0) \cup \{\pi\}$, *i.e., it cannot be* $\angle\left(B^{(r)} - A,\ B^{(r-1)} - B^{(r)}\right) = 0$.
2. *The following inequalities hold:*

$$0 < \bar{\theta}_r = \angle\left(B^{(r+1)} - B^{(r)}, B^{(r)} - A\right) \le \alpha_r = \angle\left(B^{(r+1)} - B^{(r)}, B^{(r)} - B^{(r-1)}\right) \le \varphi. \quad (1.49)$$

***Proof.*** 1) Note that by Statement 1.34

$$\theta_{r-1} = \angle\left(B^{(r)} - B^{(r-1)}, B^{(r)} - A\right) \in [-\sigma_{r-1}, \sigma_{r-1}]. \quad (1.50)$$

Assume, that $\epsilon = \angle\left(B^{(r)} - A, B^{(r-1)} - B^{(r)}\right) = 0$. Then ( $B^{(r)} \ne A$ by Statement 1.28, as $r\varphi < \pi$) $\epsilon \oplus \theta_{r-1} = \angle\left(B^{(r)} - B^{(r-1)}, B^{(r)} - B^{(r)}\right) = \pi$. Hence, using $\epsilon = 0$, we get $\theta_{r-1} = \pi$, which, under $r\varphi < \pi$, contradicts (1.50).

2) By virtue of Statement 1.9 ($B^{(r)} \ne A$, as discussed before)

$$\theta_{r-1} \oplus \angle\left(B^{(r)} - A,\ B^{(r-1)} - B^{(r)}\right) = \angle\left(B^{(r)} - B^{(r-1)},\ B^{(r-1)} - B^{(r)}\right) = \pi.$$

and due to already proofed statement 1 $\angle(B^{(r)}-A,\ B^{(r-1)}-B^{(r)}) \in (-\pi,0)\cup\{\pi\}$. Therefore, if $\angle\left(B^{(r)}-A,\ B^{(r-1)}-B^{(r)}\right)=\pi$ ( $\gamma\oplus\pi=\pi\Rightarrow\gamma=0$), then $\theta_{r-1}=0$, and if $\angle(B^{(r)}-A,\ B^{(r-1)}-B^{(r)})\in(-\pi,0)$, then, according to Statement 1.14,

$$\theta_{r-1}=-\pi-\angle(B^{(r)}-A,\ B^{(r-1)}-B^{(r)})\in(-\pi,0)$$

holds. Thus, using (1.50), we can conclude that

$$\theta_{r-1}\in(-\pi,0]\cap[-\sigma_{r-1},\sigma_{r-1}]=[-\sigma_{r-1},0]. \qquad (1.51)$$

Statement 1.31(3) gives $\bar{\theta}_r=\theta_{r-1}\oplus\alpha_r$. With equation (1.51) and $r\varphi<\pi$, this yields $\theta_{r-1}+\alpha_r\in[-r\varphi,\varphi]\subseteq(-\pi,\pi)$. Hence, by Statement 1.10,

$$\bar{\theta}_r=\theta_{r-1}\oplus\alpha_r=\theta_{r-1}+\alpha_r\leq\alpha_r\leq\varphi.$$

If we assume that $\angle(B^{(r+1)}-B^{(r)},B^{(r)}-A)\leq 0$, then

$$\angle(B^{(r)}-A,B^{(r+1)}-B^{(r)})\geq 0,$$

which contradicts the condition $\angle(B^{(r)}-A,B^{(r+1)}-B^{(r)})\in(-\pi,0)$.

So, (1.49) is proven. □

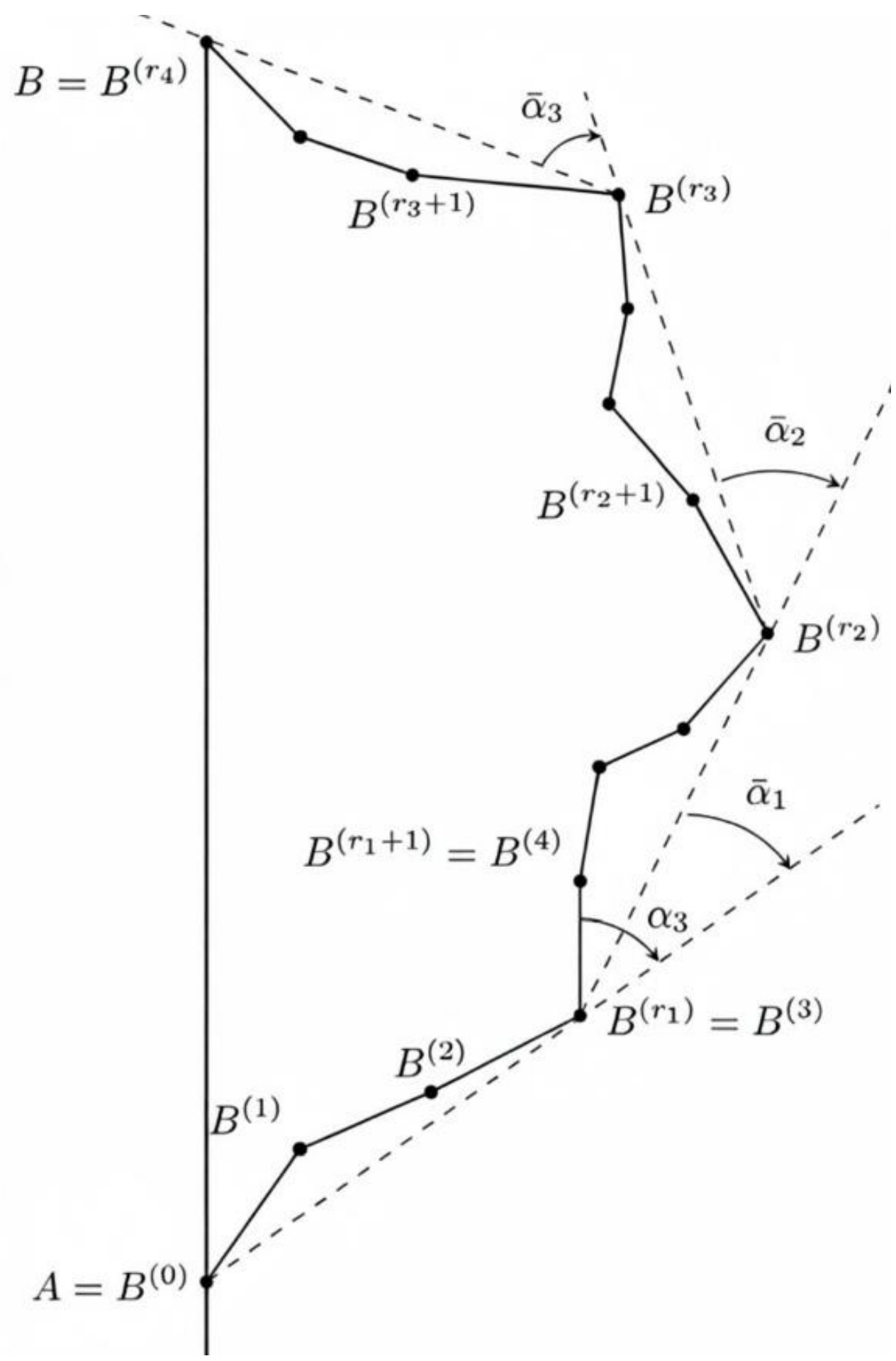

**Fig. 9. External turning points $\boldsymbol{B}^{(r_1)}, \boldsymbol{B}^{(r_2)}, \boldsymbol{B}^{(r_3)}, \boldsymbol{B}^{(r_4)}$**

**Statement 1.36.** *Suppose we are under the conditions of Statement 30, $A \neq B$, all points $B^{(1)}, \dots, B^{(n)}$ lie strictly to the right of the line $\overrightarrow{A, B}$, i.e.,*

$$\angle(B - A, B^{(i)} - A) > 0, i = 1, \dots, n, \qquad (1.52)$$

*and*

$$|\alpha_i| \leq \varphi, i = 1, \dots, n, n\varphi < \pi. \qquad (1.53)$$

*Then*

1) *All points $A, B, B^{(1)}, \dots, B^{(n)}$ are pairwise distinct.*
2) $B^{(1)}, B^{(2)}, \dots, B^{(n)} \in \breve{\boldsymbol{S}}(A, B, n\varphi) \subset \boldsymbol{S}(A, B, n\varphi)$.

*Proof.* Part 1) follows from Statement 1.28 and condition (1.53). Prove 2). Note that by (1.53), for any $i \in \{2, \dots, n\}$ we are under the conditions of Statement 1.33, and therefore inequalities (1.40) hold, by which

$$\gamma_i \in [-i\varphi, i\varphi], i = 1, \dots, n. \qquad (1.54)$$

Now find the **first external turning point** $B^{(r_1)}$ (see Fig. 1.9), $1 \le r_1 \le n$, such that all $B^{(j)}, j = 0,1, \dots, n+1$, lie not strictly to the left of the line $\overrightarrow{A, B^{(r_1)}}$, and all $B^{(j)}, j = r_1 + 1, \dots, n+1$, lie strictly to the left of it (otherwise continue the sequence). By the conditions, $r_1 = n+1$ is impossible because, for example, the point $B^{(1)}$ lies strictly to the right of $\overrightarrow{A, B} = \overrightarrow{A, B^{(n+1)}}$. We use the following conditions for choosing the point $B^{(r_1)}$, where we define $\gamma_0 = \angle(B^{(1)} - A, B^{(1)} - A) = 0$.

**(C1)**: (a) first find

$$\min\{\gamma_{j-1} = \angle(B^{(j)} - A, B^{(1)} - A) \mid j \in \{1, \dots, n+1\}\} = \omega_1$$

($\omega_1 \le 0$ since $\gamma_0 = \angle(B^{(1)} - A, B^{(1)} - A) = 0$);

(b) then find

$$r_1 = \max\{j \mid \gamma_{j-1} = \angle(B^{(j)} - A, B^{(1)} - A) = \omega_1, j \in \{1, \dots, n+1\}\}.$$

By the conditions, $\gamma_n = \angle(B^{(n+1)} - A, B^{(1)} - A) = \angle(B - A, B^{(1)} - A) > 0$, so $1 \le r_1 \le n$. Moreover, since $n\varphi < \pi$, (1.54), and $\gamma_0 = 0$, we have $\omega_1 \in (-\pi, 0]$. We will show that the found point $B^{(r_1)}$, which satisfies **(C1)**, meets all the conditions defined for the **first external turning point**. Indeed, by Statement 1.17 (set in Statement 1.17 $C = B^{(1)}, E = B^{(r_1)}, D = B^{(j)}, j \in \{1, \dots, n+1\}$, where by condition **(C1a)** $\beta = \angle\left(B^{(j)} - A, B^{(1)} - A\right) \ge \omega_1 = \angle\left(B^{(r_1)} - A, B^{(1)} - A\right)$ for $j = 1, \dots, n+1$, and by condition **(C1b)** $\beta = \angle(B^{(j)} - A, B^{(1)} - A) > \omega_1 = \angle(B^{(r_1)} - A, B^{(1)} - A)$ for $j = r_1 + 1, \dots, n+1$), all points $B^{(j)}, j \in \{1, \dots, n+1\}$ (for point $A = B^{(0)}$ it is obvious) lie not strictly to the left of the line $\overrightarrow{A, B^{(r_1)}}$, and all $B^{(j)}, j = r_1 + 1, \dots, n+1$, lie strictly to the left of it. Thus, all conditions defined for the first external turning point are satisfied.

We also show that for $r = r_1$ the conditions of Statement 1.35 are satisfied. Indeed, by (1.53) we have $r\varphi = r_1\varphi \le n\varphi < \pi$. Moreover, the point $B^{(r-1)} = B^{(r_1-1)}$ lies not strictly to the left of the line $\overrightarrow{A, B^{(r_1)}} = \overrightarrow{A, B^{(r)}}$ and the point $B^{(r+1)} = B^{(r_1+1)}$ lies strictly to the left of the line $\overrightarrow{A, B^{(r_1)}} = \overrightarrow{A, B^{(r)}}$, i.e., all conditions of Statement 1.35 are satisfied, by which

$$0 < \angle\left(B^{(r_1+1)} - B^{(r_1)}, B^{(r_1)} - A\right) \le \alpha_{r_1} =$$
$$\angle(B^{(r_1+1)} - B^{(r_1)}, B^{(r_1)} - B^{(r_1-1)}) \le \varphi. \quad (1.55)$$

After finding the **first external turning point** $B^{(r_1)}$, where $1 \le r_1 \le n$, we find the **second external turning point** $B^{(r_2)}$, where $r_1 + 1 \le r_2 \le n+1$, such that all $B^{(j)}, j = 0,1, \dots, n+1$, lie not strictly to the left of the line $\overrightarrow{B^{(r_1)}, B^{(r_2)}}$, and all $B^{(j)}, j = r_2 + 1, \dots, n+1$, lie strictly to the left of it (otherwise continue the sequence; if $r_2 = n+1$, then points of the second group are absent). We use the following conditions for choosing the point $B^{(r_2)}$:

**(C2)**: (a) first find

$$\min\{\angle(B^{(j)} - B^{(r_1)}, B^{(r_1+1)} - B^{(r_1)}) \mid j \in \{0,1,\dots,n+1\} \setminus \{r_1\}\} = \omega_2$$

( $\omega_2 \le 0$ since $\angle(B^{(r_1+1)} - B^{(r_1)}, B^{(r_1+1)} - B^{(r_1)}) = 0$ );

(b) then find

$$r_2 = \max\{j \mid \angle(B^{(j)} - B^{(r_1)}, B^{(r_1+1)} - B^{(r_1)}) = \omega_2, j \in \{0,1,\dots,n+1\} \setminus \{r_1\}\}.$$

Note that when choosing the point $B^{(r_2)}$ in accordance with condition **(C2)**, it suffices to consider only points $B^{(j)}, j \in \{0,1,\dots,n+1\} \setminus \{r_1\}$, that lie not strictly to the right of the line $\overrightarrow{B^{(r_1)}, B^{(r_1+1)}}$ (indeed, for points $B^{(j)}$ lying strictly to the left of the line $\overrightarrow{B^{(r_1)}, B^{(r_1+1)}}$ we have $\angle(B^{(r_1+1)} - B^{(r_1)}, B^{(j)} -$ $B^{(r_1)}) \in (-\pi, 0)$, hence $\angle(B^{(j)} - B^{(r_1)}, B^{(r_1+1)} - B^{(r_1)}) \in (0, \pi)$, i.e., $\angle(B^{(j)} - B^{(r_1)}, B^{(r_1+1)} - B^{(r_1)}) > 0 \geq \omega_2)$. In particular, a point lying strictly to the left of the line $\overrightarrow{B^{(r_1)}, B^{(r_1+1)}}$ is point $A$, since point $B^{(r_1+1)}$ lies strictly to the left of the line $\overrightarrow{A, B^{(r_1)}}$, and consequently, by Statement 1.20, point $A$ lies strictly to the left of the line $\overrightarrow{B^{(r_1)}, B^{(r_1+1)}}$. Accordingly, for points $B^{(j)}, j \in \{1, \dots, n + 1\} \setminus \{r_1\}$, lying not strictly to the right of the line $\overrightarrow{B^{(r_1)}, B^{(r_1+1)}}$ taking into account Statement 1.18 (set in Statement 1.18 $B = B^{(r_1)}, C = B^{(r_1+1)}, D = B^{(j)}, j \in \{1, \dots, n + 1\}$, where by (1.55) $\alpha = \angle(C - B, B - A) = \angle(B^{(r_1+1)} - B^{(r_1)}, B^{(r_1)} - A) \in (0, \varphi]$, and then by Statement 1.18 $\beta = \angle(B^{(r_1+1)} - B^{(r_1)}, B^{(j)} - B^{(r_1)}) = \angle(C - B, D - B) \in [0, \alpha], \gamma = \angle(D - B, B - A) = \angle(B^{(j)} - B^{(r_1)}, B^{(r_1)} - A) \in [0, \alpha] \subseteq [0, \varphi]$; for convenience we use the notation $\alpha, \beta, \gamma$ below) the following holds:

$$\alpha = \angle(C - B, B - A) = \angle(B^{(r_1+1)} - B^{(r_1)}, B^{(r_1)} - A) \in (0, \varphi],$$
$$\angle(B^{(r_1+1)} - B^{(r_1)}, B^{(j)} - B^{(r_1)}) = \angle(C - B, D - B) = \beta \in [0, \alpha],$$
$$\angle(B^{(j)} - B^{(r_1)}, B^{(r_1)} - A) = \gamma \in [0, \alpha] \subseteq [0, \varphi]. \qquad (1.56)$$

Hence

$$0 \geq -\beta = \angle\left(B^{(j)} - B^{(r_1)}, B^{(r_1+1)} - B^{(r_1)}\right) \geq -\alpha =$$
$$\angle(B^{(r_1)} - A, B^{(r_1+1)} - B^{(r_1)}) \geq -\varphi. \qquad (1.57)$$

Moreover, by Statement 1.9 ( $B^{(j)} \neq B^{(r_1)}$ by part 1) )

$$\angle(B^{(j)} - B^{(r_1)}, B^{(r_1)} - B^{(r_1-1)}) = \angle(B^{(j)} - B^{(r_1)}, B^{(r_1+1)} - B^{(r_1)}) \oplus$$
$$\angle(B^{(r_1+1)} - B^{(r_1)}, B^{(r_1)} - B^{(r_1-1)}) =$$
$$\angle(B^{(j)} - B^{(r_1)}, B^{(r_1+1)} - B^{(r_1)}) \oplus \alpha_{r_1}. \quad (1.58)$$

But then from (1.57), (1.58), using Statement 1.10 and taking into account that by (1.55) $\alpha_{r_1} = \angle(B^{(r_1+1)} - B^{(r_1)}, B^{(r_1)} - B^{(r_1-1)}) \in (0, \varphi]$, we obtain

$$\angle(B^{(j)} - B^{(r_1)}, B^{(r_1)} - B^{(r_1-1)}) =$$
$$\angle(B^{(j)} - B^{(r_1)}, B^{(r_1+1)} - B^{(r_1)}) + \alpha_{r_1} \in (-\varphi, \varphi]. \qquad (1.59)$$

Thus, using Statement 1.18, it is shown that (1.59) holds for all $j \in \{0,1, \dots, n+1\} \setminus \{r_1\}$ such that the point $B^{(j)}$ lies not strictly to the right of the line $\overrightarrow{B^{(r_1)}, B^{(r_1+1)}}$ (all these points, as established earlier, simultaneously lie not strictly to the left of the line $\overrightarrow{A, B^{(r_1)}}$ ).

Note that due to the condition $r_1\varphi \le n\varphi < \pi$, there are no points $B^{(j)}, j \in \{0,1, \dots, r_1 - 1\}$, lying not strictly to the right of the line $\overrightarrow{B^{(r_1)}, B^{(r_1+1)}}$. Indeed, as shown, in this case (1.59) holds, and then we could take this point $B^{(j)}$ as the point $B^{(r_1+1)}$ in the sequence of points $B^{(0)}, B^{(1)}, \dots, B^{(r_1+1)}$ satisfying Statement 1.28, i.e., set $B^{(r_1+1)} = B^{(j)}$, since by (1.59) in this case we would have

$$\angle(B^{(r_1+1)} - B^{(r_1)}, B^{(r_1)} - B^{(r_1-1)}) \in (-\varphi, \varphi],$$

which contradicts Statement 1.28 (by Statement 1.28, the equality $B^{(r_1+1)} = B^{(j)}$ is impossible for $j \in \{0,1, \dots, r_1 - 1\}$).

Thus, all points $B^{(j)}, j \in \{0,1, \dots, r_1 - 1\}$, can only lie strictly to the left of the line $\overrightarrow{B^{(r_1)}, B^{(r_1+1)}}$ and, as already shown earlier, for these points we have $\angle(B^{(j)} - B^{(r_1)}, B^{(r_1+1)} - B^{(r_1)}) > 0 \ge \omega_2$, and consequently, in **(C2)** we can everywhere replace the condition $j \in \{0,1, \dots, n+1\} \setminus \{r_1\}$ with $j \in \{r_1 + 1, \dots, n+1\}$. With this change, we have $r_2 \ge r_1 + 1$.

We will show that the found point $B^{(r_2)}$, which satisfies **(C2)**, meets the conditions defined for the **second external turning point**.

To do this, we use Statement 1.17. It was previously established that all points $B^{(j)}, j = 0,1, \dots, n+1$, lie not strictly to the left of the line $\overrightarrow{A, B^{(r_1)}}$, where $1 \le r_1 \le n$, and all points $B^{(j)}, j = r_1 + 1, \dots, n+1$, lie strictly to the left of it. Thus, the points $B^{(r_1+1)}, B^{(r_2)}$, where $r_1 + 1 \leq r_2 \leq n + 1$, lie strictly to the left of the line $\overrightarrow{A, B^{(r_1)}}$. Reversing the direction of this line, we obtain that, accordingly, all points $B^{(j)}, j = 0,1, \dots, n + 1$, lie not strictly to the right of the line $\overrightarrow{B^{(r_1)}, A}$ and the points $B^{(r_1+1)}, B^{(r_2)}$ lie strictly to the right of it. We use Statement 1.17 (set in Statement 1.17 $A = B^{(r_1)}, B = A, C = B^{(r_1+1)}, E = B^{(r_2)}, D = B^{(j)}, j \in \{1, \dots, n + 1\}$, where by condition **(C2a)** $\beta = \angle(B^{(j)} - B^{(r_1)}, B^{(r_1+1)} - B^{(r_1)}) \geq \omega_2 = \angle(B^{(r_2)} - B^{(r_1)}, B^{(r_1+1)} - B^{(r_1)}), j \in \{1, \dots, n + 1\}$, and by condition **(C2b)** in the case $r_2 \leq n$ we have $\beta = \angle(B^{(j)} - B^{(r_1)}, B^{(r_1+1)} - B^{(r_1)}) > \omega_2 = \angle(B^{(r_2)} - B^{(r_1)}, B^{(r_1+1)} - B^{(r_1)}),\ j = r_2 + 1, \dots, n + 1$), from which it follows that all points $B^{(j)}, j = 1, \dots, n + 1$, lie not strictly to the left of the line $\overrightarrow{B^{(r_1)}, B^{(r_2)}}$, and in the case $r_2 \leq n$, all points $B^{(j)}, j = r_2 + 1, \dots, n + 1$, lie strictly to the left of the line $\overrightarrow{B^{(r_1)}, B^{(r_2)}}$. Moreover, since the point $B^{(r_2)}$ lies strictly to the right of $\overrightarrow{B^{(r_1)}, A}$ (because $r_2 \geq r_1 + 1$), the point $A = B^{(0)}$ lies strictly to the left of the line $\overrightarrow{B^{(r_1)}, B^{(r_2)}}$ (see Statement 1.20).

Note also that by (1.56) (recall that condition (1.56) holds for all $B^{(j)}$ lying not strictly to the right of the line $\overrightarrow{B^{(r_1)}, B^{(r_1+1)}}$, and hence for $B^{(r_2)}$ as well, since $\angle(B^{(r_2)} - B^{(r_1)}, B^{(r_1+1)} - B^{(r_1)}) = \omega_2 \leq 0$, hence $\angle(B^{(r_1+1)} - B^{(r_1)}, B^{(r_2)} - B^{(r_1)}) = -\omega_2 \geq 0$; see Statement 1.3)

$$0 < \angle(B^{(r_2)} - B^{(r_1)}, B^{(r_1)} - A) \leq \varphi. \qquad (1.60)$$

The strict inequality on the left follows from the fact that the point $B^{(r_2)}$ lies strictly to the left of the line $\overrightarrow{A, B^{(r_1)}}$, whence (see Statements 1.3, 1.6) $\angle(B^{(r_1)} - A, B^{(r_2)} - B^{(r_1)}) \in (-\pi, 0)$, and consequently $\angle(B^{(r_2)} - B^{(r_1)}, B^{(r_1)} - A) \in (0, \pi)$.

If $r_2 = n + 1$, we conclude that we have only two external turning points: $B^{(r_1)}, B^{(r_2)} = B$, and the process of finding them ends. Otherwise (if $r_2 \leq n$), we can again use Statement 1.35. When using Statement 1.35, to avoid confusion in notation, we add a bar above the corresponding letters. Now the role of point $\bar{A} = \bar{B}^{(0)}$ will be played by point $B^{(r_1)}$ with the corresponding shift of the sequence members, i.e., the role of $\bar{B}^{(1)}$ will be played by $B^{(r_1+1)}$, etc., and accordingly $\bar{B}^{(r)} = B^{(r_2)}$, where $r = r_2 - r_1$, and since $1 \leq r_1, r_1 + 1 \leq r_2 \leq n$, we have $r\varphi = (r_2 - r_1)\varphi \leq (n-1)\varphi < \pi$. Moreover, since the point $\bar{B}^{(r-1)} = B^{(r_2-1)}$ lies not strictly to the left of the line $\overrightarrow{\bar{A}, \bar{B}^{(r)}} = \overrightarrow{B^{(r_1)}, B^{(r_2)}}$, and the point $B^{(r_2+1)}$ lies strictly to the left of this line, all conditions of Statement 1.35 are satisfied, by which, similarly to (1.49),

$$0 < \angle\left(\bar{B}^{(r+1)} - \bar{B}^{(r)}, \bar{B}^{(r)} - \bar{A}\right) \leq$$
$$\angle(\bar{B}^{(r+1)} - \bar{B}^{(r)}, \bar{B}^{(r)} - \bar{B}^{(r-1)}) \leq \varphi.$$

Then, given that $\bar{B}^{(i)} = B^{(i+r_1)}, i = 0, \dots, n - r_1 + 1, r = r_2 - r_1$, we obtain

$$0 < \angle(B^{(r_2+1)} - B^{(r_2)}, B^{(r_2)} - B^{(r_1)}) \leq \alpha_{r_2} =$$
$$\angle(B^{(r_2+1)} - B^{(r_2)}, B^{(r_2)} - B^{(r_2-1)}) \leq \varphi. \qquad (1.61)$$

Using (1.61) and Statement 1.18, we obtain a condition analogous to (1.56) that holds for all $B^{(j)}$ lying not strictly to the right of the line $\overrightarrow{B^{(r_2)}, B^{(r_2+1)}}$:

$$\angle(B^{(j)} - B^{(r_2)}, B^{(r_2)} - B^{(r_1)}) \in [0, \varphi]. \qquad (1.62)$$

If $r_2 < n + 1$, we similarly look for the **third external turning point** $B^{(r_3)}$, where $r_2 + 1 \leq r_3 \leq n + 1$, such that all $B^{(j)}, j = 0,1, \dots, n + 1$, lie not strictly to the left of the line $\overrightarrow{B^{(r_2)}, B^{(r_3)}}$, and in the case $r_3 \leq n$, all $B^{(j)}, j = r_3 +$ $1, \dots, n+1$, lie strictly to the left of it (otherwise continue the sequence). We use the following conditions for choosing the point $B^{(r_3)}$:

**(C3)**: (a) first find

$\min\{\angle(B^{(j)} - B^{(r_2)}, B^{(r_2+1)} - B^{(r_2)}) \mid j \in \{0,1, \dots, n+1\} \setminus \{r_2\}\} = \omega_3$

( $\omega_3 \le 0$ since $\angle(B^{(r_2+1)} - B^{(r_2)}, B^{(r_2+1)} - B^{(r_2)}) = 0$ );

(b) then find

$r_3 = \max\{j \mid \angle(B^{(j)} - B^{(r_2)}, B^{(r_2+1)} - B^{(r_2)}) = \omega_3, j \in \{0,1, \dots, n+1\} \setminus \{r_2\}\}.$

Using Statement 1.28 and (1.61), (1.62), we conclude that in **(C3)** (as in the case of condition **(C2)**) we can everywhere replace the condition $j \in \{0,1, \dots, n+1\} \setminus \{r_2\}$ with $j \in \{r_2+1, \dots, n+1\}$. With this change, we have $r_3 \ge r_2 + 1$.

We will show that the found point $B^{(r_3)}$, which satisfies **(C3)**, meets the conditions defined for the **third external turning point**. To do this, we again use Statement 1.17. It was previously established that all points $B^{(j)}, j = 0,1, \dots, n+1$, lie not strictly to the left of the line $\overrightarrow{B^{(r_1)}, B^{(r_2)}}$, where $r_1 + 1 \le r_2 \le n$, and all points $B^{(j)}, j = r_2+1, \dots, n+1$, lie strictly to the left of it. Thus, the points $B^{(r_2+1)}, B^{(r_3)}$, where $r_2 + 1 \le r_3 \le n+1$, lie strictly to the left of the line $\overrightarrow{B^{(r_1)}, B^{(r_2)}}$. Reversing the direction of this line, we obtain that, accordingly, all points $B^{(j)}, j = 0,1, \dots, n+1$, lie not strictly to the right of the line $\overrightarrow{B^{(r_2)}, B^{(r_1)}}$, and the points $B^{(r_2+1)}, B^{(r_3)}$ lie strictly to the right of it. We use Statement 1.17 (set in the Statement 1.17 $A = B^{(r_2)}, B = B^{(r_1)}, C = B^{(r_2+1)}, E = B^{(r_3)}, D = B^{(j)}, j \in \{0,1, \dots, n+1\}$, where by condition **(C3a)** $\beta = \angle(B^{(j)} - B^{(r_2)}, B^{(r_2+1)} - B^{(r_2)}) \ge \omega_3 = \angle(B^{(r_3)} - B^{(r_2)}, B^{(r_2+1)} - B^{(r_2)}), j \in \{0,1, \dots, n+1\}$, and by condition **(C3b)** in the case $r_3 \le n$ we have $\beta = \angle(B^{(j)} - B^{(r_2)}, B^{(r_2+1)} - B^{(r_2)}) > \omega_3 = \angle(B^{(r_3)} - B^{(r_2)}, B^{(r_2+1)} - B^{(r_2)}), j =$ $r_3+1,\dots,n+1$), from which it follows that all points $B^{(j)}, j=0,1,\dots,n+1$, lie not strictly to the left of the line $\overrightarrow{B^{(r_2)},B^{(r_3)}}$, and in the case $r_3 \le n$, all points $B^{(j)}, j=r_3+1,\dots,n+1$, lie strictly to the left of the line $\overrightarrow{B^{(r_2)},B^{(r_3)}}$.

If $r_3 = n+1$, we conclude that we have only three external turning points: $B^{(r_1)}, B^{(r_2)}, B^{(r_3)} = B$, and the process of finding them ends. Otherwise (if $r_3 \le n$), we can again use Statement 1.35. When using Statement 1.35, to avoid confusion in notation, we add a bar above the corresponding letters. Now the role of point $\bar{A} = \bar{B}^{(0)}$ will be played by point $B^{(r_2)}$ with the corresponding shift of the sequence members, i.e., the role of $\bar{B}^{(1)}$ will be played by $B^{(r_2+1)}$, etc., and accordingly $\bar{B}^{(r)} = B^{(r_3)}$. Now $r = r_3 - r_2$, and since $2 \le r_2, r_2+1 \le r_3 \le n$, we have $r\varphi = (r_3-r_2)\varphi \le (n-2)\varphi < \pi$. Moreover, since the point $\bar{B}^{(r-1)} = B^{(r_3-1)}$ lies not strictly to the left of the line $\overrightarrow{\bar{A},\bar{B}^{(r)}} = \overrightarrow{B^{(r_2)},B^{(r_3)}}$, and the point $B^{(r_3+1)}$ lies strictly to the left of this line, all conditions of Statement 1.35 are satisfied, by which, similarly to (1.49), (1.61),

$$0 < \angle(B^{(r_3+1)} - B^{(r_3)}, B^{(r_3)} - B^{(r_2)}) \le \alpha_{r_3} =$$
$$\angle(B^{(r_3+1)} - B^{(r_3)}, B^{(r_3)} - B^{(r_3-1)}) \le \varphi. \tag{1.63}$$

Note also that by (1.62) (recall that condition (1.62) holds for all $B^{(j)}$ lying not strictly to the right of the line $\overrightarrow{B^{(r_2)},B^{(r_2+1)}}$, and hence for $B^{(r_3)}$) we have

$$0 < \angle(B^{(r_3)} - B^{(r_2)}, B^{(r_2)} - B^{(r_1)}) \le \varphi. \tag{1.64}$$

The strict inequality on the left follows from the fact that since $r_3 \ge r_2+1$, the point $B^{(r_3)}$ lies strictly to the left of the line $\overrightarrow{B^{(r_1)},B^{(r_2)}}$ (see the definition of the **second external turning point** $B^{(r_2)}$).

We proceed in this way until we reach, for some $k \ge 1$, the equality $r_{k+1} = n+1$. As a result, we obtain a sequence of points

$$A = B^{(0)}, B^{(r_1)}, B^{(r_2)}, \dots, B^{(r_{k+1})} = B^{(n+1)},$$

such that, similarly to (1.60), (1.64) (set $r_0 = 0$),

$$0 < \bar{\alpha}_i = \angle(B^{(r_{i+1})} - B^{(r_i)}, B^{(r_i)} - B^{(r_{i-1})}) \le \varphi, i = 1, \dots, k, k \le n,$$

all points $A = B^{(0)}, B^{(1)}, B^{(2)}, \dots, B^{(n+1)}$ lie not strictly to the left of the lines

$$\overrightarrow{B^{(r_{i-1})}, B^{(r_i)}}, i = 1, \dots, k+1.$$

Then, by Statement 1.10 and the condition $k\varphi \le n\varphi < \pi$, we have

$$\angle(B - B^{(r_k)}, B^{(r_1)} - A) = \bar{\alpha}_1 \oplus \cdots \oplus \bar{\alpha}_k = \bar{\alpha}_1 + \cdots + \bar{\alpha}_k \in (0, n\varphi] \subset (0, \pi),$$

and moreover

$$\angle(B - B^{(r_k)}, B^{(r_1)} - A) \oplus \angle(B^{(r_1)} - A, B^{(r_k)} - B) = \pi.$$

Then by Statement 1.13 we have

$$\angle\big(B^{(r_1)} - A, B^{(r_k)} - B\big) = \pi - \angle\big(B - B^{(r_k)}, B^{(r_1)} - A\big) \in [\pi - n\varphi, \pi) \quad (1.65)$$

and accordingly

$$\angle(B^{(r_k)} - B, B^{(r_1)} - A) = -\angle(B^{(r_1)} - A, B^{(r_k)} - B) \in (-\pi, n\varphi - \pi] \subset (-\pi, 0).$$

Moreover, the points $B^{(r_1)}, B^{(r_k)}$ lie strictly to the right of the line $\overrightarrow{A, B}$, and consequently we are under the conditions of Statement 1.23, by which the lines $\overrightarrow{A, B^{(r_1)}}$ and $\overrightarrow{B, B^{(r_k)}}$ intersect at a unique point $C \in \mathbb{R}^2$ lying strictly to the right of the line $\overrightarrow{A, B}$. Further, since all points $A = B^{(0)}, B^{(1)}, B^{(2)}, \dots, B^{(n+1)}$ lie not strictly to the right of the line $\overrightarrow{A, B}$, not strictly to the right of the line $\overrightarrow{B, B^{(r_k)}} = \overrightarrow{B, C}$, and not strictly to the left of the line $\overrightarrow{A, B^{(r_1)}} = \overrightarrow{A, C}$ (or not strictly to the right of $\overrightarrow{C, A}$), then, as shown in statement 1.38 (see below), all these points belong to the triangle $\triangle(A, B, C)$ with angle (see (1.65))

$$\angle(A - C, B - C) = \angle\big(A - B^{(r_1)}, B - B^{(r_k)}\big) = \angle\big(B^{(r_1)} - A, B^{(r_k)} - B\big)$$
$$\in [\pi - n\varphi, \pi) \subset (0, \pi).$$

Moreover

$$\angle(A - C, B - C) \oplus \angle(B - C, C - A) = \angle(A - C, C - A) = \pi,$$

hence, using Statement 1.13, we obtain

$$\angle(B-C, C-A) = \pi - \angle(A-C, B-C) \in (0, n\varphi],$$

i.e.,

$$C \in \breve{\boldsymbol{S}}(A, B, n\varphi).$$

But then, by convexity of $\breve{\boldsymbol{S}}(A, B, n\varphi)$, all points of $\triangle(A, B, C) = \mathrm{Co}\,\{A, B, C\}$ (see below Statement 1.38) (except $A, B$) belong to $\breve{\boldsymbol{S}}(A, B, n\varphi)$, whence

$$B^{(1)}, B^{(2)}, \dots, B^{(n)} \in \breve{\boldsymbol{S}}(A, B, n\varphi) \subset \boldsymbol{S}(A, B, n\varphi). \ \square$$

We will need the following auxiliary statement.

***Statement 1.37.*** *Let* $E^{(i)}, C \in \mathbb{R}^n$*,* $i = 1, \dots, r, r \geq n$*, and consider the system of linear inequalities (where the vectors* $E^{(i)} \in \mathbb{R}^n$*,* $i = 1, \dots, r$*, are given, and the vector* $C$ *is a possible solution to this system):*

$$\langle E^{(i)}, C \rangle \leq t_i, i = 1, \dots, r. \tag{1.66}$$

*Then for the solution set of this system to be bounded, it is* ***sufficient*** *that there exist numbers* $\alpha_i \geq 0, i = 1, \dots, r$*, such that*

$$\textstyle\sum_{i=1}^{r} \alpha_i E^{(i)} = (0, \dots, 0) \in \mathbb{R}^n, \tag{1.67}$$

*and among the vectors* $E^{(i)}$ *with* $\alpha_i > 0, n$ *vectors are linearly independent.*

***Proof.*** For simplicity of notation, suppose $\alpha_i > 0, i = 1, \dots, n$, and the vectors $E^{(i)}, i = 1, \dots, n$, are linearly independent. Define:

$$\bar{t}_i = \textstyle\sum_{j=1}^{r} \frac{\alpha_j}{\alpha_i} t_j - t_i, i = 1, \dots, n.$$

From (1.66) and (1.67) we obtain:

$$\langle -E^{(i)}, C \rangle = \textstyle\sum_{j=1, j \neq i}^{r} \frac{\alpha_j}{\alpha_i} \langle E^{(i)}, C \rangle \leq \sum_{j=1, j \neq i}^{r} \frac{\alpha_j}{\alpha_i} t_j = \bar{t}_i, i = 1, \dots, n. \tag{1.68}$$

We will show, for example, that the solutions of system (1.66) are bounded with respect to the first variable. From the linear independence of $E^{(i)}, i = 1, \dots, n$, it follows that there exist numbers $\beta_i, \gamma_i, i = 1, \dots, n$, such that

$$(1,0, \dots, 0) = \textstyle\sum_{i=1}^{n} \beta_i E^{(i)}, (-1,0, \dots, 0) = \sum_{i=1}^{n} \gamma_i E^{(i)}. \tag{1.69}$$

Denote

$$\bar{\beta}_i = \begin{cases} \beta_i t_i, \text{if } \beta_i \geq 0, \\ -\beta_i \bar{t}_i, \text{ if } \beta_i < 0, \end{cases} \quad \bar{\gamma}_i = \begin{cases} \gamma_i t_i, \text{ if } \gamma_i \geq 0, \\ -\gamma_i \bar{t}_i, \text{ if } \gamma_i < 0. \end{cases}$$

Then, as shown, for any solution $C = (c_1, \dots, c_n)$ of system (1.66), (1.68) holds; consequently, by (1.69) we have:

$$c_1 = \langle (1,0,\dots,0), C\rangle = \sum_{i=1}^{n} \beta_i \langle E^{(i)}, C\rangle \leq \sum_{i=1}^{n} \bar{\beta}_i,$$

$$-c_1 = \langle (-1,0,\dots,0), C\rangle = \sum_{i=1}^{n} \gamma_i \langle E^{(i)}, C\rangle \leq \sum_{i=1}^{n} \bar{\gamma}_i,$$

or

$$-\sum_{i=1}^{n} \bar{\gamma}_i \leq c_1 \leq \sum_{i=1}^{n} \bar{\beta}_i, \ |c_1| \leq \max\{\left|\sum_{i=1}^{n} \bar{\beta}_i\right|, \left|\sum_{i=1}^{n} \bar{\gamma}_i\right|\}. \ \square$$

Now we can prove that the following holds:

***Statement 1.38.*** *Let $A, B, C \in \mathbb{R}^2$, $A \neq B$, and point $C \in \mathbb{R}^2$ lie strictly to the right of the line $\overrightarrow{A,B}$, i.e., $\begin{vmatrix} B-A \\ C-A \end{vmatrix} < 0$. Let $\boldsymbol{M}$ be the set of all points in $\mathbb{R}^2$ lying: (a) not strictly to the right of the line $\overrightarrow{A,B}$, (b) not strictly to the right of the line $\overrightarrow{B,C}$, (c) not strictly to the right of the line $\overrightarrow{C,A}$. Then:*

1. *Point $A$ lies strictly to the right of the line $\overrightarrow{B,C}$, and point $B$ lies strictly to the right of the line $\overrightarrow{C,A}$; hence $A, B, C \in \boldsymbol{M}$.*
2. $\boldsymbol{M} = \left\{D \in \mathbb{R}^2 \mid \begin{vmatrix} B-A \\ D-A \end{vmatrix} \leq 0, \begin{vmatrix} C-B \\ D-B \end{vmatrix} \leq 0, \begin{vmatrix} A-C \\ D-C \end{vmatrix} \leq 0\right\}$,

   *i.e., $\boldsymbol{M}$ is a nonempty convex polyhedral set (since it is defined by a finite system of linear inequalities).*
3. *Any pair of vectors: $B-A, C-B$ or $C-B, A-C$ or $A-C, B-A$ is linearly independent, and moreover*

   $$(B-A) + (C-B) + (A-C) = (0,0) \in \mathbb{R}^2. \quad (1.70)$$
4. *$\boldsymbol{M}$ is a bounded set, i.e., it is a polytope.*
5. *$\{A, B, C\}$ is the set of all extreme points of the polytope $\boldsymbol{M}$; hence $\boldsymbol{M} =$ Co $\{A, B, C\} = \triangle(A, B, C)$.*

***Proof.*** Prove 1). Since point $C$ lies strictly to the right of the line $\overrightarrow{A,B}$, by Statement 1.20 we obtain that point $B$ lies strictly to the left of the line $\overrightarrow{A,C}$, and point $A$ lies strictly to the right of the line $\overrightarrow{B,C}$. Moreover, it is obvious that point $B$ lies strictly to the right of the line $\overrightarrow{C,A}$ (see Statements 1.3, 1.12). Thus, each of the points $A, B, C$ lies strictly to the right of one of the lines $\overrightarrow{A,B}$, $\overrightarrow{B,C}$, $\overrightarrow{C,A}$ and lies on the other two lines, which implies their membership in $\boldsymbol{M}$.

Prove 2). For example, for a point $D \in \mathbb{R}^2$ to lie not strictly to the right of the line $\overrightarrow{A,B}$, it is necessary and sufficient that it either belongs to this line, which is equivalent to $\begin{vmatrix} B - A \\ D - A \end{vmatrix} = 0$ (see Statement 1.2), or lies strictly to the right of the line $\overrightarrow{A,B}$, which is equivalent to $\angle(B - A, D - A) \in (0, \pi)$ (see Statement 1.3) and $\begin{vmatrix} B - A \\ D - A \end{vmatrix} < 0$.

3) Prove the linear independence of the vectors $B - A, C - B$. Suppose $\exists t \neq 0$: $C - B = t(B - A)$. Then

$$0 > \begin{vmatrix} B - A \\ C - A \end{vmatrix} = \begin{vmatrix} B - A \\ C - B + B - A \end{vmatrix} = \begin{vmatrix} B - A \\ C - B \end{vmatrix} = t \begin{vmatrix} B - A \\ B - A \end{vmatrix} = 0,$$

a contradiction. Next, considering the vectors $C - B, A - C$, suppose $\exists t \neq 0$: $A - C = t(C - B)$. Then by Statement 1.20 we see that point A lies strictly to the right of the line $\overrightarrow{B,C}$; consequently,

$$0 > \begin{vmatrix} C - B \\ A - B \end{vmatrix} = \begin{vmatrix} C - B \\ A - C + C - B \end{vmatrix} = \begin{vmatrix} C - B \\ A - C \end{vmatrix} = t \begin{vmatrix} C - B \\ C - B \end{vmatrix} = 0,$$

a contradiction. The linear independence of the vectors $A - C, B - A$ is proved similarly. Equality (1.70) is obvious.

4) Follows from the proven Statement 3), equality (1.70), and Statement 1.37.

5) Prove that $A, B, C$ are extreme points of the polytope $\boldsymbol{M}$. Denote $l_1(D) = \begin{vmatrix} B-A \\ D-A \end{vmatrix}, l_2(D) = \begin{vmatrix} C-B \\ D-B \end{vmatrix}, l_3(D) = \begin{vmatrix} A-C \\ D-C \end{vmatrix}$.

If for any point $D \in \mathbb{R}^2$ we denote $I(D) = \{i \in \{1,2,3\} | l_i(D) = 0\}$, then, for example, $I(A) = \{1, 2\}$. As shown in, e.g., [8], a point $D \in \boldsymbol{M}$ is an extreme point of the polytope $\boldsymbol{M}$ iff the rank of the system of linear nonhomogeneous equations $l_i(D) = 0, i \in I(D)$, equals two (i.e., the dimension of $\boldsymbol{M}$). It is perfectly obvious that all three points $A, B, C$ satisfy this condition, and only they do; hence $\boldsymbol{M} = \mathrm{Co}\,\{A, B, C\} = \triangle\,(A, B, C)$ (see, e.g., [8]). □

## 2. Proof of Theorem 1

***Proof.*** Let us consider two cases.

**Case** $n = 1$**.** Condition (0.4) gives $\angle\big(B - B^{(1)}, B^{(1)} - A\big) \in [-\varphi, \varphi]$, which by definition means that $B^{(1)} \in \mathrm{cl}\,\mathbf{S}(A, B, \varphi) \backslash \{A, B\}$ (taking into account the condition $B^{(1)} \neq A, B$). The statement is proved.

**Case** $n \geq 2$**.** Subcases are possible:

**(a)** $\boldsymbol{L} \cap \{\overrightarrow{A, B}\} = \emptyset$.

**(b)** $\boldsymbol{L} \cap \{\overrightarrow{A, B}\} \neq \emptyset$, where $\{\overrightarrow{A, B}\} = \{A + t(B - A) | t \in \mathbb{R}\}$.

Consider case **(a)** ($\boldsymbol{L} \cap \{\overrightarrow{A, B}\} = \emptyset$). Then, by Property 1 (see introduction), we can assume without loss of generality that all angular points $B^{(1)}, \ldots, B^{(n)}$ of the polygonal line $\boldsymbol{L}$ lie strictly to the right of the line $\overrightarrow{A, B}$, and then by Statement 1.36 we have $B^{(1)},\ldots,B^{(n)} \in \breve{\boldsymbol{S}}(A, B, n\varphi) \subset \boldsymbol{S}(A, B, n\varphi)$. Indeed, in the case where all points $B^{(1)}, \ldots, B^{(n)}$ satisfying (0.4) lie strictly to the left of the line $\overrightarrow{A, B}$, we proceed to consider the set $\boldsymbol{S}(B, A, n\varphi)$ $(= \boldsymbol{S}(A, B, n\varphi))$ and moreover all these points lie strictly to the right of the directed line $\overrightarrow{B,A}$. Accordingly, instead of the polyline $\boldsymbol{L}$ satisfying (0.3), we can now consider the polyline

$$\bar{\boldsymbol{L}} = \left[\bar{B}^{(0)}, \bar{B}^{(1)}\right] \cup \left[\bar{B}^{(1)}, \bar{B}^{(2)}\right] \cup \cdots \cup \left[\bar{B}^{(n)}, \bar{B}^{(n+1)}\right],$$

where $\bar{B}^{(i)} = B^{(n+1-i)}$, $i = 0,1,\dots,n+1$, which satisfies properties analogous to (0.4), i.e.

$$\angle\left(\bar{B}^{(i+1)} - \bar{B}^{(i)}, \bar{B}^{(i)} - \bar{B}^{(i-1)}\right) \in [-\varphi,\varphi], i = 1,\dots,n,$$
$$\bar{B}^{(i)} \neq \bar{B}^{(i+1)}, i = 0,1,\dots,n.$$

We show that in case **(b)** ($\boldsymbol{L} \cap \left\{\overrightarrow{A,B}\right\} \neq \emptyset$), for any variant of its fulfillment, either we arrive at a contradiction with the condition $n\varphi < \pi$, or we can reduce to consideration of subproblems with a smaller number of angular points than $n$, which makes it possible to apply the proof by mathematical induction on $n$ – the number of angular points.

In this regard, we immediately discard the cases with $B^{(1)} \in \left\{\overrightarrow{A,B}\right\}$. Indeed, the case $B^{(1)} \in (A,B)$ is easily reduced by Property 2 (see introduction) to the case with the number of angular points one less, and in the cases $B^{(1)} \in (A(-\infty),A] \cup [B,A(+\infty))$ we arrive at a contradiction with Statement 1.28 or Remark 1.4 (see Section 1.2). In this regard, note that by (0.4)

$$B^{(i+1)} \in B^{(i)} + \boldsymbol{C}\left(B^{(i)} - B^{(i-1)}, \varphi\right), i = 1,\dots,n, \quad (2.1)$$

and moreover $n\varphi < \pi$, i.e., we are under the conditions of Statements 1.27,1.28 and Remark 1.4. But then, if, for example, $B^{(1)} \in [B,A(+\infty))$ (the case $B^{(1)} \in (A(-\infty),A]$ is similar), then when $B^{(1)} = B = B^{(n+1)}$ we arrive at a contradiction with Statement 1.28, and if $B^{(1)} \in (B,A(+\infty))$, then for some $t > 0$ we have $B^{(1)} = B + t(B-A)$, whence for $\bar{t} = \frac{t}{1+t} \in (0,1)$ we have $B = B^{(n+1)} = B^{(1)} + \bar{t}(A - B^{(1)})$, which contradicts Remark 1.4.

Thus, it remains to consider the cases with $B^{(1)} \notin \{\overrightarrow{A,B}\}$. In this case, without loss of generality, we can restrict ourselves to the case where the point $B^{(1)}$ is located strictly to the right of the line $\overrightarrow{A,B}$ (the opposite case is considered similarly). Let $C \in \{\overrightarrow{A,B}\} \cap \boldsymbol{L}$, where $\boldsymbol{L}$ satisfies (0.3), and let $C$ be the first point of such an intersection, i.e. (taking into account that $B^{(1)} \notin \{\overrightarrow{A,B}\}$), for some $k \in \{1, \dots, n-1\}$ we have:

$$[A, B^{(1)}] \cap \{\overrightarrow{A,B}\} = \{A\}; \text{ for } k \geq 2\ [B^{(i-1)}, B^{(i)}] \cap \{\overrightarrow{A,B}\} = \emptyset, i = 2, \dots, k,$$
$$B^{(1)} \notin \{\overrightarrow{A,B}\}, C \in [B^{(k)}, B^{(k+1)}] \cap \{\overrightarrow{A,B}\} \neq \emptyset.$$

The following cases are possible:

1) $B^{(k+1)} \in \{\overrightarrow{A,B}\}$,

2) $B^{(k+1)} \notin \{\overrightarrow{A,B}\}$.

In case 1) the following subcases are possible:

1a) $B^{(k+1)} \in \{A, B\}$,

1b) $B^{(k+1)} \in (A(-\infty), A)$,

1c) $B^{(k+1)} \in (B, A(+\infty))$,

1d) $B^{(k+1)} \in (A, B)$.

In case 1a), if $B^{(k+1)} = A = B^{(0)}$, then we arrive at a contradiction with Statement 1.28, since $k + 1 \leq n$, whence $\varphi(k+1) \leq \varphi n < \pi$ and (2.1) holds. Accordingly, in the case $B^{(k+1)} = B = B^{(n+1)}$ we also arrive at a contradiction with Statement 1.28, because for $k = n - 1$ we obtain $B^{(n)} = B^{(n+1)}$, which contradicts the condition $B^{(n)} \neq B$, and for $k \leq n - 2$ we have: if $i = n, j = k + 1$, then $i + 1 = n + 1 \geq j + 2 = k + 3$, and moreover $B^{(i+1)} = B^{(n+1)} = B = B^{(k+1)} = B^{(j)}$, which also contradicts Statement 28.

In case 1b) ($B^{(k+1)} \in (A(-\infty), A)$) denote $H = B^{(k+1)}$. Let $t \geq 0$ be a parameter. Consider the function $\Lambda_t: \mathbb{R}^2 \to \mathbb{R}^2$

$$\forall C \in \mathbb{R}^2 \;\; \Lambda_t(C) = C + t(C - H) = (1 + t)C - tH.$$

Clearly, $\Lambda_t(H) = H + t(H - H) = H, \Lambda_t(C) - H = (1 + t)(C - H)$, $\Lambda_t(C) - C = t(C - H),\; C \neq H \Rightarrow |\Lambda_t(C) - H| = (1 + t)|C - H| > |C - H|$.

It is easy to show that $\forall t \geq 0$ the function $\Lambda_t$ satisfies the following properties: $\forall D, E, F, G \in \mathbb{R}^2$, if $\bar{D} = \Lambda_t(D), \bar{E} = \Lambda_t(E), \bar{F} = \Lambda_t(F), \bar{G} = \Lambda_t(G)$, then

$$\bar{D} - \bar{E} = (1 + t)(D - E), \angle(D - E, F - G) = \angle(D - E, \bar{F} - \bar{G}) =$$
$$= \angle(\bar{D} - \bar{E}, F - G) = \angle(\bar{D} - \bar{E}, \bar{F} - \bar{G}). \qquad (2.2)$$

Choose $\bar{t} > 0: \; \Lambda_{\bar{t}}(A) = A + \bar{t}(A - H) = B$ (the number $\bar{t} > 0$ exists because $H = B^{(k+1)} \in (A(-\infty), A)$), whence

$$B - A = \Lambda_{\bar{t}}(A) - A = \bar{t}(A - H), \bar{t} = \frac{|B-A|}{|A-H|} > 0.$$

Consider the sequence:

$$\begin{gathered} \bar{B}^{(0)} = \Lambda_{\bar{t}}(A) = B, \bar{B}^{(1)} = \Lambda_{\bar{t}}\big(B^{(1)}\big), \dots, \bar{B}^{(k)} = \Lambda_{\bar{t}}\big(B^{(k)}\big), \\ \bar{B}^{(k+1)} = \Lambda_{\bar{t}}\big(B^{(k+1)}\big) = \Lambda_{\bar{t}}(H) = H = B^{(k+1)}, \qquad (2.3) \\ \bar{B}^{(k+2)} = B^{(k+2)}, \dots, \bar{B}^{(n+1)} = B^{(n+1)} = B. \end{gathered}$$

In the obtained sequence, by virtue of (2.2) and (2.3), we have

$$\angle\big(\bar{B}^{(2)} - \bar{B}^{(1)}, \bar{B}^{(1)} - \bar{B}^{(0)}\big) = \angle\big(B^{(2)} - B^{(1)}, B^{(1)} - B^{(0)}\big), \dots,$$
$$\angle\big(\bar{B}^{(k+1)} - \bar{B}^{(k)}, \bar{B}^{(k)} - \bar{B}^{(k-1)}\big) = \angle\big(B^{(k+1)} - B^{(k)}, B^{(k)} - B^{(k-1)}\big),$$
$$\angle\big(\bar{B}^{(k+2)} - \bar{B}^{(k+1)}, \bar{B}^{(k+1)} - \bar{B}^{(k)}\big) = \angle\big(B^{(k+2)} - B^{(k+1)}, B^{(k+1)} - B^{(k)}\big), \dots,$$
$$\angle\big(\bar{B}^{(n+1)} - \bar{B}^{(n)}, \bar{B}^{(n)} - \bar{B}^{(n-1)}\big) = \angle\big(B^{(n+1)} - B^{(n)}, B^{(n)} - B^{(n-1)}\big),$$

and thus, taking into account (2.1), we obtain

$$\bar{B}^{(i+1)} \in \bar{B}^{(i)} + \boldsymbol{C}\big(\bar{B}^{(i)} - \bar{B}^{(i-1)}, \varphi\big), i = 1, \dots, n, \qquad (2.4)$$

and moreover $\bar{B}^{(0)} = B = \bar{B}^{(n+1)}$, which contradicts Statement 1.28.

In case 1c) ($B^{(k+1)} \in (B, A(+\infty))$) we act similarly to case 1b), also set $H = B^{(k+1)}$ and use the same function $\Lambda_t(C)$, but now choose $\bar{t} > 0$: $\Lambda_{\bar{t}}(B) = B + \bar{t}(B - H) = A$, whence

$$A - B = \Lambda_{\bar{t}}(B) - B = \bar{t}(B - H), \bar{t} = \frac{|B-A|}{|B-H|} > 0.$$

Furthermore, instead of (2.3) we use the sequence

$$\bar{B}^{(0)} = \Lambda_{\bar{t}}(A) = B, \bar{B}^{(1)} = \Lambda_{\bar{t}}\left(B^{(1)}\right), \dots, \bar{B}^{(k)} = \Lambda_{\bar{t}}\left(B^{(k)}\right),$$
$$\bar{B}^{(k+1)} = \Lambda_{\bar{t}}\left(B^{(k+1)}\right) = \Lambda_{\bar{t}}(H) = H = B^{(k+1)},$$
$$\bar{B}^{(k+2)} = B^{(k+2)}, \dots, \bar{B}^{(n+1)} = B^{(n+1)} = B.$$

Case 1d) ($B^{(k+1)} \in (A, B)$) reduces to the consideration of two subproblems with a smaller number of angular points, which makes it possible (using Property 2) to prove by mathematical induction on $n$ – the number of angular points.

Now consider case 2) ($B^{(k+1)} \notin \{\overrightarrow{A, B}\}$). In this case, the point $B^{(k+1)}$, where $k \in \{1, \dots, n-1\}$, is located strictly to the left of the line $\overrightarrow{A, B}$ and the point $C \in \left[B^{(k)}, B^{(k+1)}\right] \cap \{\overrightarrow{A, B}\}$ is the unique point of the set $\left[B^{(k)}, B^{(k+1)}\right] \cap \{\overrightarrow{A, B}\}$.

From the condition $B^{(k)}, B^{(k+1)} \notin \{\overrightarrow{A, B}\}$ we obtain that $C \in \left(B^{(k)}, B^{(k+1)}\right)$, whence it follows that for some $t \in (0,1)$ we have:

$$C = tB^{(k)} + (1-t)B^{(k+1)} = B^{(k+1)} + t(B^{(k)} - B^{(k+1)}). \qquad (2.5)$$

Regarding the point $C$, the following subcases are also possible:

2a) $C \in \{A, B\}$,

2b) $C \in (A(-\infty), A)$,

2c) $C \in (B, A(+\infty))$,

2d) $C \in (A, B)$.

In case 2a), if $C = A = B^{(0)}$, then considering the sequence of points:

$$B^{(0)}, B^{(1)}, \dots, B^{(k)}, C(= B^{(0)})$$

and taking into account the fulfillment of (2.1), we arrive at a contradiction with Statement 28, since $k \le n-1, k\varphi \le n\varphi < \pi$,

$\angle\big(C - B^{(k)}, B^{(k)} - B^{(k-1)}\big) = \angle\big(B^{(k+1)} - B^{(k)}, B^{(k)} - B^{(k-1)}\big) \in [-\varphi, \varphi]$,

whence $C \in B^{(k)} + \boldsymbol{C}\big(B^{(k)} - B^{(k-1)}, \varphi\big)$ (i.e., the point $C = B^{(0)}$ can be taken as the next $(k$+1)-th term in the sequence $B^{(0)}, B^{(1)}, \dots, B^{(k)}, B^{(k+1)}(= C)$ satisfying the conditions of Statement 1.28). Accordingly, in the case $C = B = B^{(n+1)}$, by virtue of (2.5), we arrive at a contradiction with Remark 1.4.

Case 2b) ($C \in (A(-\infty), A)$ is considered similarly to case 1b), but in this case we use $H = C$. In this case, the following minor differences arise. We choose $\bar{t} > 0$ from the condition: $\Lambda_{\bar{t}}(A) = A + \bar{t}(A - C) = B$, i.e., now $\bar{t} =| B - A |/| A - C |> 0$. Next, we consider the sequence:

$$\bar{B}^{(0)} = \Lambda_{\bar{t}}(A) = B, \bar{B}^{(1)} = \Lambda_{\bar{t}}\big(B^{(1)}\big), \dots, , \bar{B}^{(k)} = \Lambda_{\bar{t}}\big(B^{(k)}\big),$$

$$\bar{B}^{(k+1)} = B^{(k+1)}, \dots, \bar{B}^{(n+1)} = B^{(n+1)} = B.$$

Here, taking into account (2.5), we have

$$\bar{B}^{(k)} = \Lambda_{\bar{t}}\big(B^{(k)}\big) = B^{(k)} + \bar{t}\big(B^{(k)} - C\big) = B^{(k)} + \bar{t}(1-t)\big(B^{(k)} - B^{(k+1)}\big),$$

whence

$$\bar{B}^{(k+1)} - \bar{B}^{(k)} = B^{(k+1)} - B^{(k)} - \bar{t}(1-t)\big(B^{(k)} - B^{(k+1)}\big) =$$

$$\big(B^{(k+1)} - B^{(k)}\big)(1 + \bar{t}(1-t)),$$

$$\angle\big(\bar{B}^{(k+1)} - \bar{B}^{(k)}, \bar{B}^{(k)} - \bar{B}^{(k-1)}\big) = \angle\big(B^{(k+1)} - B^{(k)}, B^{(k)} - B^{(k-1)}\big)$$

and thus, taking into account (2.1), (2.2), we obtain

$$\bar{B}^{(i+1)} \in \bar{B}^{(i)} + \boldsymbol{C}\big(\bar{B}^{(i)} - \bar{B}^{(i-1)}, \varphi\big), i = 1, \dots, n. \quad (2.4)$$

Then, similarly to Case 1b), we arrive at a contradiction with Statement 1.28.

Case 2c) ($C \in (B, A(+\infty))$) is considered similarly to case 1c), and in this case we also use $H = C$.

Case 2d) ($C \in (A, B)$) is analogous to case 1d) and reduces to the consideration of two subproblems with a smaller number of angular points, which makes it possible (using Property 2) to prove by mathematical induction on $n$ – the number of angular points. □

## 3. Proof of Lemma 1

We first present some auxiliary statements.

**Statement 2.1.** *Let* $A, B \in \mathbb{R}^2$, $A \neq B, \varphi \in (0, \pi), C \in \boldsymbol{S}(A, B, \varphi), C \in \overrightarrow{A, B}$. *Then* $C \in (A, B)$.

***Proof.*** By the definition of $\boldsymbol{S}(B, B, \varphi)$, from $C \in \boldsymbol{S}(A, B, \varphi)$ we obtain $C \neq A, C \neq B$. Assume that $C \in \overrightarrow{A, B} \setminus [A, B]$. Then either $C = A + t(A - B)$ or $C = B + t(B - A)$ with $t > 0$. In the case $C = A + t(A - B)$ we have:

$$\angle(B - C, C - A) = \angle\big(B - A - t(A - B), t(A - B)\big) =$$
$$= \angle\big((1 + t)(B - A), t(A - B)\big) = \angle(B - A, A - B) = \pi > \varphi.$$

which contradicts the condition $C \in \boldsymbol{S}(A, B, \varphi)$. In the case $C = B + t(B - A)$ we have:

$$\angle(B - C, C - A) = \angle(t(A - B), (1 + t)(B - A)) =$$
$$= \angle(A - B, B - A) = \pi > \varphi.$$

which contradicts the condition $C \in \boldsymbol{S}(A, B, \varphi)$. □

Consider also the following auxiliary problem.

**Problem O (on partitioning a general rotation into intermediate rotations).** Let pairwise distinct points $A, B, G \in \mathbb{R}^2$ be given, along with the angle $\Psi = \angle(B - G,\ G - A) \in (0, \pi)$, and a point $C = A + t(G - A)$ with $t \in$ $(0,1)$. It is required to construct $n \geq 2$ points $B^{(1)}(=C), B^{(2)}, \dots, B^{(n)}$ such that (denoting also $B^{(0)} = A, B^{(n+1)} = B$):

1. The angles between consecutive segments are equal to $\alpha = \Psi/n$:
$$\angle(B^{(2)} - B^{(1)}, B^{(1)} - A) = \alpha, \angle(B^{(k)} - B^{(k-1)}, B^{(k-1)} - B^{(k-2)}) = \alpha,$$
$$k = 2, \dots, n, \angle(B - B^{(n)}, B^{(n)} - B^{(n-1)}) = \alpha.$$
2. The lengths of the segments are equal:
$$| B^{(2)} - B^{(1)} |=| B^{(3)} - B^{(2)} |= \cdots =| B^{(n)} - B^{(n-1)} |= d.$$
3. The point $B^{(n)}$ lies on the segment $[B, G]$, and $| B^{(n)} - G |=| C - G |$.

   We introduce the following notation:
$$U = G - A, l_A =| U |, \bar{U} = U/l_A,$$
$$V = G - B, l_B =| V |, \bar{V} = V/l_B,$$

$R(\theta) = \begin{bmatrix} \cos\theta & -\sin\theta \\ \sin\theta & \cos\theta \end{bmatrix}$ – the counterclockwise rotation matrix by angle $\theta$.

**Statement 2.2.** Let $| B - G |>| C - G |$. Then

1. The points $B^{(2)}, \dots, B^{(n)}$ satisfying conditions 1)–3) exist and are given by the formulas:
$$B^{(k)} = A + tU + \lambda \sum_{j=1}^{k-1} R(j\alpha)\bar{U} \ , k = 2, \dots, n, \qquad (3.1)$$

where
$$\alpha = \frac{\Psi}{n}, \lambda = \frac{2(1-t)\cos(\Psi/2)\sin(\alpha/2)}{\sin((\Psi-\alpha)/2)}. \qquad (3.2)$$

In this case, the segment length is $d = \lambda l_A$, and the position of the point $B^{(n)}$ on the line $\overrightarrow{BG}$ is: $B^{(n)} = G + sV$, where $s = (1-t)\mu$, $\mu = l_A/l_B$; i.e.,
$$B^{(n)} = B + (1-s)(G-B), B^{(n)} - G = sV = (1-t)\mu V,$$
$$| B^{(n)} - G |= (1-t)\mu \, | V |= (1-t) \, | U |=| C - G |$$

(since $| C - G |=| A + t(G-A) - G |= (1-t) \, | A - G |= (1-t) \, | U |$).

2. The points $B^{(1)}, B^{(2)}, \dots, B^{(n)}$ belong to the triangle $\triangle ABG$.

3. For any $n \geq 1$ we have:

$$\lambda \leq \frac{2(1-t)}{n}, d = \lambda l_A \leq \frac{2(1-t)}{n} \mid G - A \mid.$$

The expression

$$\frac{2\cos{(\Psi/2)}\sin{(\Psi/2)}}{\sin{((\Psi-\alpha)/2)}}, \alpha = \frac{\Psi}{n}, \Psi \in (0, \pi)$$

has, for fixed $n \geq 2$, the supremum $2/n$, which is attained as the limit $\Psi \to 0+$. Inside the interval, the maximum is strictly less than $2/n$; i.e., there is no exact maximum in the domain, only a supremum.

We now present the ***proof*** of Lemma 1. If $n = 1, B^{(1)} \in \mathrm{cl}\, \boldsymbol{S}(A, B, \varphi)$, then the segments $[A, B^{(1)}], [B^{(1)}, A]$ form the desired polygonal line, and condition (0.4) holds by definition of the closure of the set $\boldsymbol{S}(A, B, \varphi)$.

Now let $n \geq 2$ and $B^{(1)} \in \boldsymbol{S}(A, B, n\varphi), n\varphi < \pi$. It is required to construct points $B^{(1)}, B^{(2)}, \dots, B^{(n)}$ satisfying (0.4).

**Case 1.** Consider first the simplest case when $B^{(1)} \in \overrightarrow{A, B}$. Then, by Statement 2.1, we have $B^{(1)} \in (A, B)$. A sequence of two rotations, first by angle $\varphi$ and then by $-\varphi$, leads to motion parallel to $\overrightarrow{A, B}$; the next rotation by $-\varphi$ may bring us back to the line $\overrightarrow{A, B}$. If there are exactly three rotations, by adjusting the lengths of the segments between rotations we can exactly reach point $B$ after the third rotation. If there are only two rotation points, we make the first by $\varphi/2$ and the second by $-\varphi$. Again, by adjusting the lengths we can exactly reach point $B$ after the second rotation. If there are four rotation points, proceeding as in the three-rotation case, we make the third turning point lie on $(A, B)$. Then after a rotation by $\varphi$ we reach point $B$ along the line $\overrightarrow{A, B}$. All other possible cases reduce to combinations of those described. For example, in the case of five rotations we finish the motion using the two-rotation technique, i.e., we combine the three-rotation and two-rotation approaches.

**Case 2.** Let $B^{(1)} \notin \overrightarrow{A,B}$. Without loss of generality, assume that point $B^{(1)}$ lies strictly to the right of the line $\overrightarrow{A,B}$ (the opposite case is symmetric), i.e., $\angle(B-A,\, B^{(1)}-A) > 0$, whence

$$0 > \begin{vmatrix} B-A \\ B^{(1)}-A \end{vmatrix} = \begin{vmatrix} B-A-(B^{(1)}-A) \\ B^{(1)}-A \end{vmatrix} = \begin{vmatrix} B-B^{(1)} \\ B^{(1)}-A \end{vmatrix},$$

which implies $\angle(B-B^{(1)},\, B^{(1)}-A) > 0$, and since $B^{(1)} \in \boldsymbol{S}(A,B,n\varphi)$, we have $\angle(B-B^{(1)},\, B^{(1)}-A) \in (0,n\varphi) \subset (0,\pi)$. Consider

$$G(t) = B^{(1)} + t(B^{(1)}-A), \Psi(t) = \angle(B-G(t),\, G(t)-A), t \geq 0.$$

Let

$$s = \min\{|\, B^{(1)}-B \,|\,,\ |\, B^{(1)}-A \,|\} > 0. \qquad (3.3)$$

By (3.3) and because $G(0) = B^{(1)}$, the function $\Psi(t)$ is continuous for sufficiently small $t \geq 0$ (see [2] and Section 5.3 bellow); consequently, for some $\bar{t} > 0$ we have

$$\left|G(\bar{t}) - B^{(1)}\right| < \frac{s}{2}, \bar{\Psi} = \Psi(\bar{t}) \in (0,n\varphi) \subset (0,\pi). \qquad (3.4)$$

Set $\bar{G} = G(\bar{t})$. Then using (3.3) and (3.4) we obtain

$$|\, \bar{G} - B \,| > \frac{s}{2} > |\, \bar{G} - B^{(1)} \,|. \qquad (3.5)$$

Note that for $t_1 = \frac{1}{1+\bar{t}}$ we have $B^{(1)} = A + t_1(\bar{G} - A)$, because

$$A + t_1(\bar{G} - A) = A + t_1(B^{(1)} + \bar{t}(B^{(1)} - A) - A)$$
$$= (1 - t_1(1+\bar{t}))A + t_1(1+\bar{t})B^{(1)} = B^{(1)}.$$

Thus we are in the conditions of **Problem O**, and by (3.5) Statement 2.2 applies, from which (setting in that statement $C = B^{(1)}$, $G = \bar{G}$, $t = t_1 =$ $\frac{1}{1+\bar{t}}$, $\Psi = \bar{\Psi} \in (0, n\varphi) \subset (0, \pi)$, $\alpha = \frac{\bar{\Psi}}{n} \in (0, \varphi)$) it follows that there exist (and they can be constructed explicitly according to simple formulas) points $B^{(1)}, B^{(2)}, \dots, B^{(n)}$ such that

$$\angle\left(B^{(i+1)} - B^{(i)}, B^{(i)} - B^{(i-1)}\right) \in \alpha = \frac{\bar{\Psi}}{n} \in (0, \varphi), i = 1, \dots, n,$$

where $B^{(0)} = A, B^{(n+1)} = B$; whence in particular it follows that for the polyline $\boldsymbol{L}$ composed of these points according to (0.3), conditions (0.4) and (0.5) also hold. Moreover, the statement of Lemma 1 has been significantly strengthened. In particular, it turned out that in the case $B^{(1)} \notin (A, B)$, the required points can be sought in the same part (i.e., in the set $\breve{\boldsymbol{S}}(A, B, n\varphi)$ or $\widehat{\boldsymbol{S}}(A, B, n\varphi)$) to which the point $B^{(1)}$ belongs. Indeed, as follows from the case considered in the given proof, all the found points $B^{(1)}, B^{(2)}, \dots, B^{(n)}$ l?ie strictly to the right of the line $\overrightarrow{A, B}$. □

## 4. Proof of Theorem 2

***Proof.*** Let $\left(B^{(1)}, \dots, B^{(n)}\right) \in \boldsymbol{S}^{(n)}(A, B, \varphi)$, i.e. $B^{(1)}, \dots, B^{(n)}$ – is a sequence of interior points of a polygonal chain satisfying (0.3), (0.4). Applying Theorem 1 to the polygonal chain $\boldsymbol{L}$, we obtain $B^{(1)} \in \boldsymbol{S}(A, B, n\varphi)$. Next, applying the same theorem to the part of the polygonal chain starting at $B^{(1)}$, we have $B^{(2)} \in \boldsymbol{S}\left(B^{(1)}, B, (n-1)\varphi\right)$. Furthermore, from condition (0.4) for $i = 1$ it follows that $\angle\left(B^{(2)} - B^{(1)}, B^{(1)} - A\right) \in [-\varphi, \varphi]$, i.e. $B^{(2)} \in B^{(1)} + \boldsymbol{C}\left(B^{(1)} - A, \varphi\right)$. Thus,

$$B^{(2)} \in \boldsymbol{S}\left(B^{(1)}, B, (n-1)\varphi\right) \cap \left[B^{(1)} + \boldsymbol{C}\left(B^{(1)} - A, \varphi\right)\right].$$

Continuing analogously, for each $i = 2, \dots, n-2$ we obtain

$$B^{(i+1)} \in \boldsymbol{S}\big(B^{(i)}, B, (n-i)\varphi\big) \cap \big[B^{(i)} + \boldsymbol{C}\big(B^{(i)} - B^{(i-1)}, \varphi\big)\big],$$

and for the last point

$$B^{(n)} \in \big[\mathrm{cl}\, \boldsymbol{S}\big(B^{(n-1)}, B, \varphi\big) \backslash \{A, B\}\big] \cap \big[B^{(n-1)} + \boldsymbol{C}\big(B^{(n-1)} - B^{(n-2)}, \varphi\big)\big].$$

Hence, the sequence $\big(B^{(1)}, \dots, B^{(n)}\big)$ satisfies the membership condition on the right-hand side of (0.6). □

**Converse inclusion.** We prove that for any sequential selection points $B^{(1)}, \dots, B^{(n)}$ according to the right-hand side of (0.6), each term in this direct product will be non-empty and any such sequence satisfies (0.4). Indeed, suppose we have chosen any point

$$B^{(1)} \in \boldsymbol{S}(A, B, n\varphi). \qquad (4.1)$$

By Lemma 1, for the point $B^{(1)}$ there exists a polygonal chain

$$\boldsymbol{L} = \big[B^{(0)}, B^{(1)}\big] \cup \big[B^{(1)}, B^{(2)}\big] \cup \dots \cup \big[B^{(n)}, B^{(n+1)}\big],$$

satisfying (0.3), (0.4). Then, by (0.4), according to the definition of the set $\boldsymbol{C}(E, \varphi)$, for the point $B^{(2)}$ from any such polygonal chain we have

$$B^{(2)} \in B^{(1)} + \boldsymbol{C}\big(B^{(1)} - A, \varphi\big).$$

Note that the point $B^{(2)}$ is the first turning point of the polygonal chain

$$\bar{\boldsymbol{L}} = \big[B^{(1)}, B^{(2)}\big] \cup \dots \cup \big[B^{(n)}, B^{(n+1)}\big],$$

and by Theorem 1, $B^{(2)} \in \boldsymbol{S}\big(B^{(1)}, B, (n-1)\varphi\big)$. Thus,

$$B^{(2)} \in \boldsymbol{S}\big(B^{(1)}, B, (n-1)\varphi\big) \cap \big[B^{(1)} + \boldsymbol{C}\big(B^{(1)} - A, \varphi\big)\big], \qquad (4.2)$$

and the set on the right-hand side of (4.2) is non-empty.

Suppose we have already chosen two points $B^{(1)}, B^{(2)}$, satisfying (4.1), (4.2). By Lemma 1, for the point $B^{(2)}$ there exists a polygonal chain

$$\bar{\boldsymbol{L}} = \big[B^{(1)}, B^{(2)}\big] \cup \dots \cup \big[B^{(n)}, B^{(n+1)}\big],$$

satisfying (0.4). Then, by (0.4), according to the definition of the set $\boldsymbol{C}(E, \varphi)$, for the point $B^{(3)}$ from any such polygonal chain we have

$$B^{(3)} \in B^{(2)} + \boldsymbol{C}\big(B^{(2)} - B^{(1)}, \varphi\big).$$

Note that the point $B^{(3)}$ is the first turning point of the polygonal chain

$$\bar{\bar{\boldsymbol{L}}} = \big[B^{(2)}, B^{(3)}\big] \cup \dots \cup \big[B^{(n)}, B^{(n+1)}\big],$$

and by Theorem 1, $B^{(3)} \in \boldsymbol{S}\big(B^{(2)}, B, (n-2)\varphi\big)$. Thus,

$$B^{(3)} \in \boldsymbol{S}\big(B^{(2)}, B, (n-2)\varphi\big) \cap \big[B^{(2)} + \boldsymbol{C}\big(B^{(2)} - B^{(1)}, \varphi\big)\big], \quad (4.3)$$

and the set on the right-hand side of (4.3) is non-empty. In a completely analogous way, it is proved that each subsequent term in the direct product on the right-hand side of (0.6) is non-empty for any selection of the previous terms according to (0.6). Thus, the sequence of points $B^{(1)}, \dots, B^{(n)}$ in the direct product on the right-hand side of (0.6) can always be extended for any choice of the initial terms $B^{(1)}, \dots, B^{(i)}$ (where $i \in \{0, \dots, n-1\}$). The fulfillment of (0.4) for this sequence follows from the condition

$$B^{(i)} \in B^{(i-1)} + \boldsymbol{C}\big(B^{(i-1)} - B^{(i-2)}, \varphi\big), i = 2, \dots, n+1. \square$$

## 5. On the approximation of the set $\boldsymbol{S}^{(n)}(\boldsymbol{A}, \boldsymbol{B}, \boldsymbol{\varphi})$ by the discrete set $\boldsymbol{S}^{(n)}(\boldsymbol{Q}^{(\tau)}, \boldsymbol{A}, \boldsymbol{B}, \boldsymbol{\varphi})$

### 5.1. On the continuous dependence (in the Hausdorff metric) of the set $\boldsymbol{S}^{(n)}(\boldsymbol{A}, \boldsymbol{B}, \boldsymbol{\varphi})$ on certain constraints on $\big(\boldsymbol{B}^{(1)}, \dots, \boldsymbol{B}^{(n)}\big) \in \boldsymbol{S}^{(n)}(\boldsymbol{A}, \boldsymbol{B}, \boldsymbol{\varphi})$

Let $\varphi \in (0,\ \pi), n \geq 1, A = B^{(0)} \in \mathbb{R}^2, B = B^{(n+1)} \in \mathbb{R}^2, A \neq B, s > 0$. Denote:

$$\boldsymbol{S}^{(n,s)}(A, B, \varphi) = \{\big(B^{(1)}, \dots, B^{(n)}\big) \in \boldsymbol{S}^{(n)}(A, B, \varphi) \,\big|\, \big|B^{(i)} - B^{(j)}\big| \geq s,$$
$$i, j = 0, 1, \dots, n+1, i \neq j\}.$$

**Statement 5.1.** Let $n\varphi < \pi$, and assume that for all $\big(B^{(1)}, \dots, B^{(n)}\big) \in \boldsymbol{S}^{(n)}(A, B, \varphi)$ condition (0.5) holds. Then (using the notation for $h_m,\ \alpha_m,\ \rho_m$ from [1]):

$$\lim_{s\to 0+} h_m\left(\boldsymbol{S}^{(n)}(A,B,\varphi),\boldsymbol{S}^{(n,s)}(A,B,\varphi)\right) =$$

$$= \lim_{s\to 0+} \alpha_m\left(\boldsymbol{S}^{(n)}(A,B,\varphi),\boldsymbol{S}^{(n,s)}(A,B,\varphi)\right) = 0, \qquad (5.1)$$

**Proof.** Assuming that (5.1) does not hold, we obtain that

$$\exists a > 0, W(k) = \left(B^{(1)}(k), \dots, B^{(n)}(k)\right) \in \boldsymbol{S}^{(n)}(A,B,\varphi),\ s_k > 0:$$

$$\rho_m\left(W(k),\boldsymbol{S}^{(n,s_k)}(A,B,\varphi)\right) \geq a,\ k = 1,2,\dots, \text{ and } s_k \to 0 \text{ as } k \to \infty.$$

Due to the boundedness of $\boldsymbol{S}^{(n)}(A,B,\varphi)$ (since $n\varphi < \pi$), for simplicity of notation, we assume that:

$$W(k) \to \bar{W} \in \mathrm{cl}\,\boldsymbol{S}^{(n)}(A,B,\varphi),\ \rho_m\left(\bar{W},\boldsymbol{S}^{(n,s_k)}(A,B,\varphi)\right) \geq \frac{a}{2},$$

$$k = N, N+1, N+2, \dots,$$

where $N \in \mathbb{N}$. Let $\bar{W}(v) = \left(\bar{B}^{(1)}(v), \dots, \bar{B}^{(n)}(v)\right) \in \boldsymbol{S}^{(n)}(A,B,\varphi), v = 1,2,\dots$, $\bar{W}(v) \to \bar{W} \in \mathrm{cl}\,\boldsymbol{S}^{(n)}(A,B,\varphi)$ as $v \to \infty$, and for some index $\bar{v} \in \mathbb{N}$ it holds that $\|\bar{W}(\bar{v}) - \bar{W}\| < \frac{a}{4}$. Since $\bar{W}(\bar{v}) \in \boldsymbol{S}^{(n)}(A,B,\varphi)$, then (see Statement 1 of Theorem 1, as well as the assumption that (0.5) holds for all $\left(B^{(1)}, \dots, B^{(n)}\right) \in \boldsymbol{S}^{(n)}(A,B,\varphi)$):

$$\left|\bar{B}^{(i)}(\bar{v}) - \bar{B}^{(j)}(\bar{v})\right| > 0, \text{ where } j = 0,1,\dots,n+1, i \neq j,$$

$$\bar{\alpha}_i = \angle\left(\bar{B}^{(i+1)}(\bar{v}) - \bar{B}^{(i)}(\bar{v}), \bar{B}^{(i)}(\bar{v}) - \bar{B}^{(i-1)}(\bar{v})\right) \in [-\varphi,\varphi],$$

Consequently, there exists an index $k_{\bar{v}} \geq N$, such that $\bar{W}(\bar{v}) \in \boldsymbol{S}^{\left(n,s_{k_{\bar{v}}}\right)}(A,B,\varphi)$. Then: $\frac{a}{4} > \|\bar{W} - \bar{W}(\bar{v})\| \geq \rho_m\left(\bar{W},\boldsymbol{S}^{\left(n,s_{k_{\bar{v}}}\right)}(A,B,\varphi)\right) \geq \frac{a}{2}$,

i.e., we have reached a contradiction. Thus, (5.1) is proved.

**Remark 5.1.** Statement 5.1 indicates that we can strengthen the conditions on $\left(B^{(1)},\dots,B^{(n)}\right) \in \boldsymbol{S}^{(n)}(A,B,\varphi)$ of the form $B^{(i)} \neq B^{(j)}$, where $i,j \in \{0,1,\dots,n+1\}, i \neq j$, replacing them with stronger ones used in Statement 5.1 and depending on a sufficiently small parameter $s > 0$. In addition, one more "strengthening" will be needed, namely, instead of the condition

$$\angle\left(B^{(i+1)} - B^{(i)}, B^{(i)} - B^{(i-1)}\right) \in [-\varphi, \varphi], i = 1, \dots, n,$$

we use a condition of the form

$$\angle\left(B^{(i+1)} - B^{(i)}, B^{(i)} - B^{(i-1)}\right) \in [-\varphi + \omega, \varphi - \omega], i = 1, \dots, n,$$

where the number $\omega > 0$ is sufficiently small.

We first prove the following simple statement.

**Statement 5.2.** Let $C, D, E$ be pairwise distinct points, $d, t > 0$, $|D - C| = d, |E - D| = t,$

$\varphi = \angle(E - D, D - C)) \in (0, \pi)\ (\in (-\pi, 0)), \delta = \angle(E - C, D - C)),$

$\gamma = \angle(E - D, E - C).$

Then

$$\delta \in (0, \varphi)\ (\in (-\varphi, 0)), \gamma \in (0, \varphi)\ (\in (-\varphi, 0)), \gamma = \varphi - \delta,$$

$$\delta = \operatorname{arctg}\left[\frac{\sin\varphi}{\frac{d}{t} + cos\varphi}\right] > 0\ (< 0). \quad (5.2)$$

***Proof.*** Consider the first case (the second case with negative angles $\varphi, \delta, \gamma$ is treated similarly). From the consideration of $\triangle (C, D, E)$ (see Fig. 1) we obtain $\pi - \varphi + \delta + \gamma = \pi$, whence $\gamma = \varphi - \delta$. Next, note that

$$\text{tg}\ \delta = \frac{t\sin\varphi}{d + t\cos\varphi} = \frac{\sin\varphi}{\frac{d}{t} + cos\varphi} > 0\ (< 0),$$

from which (5.2) follows.

Note also that

$$\varphi = \angle(E - D, D - C) \in (0, \pi) \Leftrightarrow \begin{vmatrix} E - D \\ D - C \end{vmatrix} < 0 \Leftrightarrow$$

$$\Leftrightarrow \begin{vmatrix} E - D \\ D - C \end{vmatrix} = \begin{vmatrix} E - D + D - C \\ D - C \end{vmatrix} = \begin{vmatrix} E - C \\ D - C \end{vmatrix} < 0 \Leftrightarrow$$

$$\Leftrightarrow \delta = \angle(E - C, D - C) \in (0, \pi),$$

$$\varphi = \angle(E - D, D - C) \in (0, \pi) \Leftrightarrow \begin{vmatrix} E - D \\ D - C \end{vmatrix} < 0 \Leftrightarrow$$

$$\Leftrightarrow \begin{vmatrix} E - D \\ D - C \end{vmatrix} = \begin{vmatrix} E - D \\ D - E + E - C \end{vmatrix} = \begin{vmatrix} E - D \\ E - C \end{vmatrix} < 0 \Leftrightarrow$$

$$\Leftrightarrow \gamma = \angle(E - D, E - C) \in (0, \pi).$$

Furthermore,

$$\varphi = \gamma \oplus \delta \Rightarrow \varphi = \gamma + \delta + 2j\pi, j \in \mathbb{Z},$$

$$\gamma, \delta \in (0, \pi) \Rightarrow \gamma + \delta \in (0, 2\pi).$$

Using the fact that $\varphi \in (0, \pi)$, we obtain:

- $j = 0$ works, because $\varphi = \gamma + \delta \in (0,2\pi) \cap (0, \pi) \neq \emptyset$;
- $j = 1$ does not work, because $\varphi = \gamma + \delta + 2\pi \in (2\pi, 4\pi) \cap (0, \pi) = \emptyset$;
- $j = -1$ does not work, because $\varphi = \gamma + \delta - 2\pi \in (-2\pi, 0) \cap (0, \pi) = \emptyset$.

The remaining cases are even less suitable. Thus, only one case remains: $j = 0, \varphi = \gamma + \delta, \gamma = \varphi - \delta$. □

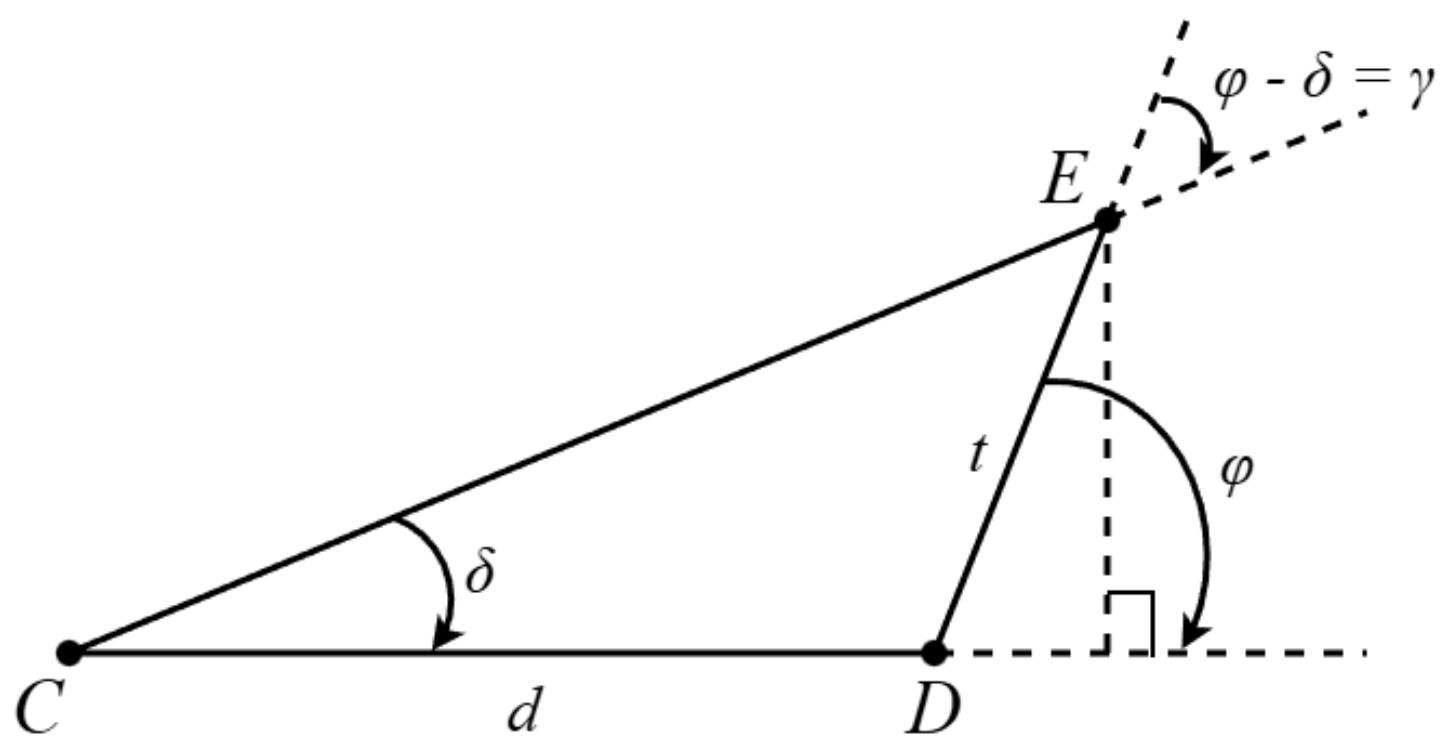

**Fig. 5.1**

For a set $\boldsymbol{M} \subset \mathbb{R}^2$, denote $\operatorname{diam}_e \boldsymbol{M} = \sup\{| C - D | : C, D \in \boldsymbol{M}\}$.

**Statement 5.3.** Let $A, B \in \mathbb{R}^2, A \neq B, \varphi \in (0, \pi/2), n\varphi < \pi, s > 0, d \geq \operatorname{diam}_e \boldsymbol{S}(A, B, n\varphi) \geq s, \kappa = d/s$. Then there exist numbers $\omega = \omega(s, d, \varphi, n) > 0, \theta = \theta(s, d, \varphi, n) > 2$ (which can be constructively computed from the values $s, d, \varphi, n$) such that for every $W = (B^{(1)}, \dots, B^{(n)}) \in \boldsymbol{S}^{(n,s)}(A, B, \varphi)$, for every $t \in (0, s/\theta]$ (i.e., $t > 0$, $s/t \geq \theta$), one can construct (by simple algorithmically realizable formulas) a set of points (depending on the parameter $t$)

$$\tilde{B}^{(i)}(t) \in [B^{(i)}, B^{(i+1)}], i = 1, \dots, n$$

(i.e., the points $\tilde{B}^{(1)}(t), \dots, \tilde{B}^{(n)}(t)$ belong to the polygonal line

$$L = [B^{(0)}, B^{(1)}] \cup [B^{(1)}, B^{(2)}] \cup \cdots \cup [B^{(n)}, B^{(n+1)}]$$

such that

$$\widetilde{W}(t) = (\tilde{B}^{(1)}(t), \dots, \tilde{B}^{(n)}(t)) \in \boldsymbol{S}^{(n,s-t\kappa)}(A, B, \varphi - \omega t),$$

$$\| W - \widetilde{W}(t) \| \leq \frac{t\kappa}{2}. \qquad (5.3)$$

***Proof.*** Let $W = (B^{(1)}, \dots, B^{(n)}) \in \boldsymbol{S}^{(n,s)}(A, B, \varphi)$. Then (by definition of this set)

$$s \leq | B^{(i)} - B^{(i-1)} | \leq \text{diam}_{\,e}\, S(A, B, n\varphi) \leq d, i = 1, \dots, n+1.$$

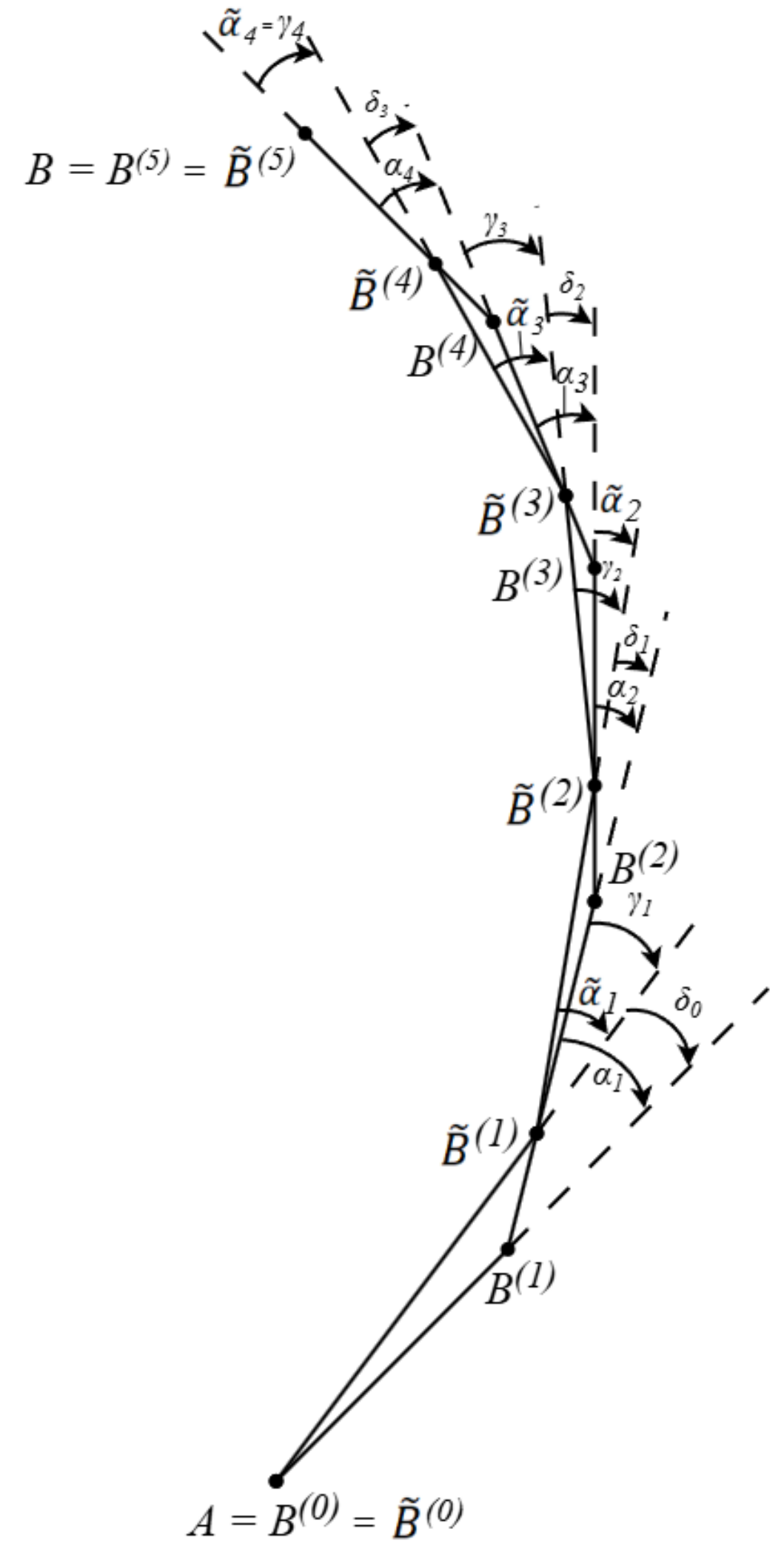

**Fig. 5.2.** The case $n = 4$

Set $A = B^{(0)}, B = B^{(n+1)}$. Consider

$$\alpha_i = \angle(B^{(i+1)} - B^{(i)},\, B^{(i)} - B^{(i-1)}) \in [-\varphi, 0) \cup (0, \varphi], i = 1, \dots, n.$$

First let $t \in (0, s/2)$ (the restrictions on $t$ will be refined later), $\widetilde{W}(t) = (\tilde{B}^{(1)}(t), \dots, \tilde{B}^{(n)}(t))$,

$$\begin{aligned}&\tilde{B}^{(0)}(t) = B^{(0)}, \tilde{B}^{(i)}(t) = B^{(i)} + t\omega_{i-1}\frac{B^{(i+1)}-B^{(i)}}{\left|B^{(i+1)}-B^{(i)}\right|}, i = 1, \dots, n,\\ &\tilde{B}^{(n+1)}(t) = B^{(n+1)},\end{aligned} \tag{5.4}$$

where $\omega_i > 0$, $i = 0, \dots, n-1$, are parameters of the method. Then (using $\varphi \in (0,\ \pi/2) \Rightarrow \left|\tilde{B}^{(i+1)}(t) - \tilde{B}^{(i)}(t)\right| \geq \left|B^{(i+1)} - \tilde{B}^{(i)}(t)\right|$)

$$\left|\tilde{B}^{(i)}(t) - B^{(i)}\right| = t\omega_{i-1}, i = 1, \dots, n,$$

$$\left|\tilde{B}^{(i+1)}(t) - \tilde{B}^{(i)}(t)\right| \geq \left|B^{(i+1)} - \tilde{B}^{(i)}(t)\right| \geq \left|B^{(i+1)} - B^{(i)}\right| -$$

$$-\left|B^{(i)} - \tilde{B}^{(i)}(t)\right| = \left|B^{(i+1)} - B^{(i)}\right| - t\omega_{i-1} \geq$$

$$\geq s - t\omega_{i-1}, i = 1, \dots, n, \qquad (5.5)$$

$$\left|\tilde{B}^{(1)}(t) - \tilde{B}^{(0)}(t)\right| = \left|\tilde{B}^{(1)}(t) - B^{(0)}\right| \geq \left|B^{(1)} - B^{(0)}\right| \geq s. \qquad (5.6)$$

Denote:

$$k_i = \frac{\left|B^{(i)} - B^{(i-1)}\right|}{s} \in [1, \kappa], i = 1, \dots, n+1. \qquad (5.7)$$

Set

$$\omega_0 = 1, \omega_i = \frac{k_{i+1}\,\omega_{i-1}}{\omega_{i-1} + k_i}, i = 1, \dots, n. \qquad (5.8)$$

Then (since $k_2/\omega_1 = 1 + k_1$, $k_3/\omega_2 = 1 + k_2/\omega_1 = 2 + k_1, \dots$)

$$\frac{k_{i+1}}{\omega_i} = \frac{\omega_{i-1} + k_i}{\omega_{i-1}} = 1 + \frac{k_i}{\omega_{i-1}} = i + \frac{k_1}{\omega_0} = i + k_1, i = 0, \dots, n, \quad (5.9)$$

$$\frac{k_{i+1}}{\omega_i} - \frac{k_i}{\omega_{i-1}} = 1, i = 1, \dots, n, \qquad (5.10)$$

$$1 \leq i + 1 \leq \frac{k_{i+1}}{\omega_i} = i + k_1 \leq n + \kappa, i = 0, \dots, n, \qquad (5.11)$$

$$\frac{1}{i+\kappa} \leq \omega_i = \frac{k_{i+1}}{i + k_1} \leq \frac{\kappa}{i+1} \leq \frac{\kappa}{2}, \frac{2}{\kappa} \leq \frac{1}{\omega_i} \leq i + \kappa, i = 1, \dots, n, \quad (5.12)$$

$$\frac{1}{\kappa} \leq \frac{\frac{1}{i-1+\kappa} + 1}{\kappa} \leq \frac{\omega_{i-1}}{\omega_i} = \frac{\omega_{i-1} + k_i}{k_{i+1}} \leq \frac{\frac{\kappa}{i} + \kappa}{1} = \kappa\frac{(1+i)}{i} \leq \frac{3}{2}\kappa, i = 2, \dots, n, \quad (5.13)$$

Using (5.4), (5.12) we obtain

$$t\frac{\kappa}{2} \geq t \geq t\omega_0 = \mid \tilde{B}^{(1)}(t) - B^{(1)} \mid \geq \parallel \tilde{B}^{(1)}(t) - B^{(1)} \parallel,$$

$$t\frac{\kappa}{2} \geq t\frac{\kappa}{i} \geq t\omega_{i-1} = \mid \tilde{B}^{(i)}(t) - B^{(i)} \mid \geq \parallel \tilde{B}^{(i)}(t) - B^{(i)} \parallel, i = 2, \dots, n,$$

hence

$$\parallel \widetilde{W}(t) - W \parallel = \max\left\{\parallel \tilde{B}^{(i)}(t) - B^{(i)} \parallel : i = 1, \dots, n\right\} \leq t\omega_i \leq t\,\frac{\kappa}{2}. \qquad (5.14)$$

Consider also the angles (see Fig. 5.2, where, for brevity, expressions of the form $(t)$ are omitted)

$$\tilde{\alpha}_i(t) = \angle(\tilde{B}^{(i+1)}(t) - \tilde{B}^{(i)}(t), \tilde{B}^{(i)}(t) - \tilde{B}^{(i-1)}(t)), i = 1, \dots, n,$$

$$\delta_i(t) = \angle(\tilde{B}^{(i+1)}(t) - \tilde{B}^{(i)}(t), B^{(i+1)} - \tilde{B}^{(i)}(t)), i = 0, \dots, n-1,$$

$$\gamma_i(t) = \angle(B^{(i+1)} - \tilde{B}^{(i)}(t), \tilde{B}^{(i)}(t) - \tilde{B}^{(i-1)}(t)), i = 1, \dots, n.$$

We first consider the simplest cases, which we will gradually make more complex. Consider first the simplest case (see Fig. 5.2)

$$\alpha_i = \angle\left(B^{(i+1)} - B^{(i)}, B^{(i)} - B^{(i-1)}\right) = \varphi, i = 1, \dots, n. \qquad (5.15)$$

Then, using (5.7), (5.8), by Statement 5.2 we have:

$$\tan\delta_0(t) = \frac{t\omega_0 \sin\varphi}{sk_1 + t\omega_0\cos\varphi} = \frac{\sin\varphi}{\frac{sk_1}{t} + \cos\varphi}, \qquad (5.16)$$

$$\tan\delta_1(t) = \frac{t\omega_1 \sin\varphi}{sk_2 - t\omega_0 + t\omega_1\cos\varphi} = \frac{\sin\varphi}{\frac{sk_2}{t\omega_1} - \frac{1}{\omega_1} + \cos\varphi}, \qquad (5.17)$$

and using (5.9), (5.12) we obtain

$$\frac{k_2}{\omega_1} - k_1 = 1 + k_1 - k_1 = 1,$$

$$\frac{2}{\kappa} \le \frac{1}{\omega_1} = \frac{1 + k_1}{k_2} \le \kappa + 1,$$

$$\frac{k_2}{\omega_1} = 1 + k_1 \le \kappa + 1,$$

and consequently,

$$\tan\ \delta_0(t) - \tan\ \delta_1(t) = \left[\frac{\frac{sk_2}{t\omega_1} - \frac{1}{\omega_1} - \frac{sk_1}{t}}{\left(\frac{sk_1}{t} + \cos\ \varphi\right)\left(\frac{sk_2}{t\omega_1} - \frac{1}{\omega_1} + \cos\ \varphi\right)}\right] \sin\ \varphi$$

$$\geq \left[\frac{\frac{s}{t}\left(\frac{k_2}{\omega_1} - k_1\right) - (\kappa+1)}{\left(\frac{sk_1}{t} + \cos\ \varphi\right)\left(\frac{sk_2}{t\omega_1} - \frac{1}{\omega_1} + \cos\ \varphi\right)}\right] \sin\ \varphi$$

$$\geq \left[\frac{\frac{s}{t} - (\kappa+1)}{\left(\frac{sk_1}{t} + 1\right)\left(\frac{sk_2}{t\omega_1} + 1\right)}\right] \sin\ \varphi$$

$$\geq \frac{t}{s}\ \frac{\left(1 - \frac{t}{s}(\kappa+1)\right)}{\left(\kappa + \frac{t}{s}\right)\left(\kappa + 1 + \frac{t}{s}\right)} \sin\ \varphi.$$

One can show that for $0 < x \leq 0.17\,\frac{1}{\kappa}$, $x \leq 1/5$ we have

$$\frac{1 - (\kappa+1)x}{(\kappa+x)(\kappa+1+x)} \geq \frac{1}{2\kappa(\kappa+1)},$$

and consequently, for $\frac{s}{t} \geq \max\{5;\ \kappa/0.17\}$ the following holds:

$$\tan\ \delta_0(t) - \tan\ \delta_1(t) \geq \frac{t}{s}\ \frac{1}{2\kappa(\kappa+1)} \sin\ \varphi.$$

But then for $s \geq \max\{5;\ \kappa/0.17\}\, t$ we have

$$\tan\ \delta_0(t) - \tan\ \delta_1(t) = (\delta_0(t) - \delta_1(t)) \tan \xi_1(t) \geq \frac{t}{s}\ \frac{1}{2\kappa(\kappa+1)} \sin\ \varphi,$$

where $\xi_1(t) \in [\delta_1(t), \delta_0(t)]$. Since the function $y = \tan x = 1/\cos^2 x$ is monotonically increasing for $x \in [0, \pi/2)$, and $\delta_0(t) \geq \xi_1(t) \geq \delta_1(t)$, we get

$$(\delta_0(t) - \delta_1(t)) \tan \delta_0(t) \geq (\delta_0(t) - \delta_1(t)) \tan \xi_1 \geq \frac{t}{s}\ \frac{1}{2\kappa(\kappa+1)} \sin\ \varphi,$$

where

$$\tan \delta_0(t) = \frac{1}{\cos^2 \delta_0(t)} = 1 + \tan^2 \delta_0 \leq 1 + \frac{1}{25} = 1.04,$$

since for $s \geq 5t$ (taking into account that $1 \leq k_1$)

$$\tan\ \delta_0(t) = \frac{\sin \varphi}{\frac{s}{t} k_1 + \cos\ \varphi} \leq \frac{1}{\frac{s}{t} k_1} \leq \frac{1}{5k_1} \leq \frac{1}{5},$$

hence

$$\tilde{\alpha}_1(t) - \alpha_1 = \delta_0(t) - \delta_1(t) \geq \frac{t}{4s \cdot 1.04} \; \frac{1}{2\kappa(\kappa+1)} \sin\varphi$$

$$= \frac{t}{4.16\, s} \; \frac{1}{2\kappa(\kappa+1)} \sin\varphi. \qquad (5.18)$$

Now let $i \in \{2, \dots, n-1\}$. Then

$$sk_i = \left|B^{(i)} - B^{(i-1)}\right|,$$

$$\tan\ \delta_{i-1}(t) = \frac{t\omega_{i-1}\sin\ \varphi}{sk_i - t\omega_{i-2} + t\cos\ \varphi} = \frac{\sin\varphi}{\frac{sk_i}{t\omega_{i-1}} - \frac{\omega_{i-2}}{\omega_{i-1}} + \cos\ \varphi},$$

$$\tan\ \delta_i(t) = \frac{t\omega_i\sin\ \varphi}{sk_{i+1} - t\omega_{i-1} + t\omega_i\cos\ \varphi} = \frac{\sin\varphi}{\frac{sk_{i+1}}{t\omega_i} - \frac{\omega_{i-1}}{\omega_i} + \cos\ \varphi},$$

and using (5.8)–(5.13) we obtain

$$\tan \delta_{i-1}(t) - \tan \delta_i(t) = \left[\frac{\frac{sk_{i+1}}{t\omega_i} - \frac{\omega_{i-1}}{\omega_i} - \frac{sk_i}{t\omega_{i-1}} + \frac{\omega_{i-2}}{\omega_{i-1}}}{\left(\frac{sk_i}{t\omega_{i-1}} - \frac{\omega_{i-2}}{\omega_{i-1}} + \cos\varphi\right)\left(\frac{sk_{i+1}}{t\omega_i} - \frac{\omega_{i-1}}{\omega_i} + \cos\varphi\right)}\right] \sin\varphi \geq$$

$$\geq \left[\frac{\frac{s}{t}\left(\frac{k_{i+1}}{\omega_i} - \frac{k_i}{\omega_{i-1}}\right) - \frac{\omega_{i-1}}{\omega_i} + \frac{\omega_{i-2}}{\omega_{i-1}}}{\left(\frac{sk_i}{t\omega_{i-1}} + \cos\varphi\right)\left(\frac{sk_{i+1}}{t\omega_i} + \cos\varphi\right)}\right] \sin\varphi \geq$$

$$\geq \left[\frac{\frac{s}{t} - \frac{3}{2}\kappa + \frac{1}{\kappa}}{\left(\frac{sk_i}{t\omega_{i-1}} + 1\right)\left(\frac{sk_{i+1}}{t\omega_i} + 1\right)}\right] \sin\varphi \geq \left[\frac{\frac{s}{t} - \frac{3}{2}\kappa}{\left(\frac{sk_i}{t\omega_{i-1}} + 1\right)\left(\frac{sk_{i+1}}{t\omega_i} + 1\right)}\right] \sin\varphi \geq$$

$$\geq \frac{t}{s}\left[\frac{1 - \frac{3}{2}\kappa\frac{t}{s}}{\left(\frac{k_i}{\omega_{i-1}} + \frac{t}{s}\right)\left(\frac{k_{i+1}}{\omega_i} + \frac{t}{s}\right)}\right] \sin\varphi \geq \frac{t}{s}\left[\frac{1 - \frac{3}{2}\kappa\frac{t}{s}}{\left(n + \kappa + \frac{t}{s}\right)\left(n + \kappa + \frac{t}{s}\right)}\right] \sin\varphi,$$

$i = 2, \dots, n-1$.

It can be shown that for $0 < x \leq \frac{1}{3(n+\kappa)}$, $x \leq \frac{1}{5}$ the following holds:

$$\frac{1 - \frac{3}{2}\kappa x}{(n+\kappa+x)(n+\kappa+x)} \geq \frac{1}{2(n+\kappa)^2},$$

and consequently, for $\frac{s}{t} \geq \max\{5;\ 3(n+\kappa)\}$ we have

$$\tan\ \delta_{i-1}(t) - \tan\ \delta_i(t) \geq \frac{t}{s}\ \frac{1}{2(n+\kappa)^2}\sin\varphi.$$

But then for $s \geq \max\{5;\ 3(n+\kappa)\}\,t$ we obtain

$$\tan \delta_{i-1}(t) - \tan \delta_i(t) = (\delta_{i-1}(t) - \delta_i(t))\dot{\tan}\ \xi_i(t) \geq \frac{t}{s}\frac{1}{2(n+\kappa)^2}\sin\varphi,$$

where $\xi_i(t) \in [\delta_i(t), \delta_{i-1}(t)]$. Since the function $y = \dot{\tan}\, x = 1/\cos^2 x$ is monotonically increasing for $x \in [0, \pi/2)$ and $\delta_{i-1}(t) \geq \xi_i(t) \geq \delta_i(t)$, we get

$$(\delta_{i-1}(t) - \delta_i(t))\ \dot{\tan}\,\delta_{i-1}(t)\ \geq (\delta_{i-1}(t) - \delta_i(t))\,\dot{\tan}\,\xi_i(t) \geq \frac{t}{s}\ \frac{1}{2(n+\kappa)^2}\sin\ \varphi,$$

where

$$\dot{\tan}\,\delta_{i-1}(t) = \frac{1}{\cos^2\delta_{i-1}(t)} = 1 + \tan^2\delta_{i-1}(t) \leq 1 + \frac{1}{100} = 1.01.$$

Indeed, for $s \geq \max\{5;\ 3(n+\kappa)\}\,t$ (taking into account that $1 \leq k_1$) by virtue of (5.11), (5.13) we have (since $k_i/\omega_{i-1} = i - 1 + k_1$, $\omega_{i-1}/\omega_i \leq \frac{3}{2}\kappa$, $i \geq 2$, $\kappa \geq 1$, $n \geq 1$):

$$\frac{sk_i}{t\omega_{i-1}} - \frac{\omega_{i-2}}{\omega_{i-1}} + \cos\ \varphi \geq \frac{s}{t}(i - 1 + k_1) - \frac{3}{2}\kappa$$

$$\geq 3(n+\kappa)(i - 1 + k_1) - \frac{3}{2}\kappa \geq 6(n+\kappa) - \frac{3}{2}\kappa = \frac{9}{2}\kappa + 6n \geq 10.5,$$

hence

$$\tan\ \delta_{i-1}(t) = \frac{\sin\varphi}{\frac{sk_i}{t\omega_{i-1}} - \frac{\omega_{i-2}}{\omega_{i-1}} + \cos\ \varphi} \leq \frac{1}{10.5} < \frac{1}{10},$$

$$\tilde{\alpha}_i(t) - \alpha_i = \delta_{i-1}(t) - \delta_i(t) \geq \frac{t}{s\cdot 1.01}\ \frac{1}{2(n+\kappa)^2}\sin\ \varphi = \frac{t}{2.02\,s}\ \frac{1}{2(n+\kappa)^2}\sin\ \varphi, i = 2, \ldots, n-1. \qquad (5.19)$$

Moreover, the case $i = n$ remains. Then, if we introduce the quantities

$$k_{n+1} = \frac{\left|B^{(n+1)} - B^{(n)}\right|}{s} \in [1, \kappa],\ \omega_n = \frac{k_{n+1}}{n+k_1},\ \delta_n(t)\ = \arctan\frac{\sin\varphi}{\frac{sk_{n+1}}{t\omega_n} - \frac{\omega_{n-1}}{\omega_n} + \cos\varphi},$$

we again obtain a case similar to that described above for $i = 2, \dots, n-1$, i.e., formula (5.19) will also hold in this case, from which it follows that for $s \geq \max\{5;\ 3(n+\kappa)\}\, t$ we have

$$\tilde{\alpha}_n(t) - \alpha_n \geq \delta_{n-1}(t) - \delta_n(t) \geq \frac{t}{s\cdot 1.01}\ \frac{1}{2(n+\kappa)^2}\sin\varphi = \frac{t}{2.02\, s}\ \frac{1}{(n+\kappa)^2}\sin\varphi. \quad (5.20)$$

Now note that by virtue of (5.14), (5.18), (5.19), (5.20) we have:

$$\|\widetilde{W} - W\| = \left\|\left(B^{(1)}, \dots, B^{(n)}\right) - \left(\tilde{B}^{(1)}(t), \dots, \tilde{B}^{(n)}(t)\right)\right\| \leq t\frac{\kappa}{2},$$

$$\alpha_i - \frac{t}{s}\frac{1}{\max\{2.02(n+\kappa)^2; 8.32\kappa(\kappa+1)\}}\sin\varphi \geq \tilde{\alpha}_i,\ i = 1, \dots, n. \quad (5.21)$$

Observe that by (5.5), (5.6), (5.12)

$$\left|\tilde{B}^{(i+1)}(t) - \tilde{B}^{(i)}(t)\right| \geq s - t\omega_{i-1} \geq s - t\kappa, i = 1, \dots, n,$$

$$\left|\tilde{B}^{(1)}(t) - \tilde{B}^{(0)}(t)\right| = \left|\tilde{B}^{(1)}(t) - B^{(0)}\right| \geq s.$$

Furthermore, by (5.4), (5.5) we have (also assuming $\kappa \geq 2 \Rightarrow 1 + \kappa/2 \leq \kappa$):

$$\left|B^{(i)} - B^{(j)}\right| = \left|B^{(i)} - \tilde{B}^{(i)}(t) + \tilde{B}^{(i)}(t) - \tilde{B}^{(j)}(t) + \tilde{B}^{(j)}(t) - B^{(j)}\right| \leq$$

$$\leq \left|\tilde{B}^{(i)}(t) - B^{(i)}\right| + \left|\tilde{B}^{(i)}(t) - \tilde{B}^{(j)}(t)\right| + \left|B^{(j)} - \tilde{B}^{(j)}(t)\right| \Rightarrow$$

$$\Rightarrow \left|\tilde{B}^{(i)}(t) - \tilde{B}^{(j)}(t)\right| \geq s - t\omega_{i-1} - t\omega_{j-1} = s - t(\omega_{i-1} + \omega_{j-1}) \geq$$

$$\geq s - t\kappa, i, j = 1, \dots, n+1, i \neq j, \quad (5.22)$$

$$\left|\tilde{B}^{(0)}(t) - \tilde{B}^{(i)}(t)\right| = \left|B^{(0)} - \tilde{B}^{(i)}(t)\right| \geq s - t\omega_{i-1} \geq$$

$$\geq s - t\kappa, i = 1, \dots, n, \quad (5.23)$$

$$\left|\tilde{B}^{(n+1)}(t) - \tilde{B}^{(i)}(t)\right| = \left|B^{(n+1)} - \tilde{B}^{(i)}\right| \geq s - t\omega_{i-1} \geq$$

$$\geq s - t\kappa,\ i = 1, \dots, n. \quad (5.24)$$

Thus (using $\varphi \in (0, \pi/2)$), by virtue of (5.21)–(5.24) we obtain that in the case where (5.4), (5.8) hold for $\left(\tilde{B}^{(1)}(t), \dots, \tilde{B}^{(n)}(t)\right)$, there exist constants

$$\omega = \omega(s, d, \varphi, n) = \frac{1}{\max\{2.02(n+\kappa)^2;\ 8.32\kappa(\kappa+1)\}}\sin\varphi > 0, \quad (5.25)$$

$$\theta = \theta(s, d, \varphi, n) = \max\{5;\ 3(n+\kappa);\ \kappa/0.17\},$$

such that for all $t \in (0, s/\theta]$ (i.e., for $s \geq \max\{5;\ 3(n+\kappa);\ \kappa/0.17\}\, t$) we have:

$$\left(\tilde{B}^{(1)},\dots,\tilde{B}^{(n)}\right)\in \boldsymbol{S}^{(n,s-t\kappa)}(A,B,\varphi-\omega t), \tag{5.26}$$

$$\left\|\left(B^{(1)},\dots,B^{(n)}\right)-\left(\tilde{B}^{(1)},\dots,\tilde{B}^{(n)}\right)\right\|\le t\frac{\kappa}{2}. \tag{5.27}$$

Let us now consider the cases where the condition

$$\alpha_i=\angle(B^{(i+1)}-B^{(i)},\,B^{(i)}-B^{(i-1)})=\varphi, i=1,\dots,n,$$

is violated.

First, let

$$\alpha_i\in(0,\varphi], i=1,\dots,n. \tag{5.28}$$

Suppose (5.4) and (5.8) hold. Again we use the previously defined angles $\tilde{\alpha}_i(t), \delta_{i-1}(t), \gamma_i(t), i=1,\dots,n$. Then

$$\alpha_1=\varphi-\sigma_1, \alpha_2=\varphi-\sigma_2, \sigma_1,\sigma_2\in[0,\varphi).$$

Taking into account the dependence of the angle $\tilde{\alpha}_1(t)$ on $\sigma_1,\sigma_2$ and of the angles $\delta_0(t),\delta_1(t)$ on $\sigma_1$ and $\sigma_2$ respectively, we now write:

$$\tilde{\alpha}_1(t,\sigma_1,\sigma_2), \delta_0(t,\sigma_1), \delta_1(t,\sigma_2).$$

As before, using Statement 5.2, we have:

$$\gamma_1(t)=\angle(B^{(2)}-\tilde{B}^{(1)}(t),\,\tilde{B}^{(1)}(t)-A)=\alpha_1-\delta_0(t,\sigma_1),$$

$$\tilde{\alpha}_1(t,\sigma_1,\sigma_2)=\gamma_1(t)+\delta_1(t,\sigma_2)=\varphi-\sigma_1-\delta_0(t,\sigma_1)+\delta_1(t,\sigma_2),$$

where (analogously to (5.16), (5.17))

$$\tan\delta_0(t,\sigma_1)=\frac{t\omega_0\sin(\varphi-\sigma_1)}{sk_1+\mathrm{t}\cos(\varphi-\sigma_1)}=\frac{\sin(\varphi-\sigma_1)}{\frac{sk_1}{t}+\cos(\varphi-\sigma_1)}, \tag{5.29}$$

$$\tan\delta_1(t,\sigma_2)=\frac{t\omega_1\sin(\varphi-\sigma_2)}{sk_2-t\omega_0+t\omega_1\cos(\varphi-\sigma_2)}$$

$$=\frac{\sin(\varphi-\sigma_2)}{\frac{sk_2}{t\omega_1}-\frac{1}{\omega_1}+\cos(\varphi-\sigma_2)}. \tag{5.30}$$

We find the solution of the problem

$$\tilde{\alpha}_1(t,\sigma_1,\sigma_2)=\varphi-\sigma_1-\delta_0(t,\sigma_1)+\delta_1(t,\sigma_2)\to\max;\sigma_1,\sigma_2\in[0,\varphi], \tag{5.31}$$

Where $\delta_0(t,\sigma_1),\delta_1(t,\sigma_2)$ satisfy (5.29), (5.30); i.e., we determine the "worst" case, namely the maximum possible value of the angle $\tilde{\alpha}_1(\sigma_1,\sigma_2)$ under the same strategy with parameters $t,\omega_i,s$ satisfying (5.4), (5.8) and the inequality $s\ge\max\{5;\ 3(n+\kappa);\ \kappa/0.17\}\,t$, and with $\sigma_1,\sigma_2\in[0,\varphi]$. As will be shown, under these conditions

$$\tilde{\alpha}_1(t,\sigma_1,\sigma_2)\le\varphi-(\delta_0(t,0)-\delta_1(t,0))=\varphi-(\delta_0(t)-\delta_1(t)),$$

where $\delta_0(t), \delta_1(t)$ satisfy (5.16), (5.17); i.e., the "worst" case corresponds to the previously considered case $\alpha_1 = \varphi, \alpha_2 = \varphi$.

In problem (5.31), a decomposition into two problems (**Problem 1** and **Problem 2**) is possible:

**Problem 1:**

$$\vartheta(t,\sigma_1) = \sigma_1 + \delta_0(t,\sigma_1) = \sigma_1 + \arctan\left(\frac{\sin(\varphi - \sigma_1)}{\frac{sk_1}{t} + \cos(\varphi - \sigma_1)}\right) \to \min\,;$$

$\sigma_1 \in [0, \varphi]$. (5.32)

Note that the function $y = x + \arctan\frac{\sin(\varphi - x)}{z + \cos(\varphi - x)}$, where $\varphi \in (0, \pi/2), z > 0$ are parameters, has the following derivative (with respect to $x$) for $x \in [0, \varphi], \varphi \in (0, \pi/2), z > 0$ ($z, \varphi$ are parameters, $\varphi - x \in [0, \varphi] \subseteq [0, \pi/2)$):

$$\frac{dy}{dx} = \frac{z^2 + z\cos(\varphi - x)}{z^2 + 2z\cos(\varphi - x) + 1} \geq \frac{z^2 + z\cos\varphi}{(z+1)^2} > 0,$$

which implies that the objective function $\vartheta(t,\sigma_1)$ in problem (5.32) is monotonically increasing for $\sigma_1 \in [0, \varphi]$. Consequently, the minimum value of this function is attained at $\sigma_1 = 0$ and equals $\delta_0(t, 0) = \delta_0(t)$.

Thus, if $\sigma_1 \in [0, \varphi]$, the contribution of the term $-\vartheta(t,\sigma_1)$ to the objective function $\tilde{\alpha}_1(t, \sigma_1, \sigma_2) = -\vartheta(t,\sigma_1) + \delta_1(t, \sigma_2)$ is no greater than $-\delta_0(t) = -\min\{\vartheta(t,\sigma_1) = \sigma_1 + \delta_0(t,\sigma_1) \mid \sigma_1 \in [0,\varphi]\}$, which corresponds to the case $\sigma_1 = 0, \alpha_1 = \varphi - \sigma_1 = \varphi$.

**Problem 2:**

$$\delta_1(t,\sigma_2) = \arctan\left(\frac{\sin(\varphi - \sigma_2)}{\frac{sk_2}{t\omega_1} - \frac{1}{\omega_1} + \cos(\varphi - \sigma_2)}\right) \to \max\,; \sigma_2 \in [0, \varphi]. \quad (5.33)$$

Note that (since $s \geq 5t \Rightarrow \frac{sk_2}{t\omega_1} - \frac{1}{\omega_1} \geq \frac{4}{\omega_1} > 0$) the function $y = \arctan \frac{\sin(\varphi - x)}{z + \cos(\varphi - x)}$ has the following derivative (with respect to $x$) for $x \in [0, \varphi], \varphi \in (0, \pi/2), z > 0$ ($z, \varphi$ are parameters):

$$\frac{dy}{dx} = \frac{-z\cos(\varphi - x) - 1}{z^2 + 2z\cos(\varphi - x) + 1} \leq -\frac{z\cos\varphi + 1}{(z+1)^2} < 0,$$

which implies that the objective function $\delta_1(t, \sigma_2)$ in problem (5.33) is monotonically decreasing for $\sigma_2 \in [0, \varphi]$. Consequently, the maximum value of this function is attained at $\sigma_2 = 0$ and equals $\delta_1(t, 0) = \delta_1(t)$.

Thus, if $\sigma_2 \in [0, \varphi]$, the contribution of the term $\delta_1(\sigma_2)$ to the objective function $\tilde{\alpha}_1(\sigma_1, \sigma_2) = \vartheta(\sigma_1) + \delta_1(\sigma_2)$ is no greater than $\delta_1 = \min\{\delta_1(\sigma_2) \mid \sigma_2 \in [0, \varphi]\}$, which corresponds to the case $\sigma_2 = 0, \alpha_2 = \varphi - \sigma_2 = \varphi$.

Therefore, the maximum value in problem (5.31) is attained at $\sigma_1 = 0, \sigma_2 = 0$ and equals $\varphi - (\delta_0(t) - \delta_1(t))$; i.e., it coincides with the value of $\tilde{\alpha}_1(t)$ for the case $\alpha_1 = \varphi, \alpha_2 = \varphi$. An analogous result is obtained for the remaining $\alpha_i, \tilde{\alpha}_i, i = 1, \dots, n$. Hence, the previously obtained results also hold in the case (5.28).

The method also remains valid (see the second case in Statement 1.2) if $\alpha_i \in [-\varphi, 0), i = 1, \dots, n$.

Consider now the case where the signs of $\alpha_i$ alternate:

$$\alpha_i \in [-\varphi, 0) \cup (0, \varphi], i = 1, \dots, n.$$

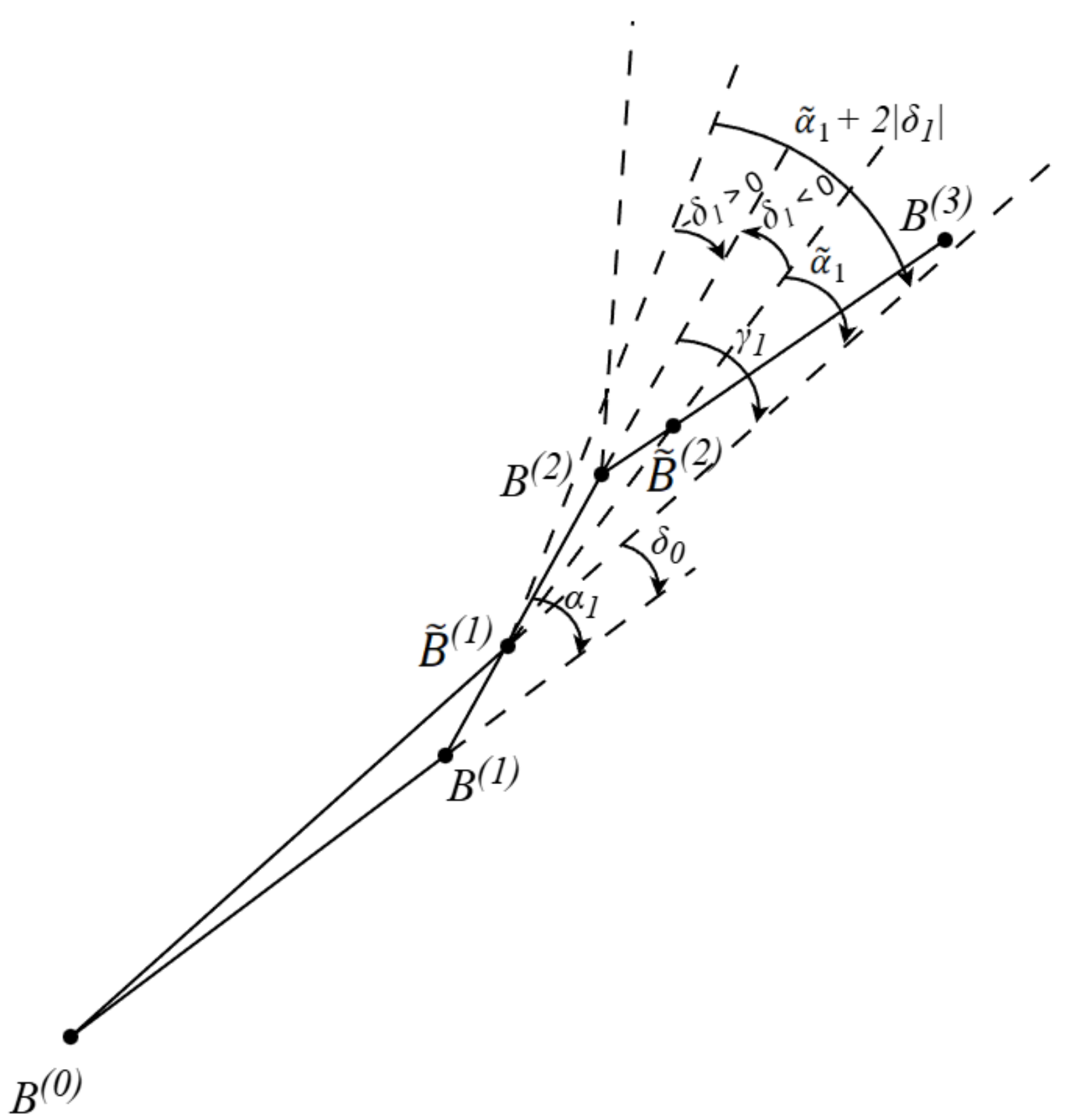

**Fig. 5.3**

As shown in Fig. 5.3 (where, for brevity, expressions of the form $(t)$ are omitted), the alternation of signs of the quantities $\alpha_i, i = 1, \dots, n$, leads, in the case $\alpha_i > 0, \alpha_{i+1} < 0$, to an increase in the value of $\alpha_i - \tilde{\alpha}_i > 0$ compared to the case $\alpha_i > 0, \alpha_{i+1} > 0$. Correspondingly, in the case $\alpha_i < 0, \alpha_{i+1} > 0$, we obtain an increase in the value of $\tilde{\alpha}_i - \alpha_i > 0$ compared to the case $\alpha_i < 0, \alpha_{i+1} < 0$. For example, if (see Fig. 5.3)

$$\alpha_1 > 0, \alpha_2 = \angle(B^{(3)} - B^{(2)}, B^{(2)} - B^{(1)}) = \angle(\tilde{B}^{(2)} - B^{(2)}, B^{(2)} - \tilde{B}^{(1)}) < 0,$$

then the sign of the quantity

$$\delta_1 = \angle(\tilde{B}^{(2)} - \tilde{B}^{(1)}, B^{(2)} - \tilde{B}^{(1)})$$

becomes negative (see Statement 5.2; instead of positive in the case $\alpha_2 > 0$), but the magnitude $|\ \delta_1\ |$ remains the same. Thus, instead of

$$\alpha_1 > 0, \alpha_2 > 0, \delta_0 > 0, \delta_1 > 0, \alpha_1 - \tilde{\alpha}_1 = \delta_0 - \delta_1 = \delta_0 - |\delta_1| > 0,$$

we have

$$\alpha_1 > 0, \alpha_2 < 0, \delta_0 > 0, \delta_1 < 0, \alpha_1 - \tilde{\alpha}_1 = \delta_0 - \delta_1 = \delta_0 + |\delta_1| > 0,$$

i.e., the difference $\alpha_1 - \tilde{\alpha}_1$ remains positive but increases by $2|\delta_1| > 0$. A similar situation occurs at each alternation of signs of the quantities $\alpha_i, i = 1, \dots, n$.

Let us make one more remark concerning the angles $\alpha_i, \tilde{\alpha}_i(t), i = 1, \dots, n$. Suppose we are tasked with ensuring the inequalities

$$\tilde{\alpha}_i(t) = \angle\left(\tilde{B}^{(i+1)}(t) - \tilde{B}^{(i)}(t), \tilde{B}^{(i)}(t) - \tilde{B}^{(i-1)}(t)\right) \in (0, \varphi - \omega t], \quad i = 1, \dots, n, \tag{5.34}$$

for all sufficiently small $t > 0$, where, for example, $\omega = \omega(s, d, \varphi, n) > 0$ satisfies (5.25). Then, in the case where for several consecutive indices we have $\alpha_i(t) \in (\varphi - \omega t, \varphi]$, we proceed according to the described method, aiming to achieve (5.34). Otherwise, we begin applying the method starting from the first index $i_0 \in \{1, \dots, n\}$ for which $\alpha_{i_0} \in (\varphi - \omega t, \varphi]$, and stop its application at the first index $i_1 \in \{i_0 + 1, \dots, n\}$ for which $\alpha_{i_1} \in (\varphi - \omega t, \varphi]$, setting $\omega_i = 0$ for $i = i_1, \dots, i_2$, where $\alpha_i \in (\varphi - \omega t, \varphi], i = i_1, \dots, i_2$, and $\alpha_{i_2+1} \in (0, \varphi - \omega t]$. After each such occurrence of $\alpha_i \in (\varphi - \omega t, \varphi]$, if later we again encounter a case with $\alpha_j \in (\varphi - \omega t, \varphi]$ where $j \in \{i + 1, \dots, n\}$, we apply the method "anew", i.e., using the values $\omega_i$ in the same order as if the equality $i = 1$ held (that is, for any intermediate $i \in \{1, \dots, n\}$ such that $\alpha_i \in (0, \varphi - \omega t], \alpha_{i+1} \in (\varphi - \omega t, \varphi]$, we apply the described method, namely we set $\tilde{B}^{(i)}(t)$ according to (5.4), but using the value $\omega_1$ instead of $\omega_i$, then using $\omega_2$ instead of $\omega_{i+1}$, and so on). □

## 5.2. Some properties of the angle $\Psi(C) = \angle(B - C, C - A)$, where $A = (0, -1), B = (0, 1)$

For simplicity of notation, let $A = (0,-1), B = (0,1)$ throughout this section. We will be interested in some properties of the angle $\Psi(C) = \angle(B - C, C - A)$, where $C \in \mathbb{R}^2 \setminus \{A, B\}$. Note that

$$u(C) = \cos\ \Psi(C) = \frac{\langle B - C,\ C - A\rangle}{\mid B - C \mid \cdot \mid C - A \mid} =$$

$$\frac{\langle(-c_1, 1 - c_2),\ (c_1, c_2 + 1)\rangle}{\sqrt{c_1^2 + (1 - c_2)^2}\ \sqrt{c_1^2 + (1 + c_2)^2}} = \frac{1 - c_1^2 - c_2^2}{\sqrt{(1 + c_1^2 + c_2^2)^2 - 4c_2^2}},$$

and taking into account the orientation of the angles,

$$\Psi(C) = \text{sign}\ (c_1)\arccos\ u(C),$$

where

$$\text{sign}\ (c_1) = \begin{cases} 1, & \text{if } c_1 \geq 0, \\ -1, & \text{if } c_1 < 0. \end{cases}$$

It is easy to show that

$$\frac{\partial\Psi(C)}{\partial c_1} = \frac{2(1+c_1^2-c_2^2)}{(1+c_1^2+c_2^2)^2-4c_2^2},\ \frac{\partial\Psi(C)}{\partial c_2} = \frac{4c_1c_2}{(c_1^2+c_2^2+1)^2-4c_2^2},$$

i.e., these functions are infinitely differentiable on the whole space $\mathbb{R}^2$ except at the points $A = (0,-1), B = (0,1)$ (i.e., on $\mathbb{R}^2 \setminus \{A, B\}$), where their common denominator vanishes.

By virtue of the smoothness of the function $\Psi(C)$ defined on the set $\mathbb{R}^2 \setminus \{A, B\}$, we conclude that for a given $\varphi \in (0, \pi)$ this function has bounded partial derivatives with respect to $c_1, c_2$ on the bounded closed set

$$\boldsymbol{S}_\varphi^{(1)}(s) = \{C \in \boldsymbol{S}_\varphi^{(1)} \mid\ \mid C - A \mid \geq s,\ \mid C - B \mid \geq s\}$$

for any fixed $s > 0$, where [1]

$$\boldsymbol{S}_\varphi^{(1)} = \{C \in \mathbb{R}^2 \mid \Psi(C) = \angle(B - C, C - A) \in [-\varphi, \varphi]\} = \widehat{\boldsymbol{S}}_\varphi^{(1)} \cup \widehat{\boldsymbol{S}}_\varphi^{(1)},$$

$$\widehat{\boldsymbol{S}}_\varphi^{(1)} = \{C \in \mathbb{R}^2 \mid \Psi(C) = \angle(B - C, C - A) \in [0, \varphi]\},$$

$$\widehat{\boldsymbol{S}}_\varphi^{(1)} = \{C \in \mathbb{R}^2 \mid \Psi(C) = \angle(B - C, C - A) \in [-\varphi, 0]\},$$

$$\breve{\boldsymbol{S}}_{\varphi}^{(1)} = \left\{(x, y) \in \mathbb{R}^2 \mid 0 \le x \le \tan\frac{\varphi}{2},\ -\sqrt{\frac{1}{\sin^2\varphi} - (x + \cot\varphi)^2} \le y \right.$$

$$\left. \le \sqrt{\frac{1}{\sin^2\varphi} - (x + \cot\varphi)^2}\right\},$$

$$\widehat{\boldsymbol{S}}_{\varphi}^{(1)} = \left\{(x, y) \in \mathbb{R}^2 \mid -\tan\frac{\varphi}{2} \le x \le 0,\ -\sqrt{\frac{1}{\sin^2\varphi} - (x - \cot\varphi)^2} \le y \right.$$

$$\left. \le \sqrt{\frac{1}{\sin^2\varphi} - (x - \cot\varphi)^2}\right\}.$$

The convex sets $\breve{\boldsymbol{S}}_{\varphi}^{(1)}$ and $\widehat{\boldsymbol{S}}_{\varphi}^{(1)}$ will be called the right and left parts of the set $\boldsymbol{S}_{\varphi}^{(1)}$, respectively. Denote

$$\Psi_{c_i}(C) = \frac{\partial\Psi(C)}{\partial c_i}, i = 1,2,$$

$$L_{\Psi}(s) = \max\left\{\left|\Psi_{c_1}(C)\right| + \left|\Psi_{c_2}(C)\right| \,\middle|\, C \in \boldsymbol{S}_{\varphi_0}^{(1)}(s)\right\}.$$

For any vector $C = (c_1, c_2) \in \mathbb{R}^2$, we will use the notation:

$$|C| = \sqrt{c_1^2 + c_2^2}, \|C\| = \max\{|c_1|, |c_2|\}.$$

**Statement 5.4.** *Let* $C = (c_1, c_2), f(C)$ *be a function defined and continuously differentiable on some open set* $\boldsymbol{S} \subseteq \mathbb{R}^2, f_{c_i}(C) = \frac{\partial f(C)}{\partial c_i}, i = 1,2,$ *and let* $\tilde{\boldsymbol{S}}$ *be an arbitrary nonempty subset of* $\boldsymbol{S}$ *such that for some number* $M > 0$ *we have*

$$\left|f_{c_1}(C)\right| + \left|f_{c_2}(C)\right| \le M, \forall C \in \tilde{\boldsymbol{S}}.$$

Then for all $C, \bar{C} = (\bar{c}_1, \bar{c}_2) \in \tilde{\boldsymbol{S}}$, if the segment $[C, \bar{C}] \subseteq \tilde{\boldsymbol{S}}$, then

$$|f(C) - f(\bar{C})| \le M\|C - \bar{C}\|.$$

***Proof.*** Note that there exists $\xi \in [C, \bar{C}]$ such that

$$f(C)-f(\bar{C})=\langle (f_{c_1}(\xi),f_{c_2}(\xi)),\,C-\bar{C}\rangle = f_{c_1}(\xi)(c_1-\bar{c}_1)+f_{c_2}(\xi)(c_2-\bar{c}_2),$$

hence

$$|f(C)-f(\bar{C})| \le \left|f_{c_1}(\xi)\right| |c_1-\bar{c}_1| + \left|f_{c_2}(\xi)\right| |c_2-\bar{c}_2|$$

$$\le \left(\left|f_{c_1}(\xi)\right| + \left|f_{c_2}(\xi)\right|\right)\|C-\bar{C}\| \le M\|C-\bar{C}\|. \square$$

Corollary of Statement 5.4 (using the continuous differentiability of the Function $\Psi(C)$ on the open set $\mathbb{R}^2 \setminus \{A,B\}$) is

**Statement 5.5.** *Let for some* $\varphi \in (0,\pi)$, $s>0$, $t\in (0, \min\{\frac{1}{\sqrt{2}};\ \sin(\pi-\varphi)\}\,s)$, *and* $C,\bar{C}\in \boldsymbol{S}_\varphi^{(1)}$ *the following hold*: $C\in \boldsymbol{S}_\varphi^{(1)}(s)$, $\| C-\bar{C} \|\le t$. *Then* $[C,\bar{C}] \subseteq \boldsymbol{S}_\varphi^{(1)}(s-\sqrt{2}\,t)$, *and by Statement 5.4,*

$$|\ \Psi(C)-\Psi(\bar{C})\ | \le L_\Psi(s-\sqrt{2}\,t)\|C-\bar{C}\|.$$

***Proof.*** If $C,\bar{C}\in \boldsymbol{S}_\varphi^{(1)}$ belong to the same part, then by convexity of each part, $[C,\bar{C}]\subseteq \boldsymbol{S}_\varphi^{(1)}$, $|C-\bar{C}|\le \sqrt{2}\|C-\bar{C}\|\le \sqrt{2}t$, $C\in \boldsymbol{S}_\varphi^{(1)}(s)$, hence

$$[C,\bar{C}] \subseteq \boldsymbol{S}_\varphi^{(1)}(s-\sqrt{2}\,t). \qquad (5.35)$$

In the case $\varphi\in(0,\pi/2]$, the set $\boldsymbol{S}_\varphi^{(1)}$ itself is convex, so (5.35) also holds. Now let $\varphi \in (\pi/2,\pi)$, and let the points $C,\bar{C}$ belong to different parts, and suppose, for example,

$$c_1>0, \bar{c}_1<0. \qquad (5.36)$$

Then

$$\| C-\bar{C} \|\le t \ \Rightarrow\ c_1-\bar{c}_1\le t \ \Rightarrow\ c_1<t. \qquad (5.37)$$

Note further that for $[C,\bar{C}]\subseteq \boldsymbol{S}_\varphi^{(1)}$ not to hold, one of the following two conditions must be satisfied:

$$\text{(a)}\ \max\{c_2,\bar{c}_2\}>1, \quad \text{(b)}\ \min\{c_2,\bar{c}_2\}<1.$$

Without loss of generality, we may assume that the point $C$ lies on the boundary of $\breve{\boldsymbol{S}}_\varphi^{(1)}$ and the point $\bar{C}$ lies on the boundary of $\widehat{\boldsymbol{S}}_\varphi^{(1)}$. In this case, using the monotonicity properties of the curve bounding the set $\widehat{\boldsymbol{S}}_{\varphi}^{(1)}$ to the left of the line $\overrightarrow{A,B}$, as well as the monotonicity of the curve bounding the set $\breve{\boldsymbol{S}}_{\varphi}^{(1)}$ to the right of the line $\overrightarrow{A,B}$, it is easy to show that one of the following conditions now holds:

$$\text{(c) } c_2, \bar{c}_2 > 1, \quad \text{(d) } c_2, \bar{c}_2 < 1.$$

Suppose, for example, case (c) holds. Denote $\gamma = \angle(B - A,\ C - B)$. It is easy to show (see Fig. 4) that

$$\gamma > \pi - \varphi > 0 \qquad\qquad (5.38)$$

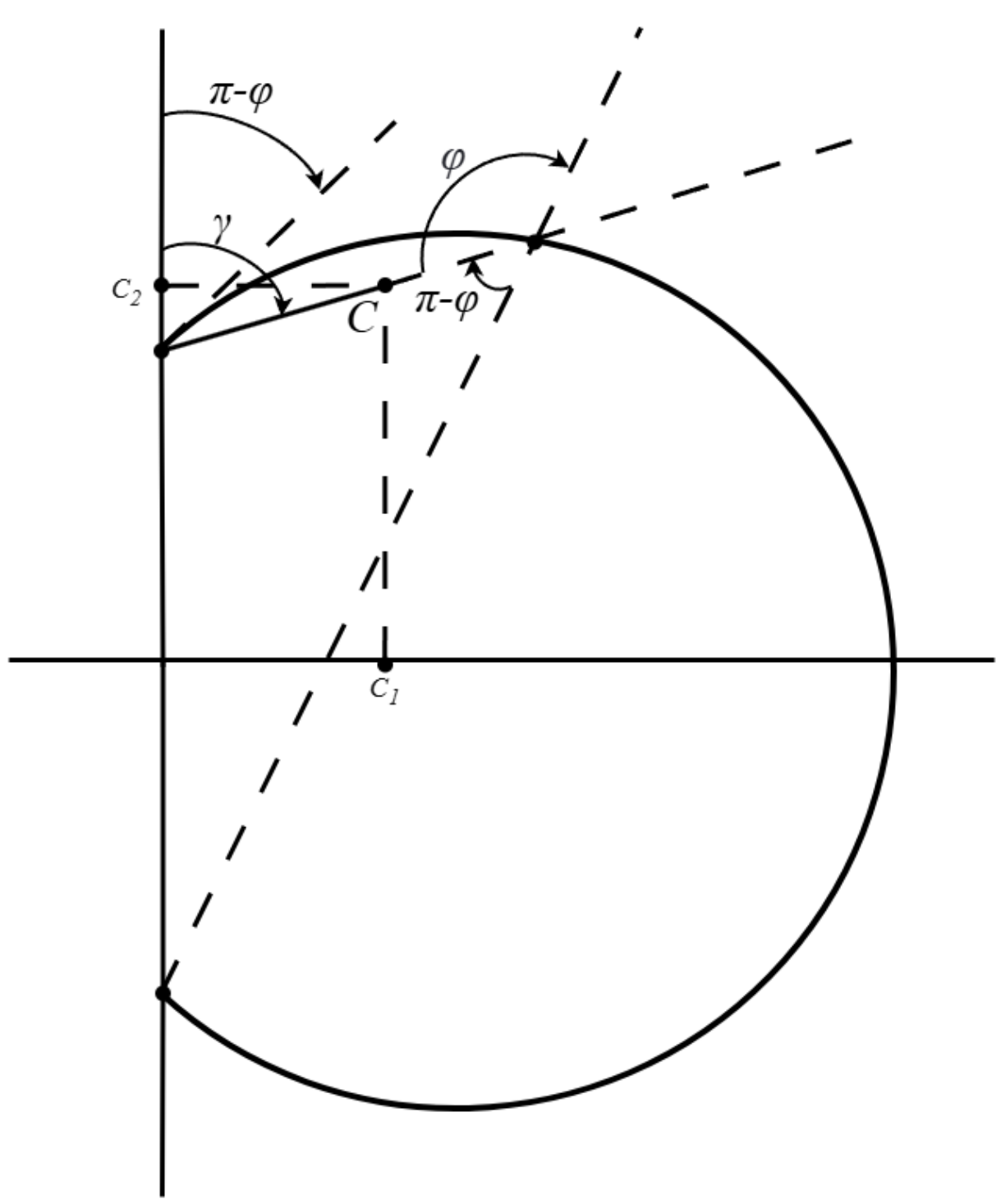

**Fig. 5.4**

(since $\pi - \varphi$ is the angle between the vector $B - A = (0,2)$ and the vector $D - B = (d_1, d_2 - 1)$, where $D = (d_1, d_2)$ is any point lying on the tangent to the circle given by the equation $(x + \cot\ \varphi)^2 + y^2 = 1/\sin^{\ 2} \varphi$, drawn at point $B$ and located strictly to the right of the line $\overrightarrow{A,B}$). Moreover, using (c) and the fact that $c_2 > 1, \tan\gamma = c_1/(c_2 - 1) > 0$, we obtain

$$\gamma \in \left(0, \frac{\pi}{2}\right). \qquad (5.39)$$

Note that $\mid C - B \mid \geq s$ and by virtue of (5.37), (5.38), (5.39)

$$\sin(\pi - \varphi) < \sin\gamma = \frac{c_1}{\mid C - B \mid} < \frac{t}{s}. \qquad (5.40)$$

On the other hand, in the case $t \in (0, \min\{1/\sqrt{2};\ \sin(\pi - \varphi)\} s)$, we have $t/s \leq \sin(\pi - \varphi)$, and consequently condition (5.40), and hence condition (5.36), cannot hold; i.e., $C, \bar{C}$ belong to the same part, and (5.35) holds. □

## 5.3. On the nature of the dependence of the angle between vectors on perturbations of the endpoints of these vectors

In this section, it will be shown that the change in the angle between vectors is of order $O(\tau)$ as $\tau \to 0+$, where $\tau$ is the magnitude of the maximum deviation of the corresponding endpoints of these vectors.

Some notation and statements will be needed in what follows.

Let $A, B, C \in \mathbb{R}^2, A \neq B$. Denote $|C| = \sqrt{c_1^2 + c_2^2}, \|C\| = \max\{|c_1|, |c_2|\}$ (where $C = (c_1, c_2)$ ),

$$W_{(A,B)} = [E^{(1)} \ \vdots \ E^{(2)}], E^{(2)} = \frac{B-A}{|B-A|} = (e_1^{(2)}, e_2^{(2)}), E^{(1)} = (e_2^{(2)}, -e_1^{(2)}),$$

i.e., $e_1^{(1)} = e_2^{(2)}$, $e_2^{(1)} = -e_1^{(2)}$. Then $W_{(A,B)}$ is an orthogonal matrix, $W_{(A,B)}^{-1} = W^{\top}{}_{(A,B)}, \left|W_{(A,B)}\right| = 1$,

$$A = \frac{|B-A|}{2} W_{(A,B)} \begin{bmatrix} 0 \\ -1 \end{bmatrix} + \frac{1}{2}(A+B) = -\frac{1}{2}(B-A) + \frac{1}{2}(A+B) = A, \ (5.41)$$

$$B = \frac{|B-A|}{2} W_{(A,B)} \begin{bmatrix} 0 \\ 1 \end{bmatrix} + \frac{1}{2}(A+B) = \frac{1}{2}(B-A) + \frac{1}{2}(A+B) = B. \qquad (5.42)$$

Find the unique point $\bar{C} \in \mathbb{R}^2$ satisfying $C = \frac{|B-A|}{2} W_{(A,B)} \bar{C} + \frac{1}{2}(A+B)$:

$$\bar{C} = \frac{2}{|B-A|} W_{(A,B)}^{\top} \left[C - \frac{1}{2}(A+B)\right]. \qquad (5.43)$$

Correspondingly,

$$\begin{bmatrix}0\\-1\end{bmatrix} = \tfrac{2}{|B-A|} W_{(A,B)}^{\top} \left[A - \tfrac{1}{2}(A+B)\right], \begin{bmatrix}0\\1\end{bmatrix} = \tfrac{2}{|B-A|} W_{(A,B)}^{\top} \left[B - \tfrac{1}{2}(A+B)\right].$$

Moreover, for the angle $\alpha = \angle(B-C, C-A)$ we have:

$$\alpha = \angle(B-C, C-A) < 0\ (>0) \Leftrightarrow \angle\left(\begin{bmatrix}0\\1\end{bmatrix} - \bar{C}, \bar{C} - \begin{bmatrix}0\\-1\end{bmatrix}\right) < 0 (>0).$$

Indeed (considering the first case for definiteness), using (5.41)–(5.43) we obtain:

$$\angle(B-C, C-A) < 0 \Leftrightarrow \begin{vmatrix}B-C\\C-A\end{vmatrix} > 0 \Leftrightarrow \begin{vmatrix}B-C\\C-A\end{vmatrix} =$$

$$= \tfrac{|B-A|}{2} \left|W_{(A,B)}\right| \cdot \begin{vmatrix}(0,1)-\bar{C}\\ \bar{C}-(0,-1)\end{vmatrix} > 0 \Leftrightarrow \begin{vmatrix}(0,1)-\bar{C}\\ \bar{C}-(0,-1)\end{vmatrix} > 0 \Leftrightarrow$$

$$\Leftrightarrow \angle\left(\begin{bmatrix}0\\1\end{bmatrix} - \bar{C}, \bar{C} - \begin{bmatrix}0\\-1\end{bmatrix}\right) < 0. \qquad (5.44)$$

At the same time, by the orthogonality of $W_{(A,B)}$ we obtain:

$$\cos\angle(B-C, C-A) = \frac{\langle B-C, C-A\rangle}{|B-C|\cdot|C-A|} =$$

$$= \frac{\left(\frac{|B-A|}{2}\right)^2 \langle W_{(A,B)}\left(\begin{bmatrix}0\\1\end{bmatrix} - \bar{C}\right), W_{(A,B)}\left(\bar{C} - \begin{bmatrix}0\\-1\end{bmatrix}\right)\rangle}{\left(\frac{|B-A|}{2}\right)^2 \left|W_{(A,B)}\left(\begin{bmatrix}0\\1\end{bmatrix} - \bar{C}\right)\right| \cdot \left|W_{(A,B)}\left(\bar{C} - \begin{bmatrix}0\\-1\end{bmatrix}\right)\right|} =$$

$$= \frac{\langle\left(\begin{bmatrix}0\\1\end{bmatrix} - \bar{C}\right), \left(\bar{C} - \begin{bmatrix}0\\-1\end{bmatrix}\right)\rangle}{\left|\left(\begin{bmatrix}0\\1\end{bmatrix} - \bar{C}\right)\right| \cdot \left|\left(\bar{C} - \begin{bmatrix}0\\-1\end{bmatrix}\right)\right|} = \cos\angle\left(\begin{bmatrix}0\\1\end{bmatrix} - \bar{C}, \bar{C} - \begin{bmatrix}0\\-1\end{bmatrix}\right). \qquad (5.45)$$

From (5.44), (5.45) it follows (see Section 5.2) that

$$\alpha = \angle(B-C, C-A) = \Psi(\bar{C}) = \angle\left(\begin{bmatrix}0\\1\end{bmatrix} - \bar{C}, \bar{C} - \begin{bmatrix}0\\-1\end{bmatrix}\right). \qquad (5.46)$$

Together with the points $A, B, C$, we shall also consider points $\tilde{A}, \tilde{B}, \tilde{C} \in \mathbb{R}^2$ such that

$$\left\|A - \tilde{A}\right\| \leq \tau, \left\|B - \tilde{B}\right\| \leq \tau, \left\|C - \tilde{C}\right\| \leq \tau, \qquad (5.47)$$

where $\tau > 0$, and also the angle $\tilde{\alpha} = \angle(\tilde{B} - \tilde{C}, \tilde{C} - \tilde{A})$. The aim of this section is to prove that in the case $C \neq A$, $C \neq B$ we have:

$$|\tilde{\alpha} - \alpha| = O(\tau) \text{ при } \tau \to 0+. \tag{5.48}$$

We shall use the matrix norm for $M = \begin{bmatrix} m_{11} & m_{12} \\ m_{21} & m_{22} \end{bmatrix}$:

$$\|M\| = \max\{|m_{11}| + |m_{12}|, |m_{21}| + |m_{22}|\}.$$

It is consistent with the vector norm $\|A\| = \max\{|a_1|, |a_2|\}$ for $A = (a_1, a_2) \in \mathbb{R}^2$:

$$\|MA\| \leq \|M\| \cdot \|A\|.$$

We shall also need the following inequality:

$$\|M\| = \left\| \begin{bmatrix} m_{11} & m_{12} \\ m_{21} & m_{22} \end{bmatrix} \right\| = \max\{|m_{11}| + |m_{12}|, |m_{21}| + |m_{22}|\} \leq$$
$$\leq \max\{2\max\{|m_{11}|, |m_{12}|\}, 2\max\{|m_{21}|, |m_{22}|\}\} =$$
$$= 2\max\{\max\{|m_{11}|, |m_{12}|\}, \max\{|m_{21}|, |m_{22}|\}\} =$$
$$= 2\max\{|m_{11}|, |m_{12}|, |m_{21}|, |m_{22}|\}. \tag{5.49}$$

We first prove the following statements.

**Statement 5.6.** *Let* $A, B, C, \tilde{A}, \tilde{B}, \tilde{C} \in \mathbb{R}^2$, $d, s, \tau > 0$, $s > 2\sqrt{2}\,\tau$, *let (5.47) hold, and also*

$$d \geq |A - B| \geq s,\ d \geq |B - C| \geq s,\ d \geq |C - A| \geq s \tag{5.50}$$

(e.g., $d = \max\{|A - B|, |B - C|, |C - A|\}$, $s = \min\{|A - B|, |B - C|, |C - A|\}$). *Let further*

$$W = W_{(A,B)} = [E^{(1)} \ \vdots\ E^{(2)}],\ \widetilde{W} = W_{(\tilde{A},\tilde{B})} = [\tilde{E}^{(1)} \ \vdots\ \tilde{E}^{(2)}],$$
$$b = |A - B|,\ \tilde{b} = |\tilde{A} - \tilde{B}|.$$

*Then*

$$|b - \tilde{b}| \leq 2\sqrt{2}\tau,\ \|\widetilde{W}^{\mathrm{T}} - W^{\mathrm{T}}\| \leq \frac{4\tau}{s} + (d + 2\tau)\frac{4\sqrt{2}\tau}{s\cdot(s - 2\sqrt{2}\tau)} = O(\tau)$$

*as* $\tau \to 0+$.

***Proof.*** Note that

$$|b - \tilde{b}| = \left||A - B| - |\tilde{A} - \tilde{B}|\right| \le \left||A - B| - |A - \tilde{B}|\right| +$$

$$+ \left||A - \tilde{B}| - |\tilde{A} - \tilde{B}|\right| \le |B - \tilde{B}| + |A - \tilde{A}| \le$$

$$\le \sqrt{2}\|B - \tilde{B}\| + \sqrt{2}\|A - \tilde{A}\| \le 2\sqrt{2}\tau. \qquad (5.51)$$

Note further that by virtue of (5.49), (5.51) and the equality $\|\tilde{E}^{(1)} - E^{(1)}\| = \|\tilde{E}^{(2)} - E^{(2)}\|$, we have:

$$\|\widetilde{W}^{\mathrm{T}} - W^{\mathrm{T}}\| = \left\|\begin{matrix}\tilde{E}^{(1)} - E^{(1)} \\ \tilde{E}^{(2)} - E^{(2)}\end{matrix}\right\| \le 2\max\{\|\tilde{E}^{(1)} - E^{(1)}\|, \|\tilde{E}^{(2)} - E^{(2)}\|\} =$$

$$= 2\|\tilde{E}^{(2)} - E^{(2)}\| = 2\left\|\frac{B - A}{|B - A|} - \frac{\tilde{B} - \tilde{A}}{|\tilde{B} - \tilde{A}|}\right\| \le$$

$$\le 2\left\|\frac{B - A}{|B - A|} - \frac{\tilde{B} - \tilde{A}}{|B - A|}\right\| + 2\left\|\frac{\tilde{B} - \tilde{A}}{|B - A|} - \frac{\tilde{B} - \tilde{A}}{|\tilde{B} - \tilde{A}|}\right\| \le$$

$$\le \frac{4\tau}{|A - B|} + 2\|\tilde{B} - \tilde{A}\|\left(\frac{\left||\tilde{B} - \tilde{A}| - |B - A|\right|}{|\tilde{B} - \tilde{A}| \cdot |B - A|}\right) \le$$

$$\le \frac{4\tau}{|A - B|} + (|A - B| + 2\tau)\frac{4\sqrt{2}\tau}{|A - B| \cdot \left(|A - B| - 2\sqrt{2}\tau\right)} \le$$

$$\le \frac{4\tau}{s} + (d + 2\tau)\frac{4\sqrt{2}\tau}{s \cdot \left(s - 2\sqrt{2}\tau\right)}. \ \square$$

(Here we used:

$$\| \tilde{B} - \tilde{A} \| \le \| B - A \| + 2\tau \le | A - B | + 2\tau,$$

$$| A - \tilde{A} | \le \sqrt{2} \ \| A - \tilde{A} \| \le \sqrt{2}\,\tau, | \tilde{B} - \tilde{A} | \ge | A - B | - | A - \tilde{A} | - | B - \tilde{B} |$$

$$\ge | A - B | - 2\sqrt{2}\,\tau).$$

**Remark 5.2.** Note that for any $2 \times 2$ orthogonal matrix $W$ we have $\| W \| \le \sqrt{2}$.

**Statement 5.7.** *Let the conditions of Statement 5.6 hold. Let further*

$b = |A - B|, \tilde{b} = |\tilde{A} - \tilde{B}|,$

$$\bar{C} = \frac{2}{|B-A|} W_{(A,B)}^{\top} \left[C - \frac{1}{2}(A+B)\right], \hat{C} = \frac{2}{|\tilde{B}-\tilde{A}|} \widetilde{W}^{\mathrm{T}} \left[\tilde{C} - \frac{1}{2}(\tilde{A}+\tilde{B})\right]. \quad (5.52)$$

*Then as* $\tau \to 0+$

$$\|\bar{C} - \hat{C}\| \le \frac{4\tau d}{s}\left[1 + \frac{d}{s-2\tau} + \frac{2\sqrt{2}}{s-2\tau}\right] + \frac{8\tau^2}{s}\left[\frac{d}{s-2\tau} + \sqrt{2} + \frac{2\sqrt{2}}{s-2\tau}\right] = O(\tau). \quad (5.53)$$

***Proof.*** Note that

$|\tilde{b}| = |\tilde{A} - \tilde{B}| \ge |A - B| - |A - \tilde{A}| - |B - \tilde{B}| \ge |A - B| - |A - \tilde{A}| - |B - \tilde{B}| \ge$

$\ge |A - B| - 2\sqrt{2}\tau \ge s - 2\sqrt{2}\tau,$

$\left\|C - \frac{1}{2}(A+B) - \left[\tilde{C} - \frac{1}{2}(\tilde{A}+\tilde{B})\right]\right\| \le 2\tau,$

$\left\|C - \frac{1}{2}(A+B)\right\| = \left\|\frac{1}{2}(C-A) + \frac{1}{2}(C-B)\right\| \le d, \left\|\tilde{C} - \frac{1}{2}(\tilde{A}+\tilde{B})\right\| \le d + 2\tau,$

Hence (using (5.51), Statement 5.6 and Remark 5.2)

$$\|\hat{C} - \bar{C}\| = \left\|\frac{2}{b} W^{\mathrm{T}}\left[C - \frac{1}{2}(A+B)\right] - \frac{2}{\tilde{b}} \widetilde{W}^{\mathrm{T}}\left[\tilde{C} - \frac{1}{2}(\tilde{A}+\tilde{B})\right]\right\| \le$$

$$\le \left\|\frac{2}{b} W^{\mathrm{T}}\left[C - \frac{1}{2}(A+B)\right] - \frac{2}{b} \widetilde{W}^{\mathrm{T}}\left[C - \frac{1}{2}(A+B)\right]\right\| +$$

$$+ \left\|\frac{2}{b} \widetilde{W}^{\mathrm{T}}\left[C - \frac{1}{2}(A+B)\right] - \frac{2}{b} \widetilde{W}^{\mathrm{T}}\left[\tilde{C} - \frac{1}{2}(\tilde{A}+\tilde{B})\right]\right\| +$$

$$+ \left\|\frac{2}{b} \widetilde{W}^{\mathrm{T}}\left[\tilde{C} - \frac{1}{2}(\tilde{A}+\tilde{B})\right] - \frac{2}{\tilde{b}} \widetilde{W}^{\mathrm{T}}\left[\tilde{C} - \frac{1}{2}(\tilde{A}+\tilde{B})\right]\right\| \le$$

$$\le \frac{2}{b} \cdot \|\widetilde{W}^{\mathrm{T}} - W^{\mathrm{T}}\| \cdot \left\|C - \frac{1}{2}(A+B)\right\| + \frac{2}{b} \cdot \|\widetilde{W}^{\mathrm{T}}\| \cdot 2\tau +$$

$$+2\|\widetilde{W}^{\mathrm{T}}\| \cdot \left\|\tilde{C} - \frac{1}{2}(\tilde{A}+\tilde{B})\right\| \cdot \left|\frac{1}{b} - \frac{1}{\tilde{b}}\right| \le$$

$$\le \frac{2d}{s} \cdot \left[\frac{4\tau}{s} + (d+2\tau)\frac{4\sqrt{2}\tau}{s(s-2\sqrt{2}\tau)}\right] + \frac{2\sqrt{2}}{s} \cdot 2\tau + 2\sqrt{2}(d+2\tau)\frac{|b - \tilde{b}|}{s(s-2\sqrt{2}\tau)} \le$$

$$\le \frac{2}{s}\left[\frac{4\tau}{s} + (d+2\tau)\frac{4\sqrt{2}\tau}{s(s-2\sqrt{2}\tau)} + \frac{2\sqrt{2}}{s} \cdot 2\tau + 2\sqrt{2}(d+2\tau)\frac{2\sqrt{2}\tau}{s(s-2\sqrt{2}\tau)}\right] \le$$

$$\le \frac{8\tau}{s^2}\left[\tau + (d+2\tau)\frac{\sqrt{2}}{s-2\sqrt{2}\tau} + \sqrt{2} + 2(d+2\tau)\frac{1}{s-2\sqrt{2}\tau}\right] \le$$

$$\leq \frac{8\tau d}{s^2}\left[\sqrt{2}+\frac{\sqrt{2}d}{s-2\sqrt{2}\tau}+\frac{2d}{s-2\sqrt{2}\tau}\right]+\frac{8\tau^2}{s^2}\left[1+\frac{2\sqrt{2}+4}{s-2\sqrt{2}\tau}\right]=O(\tau). \square$$

**Statement 5.8.** *Let the conditions of Statements 5.6 and 5.7 hold. Then (5.48) holds.*

***Proof.*** Note that similarly to (5.46) we have

$$\tilde{\alpha}=\angle\left(\tilde{B}-\tilde{C},\tilde{C}-\tilde{A}\right)=\Psi\left(\hat{C}\right)=\angle\left(\begin{bmatrix}0\\1\end{bmatrix}-\hat{C},\hat{C}-\begin{bmatrix}0\\-1\end{bmatrix}\right), \tag{5.54}$$

and consequently,

$$|\tilde{\alpha}-\alpha|=\left|\Psi(\bar{C})-\Psi(\hat{C})\right|,$$

where $\bar{C},\hat{C}$ satisfy (5.52). By Statement 5.7 (see (5.53)) there exists a number $\tau_0>0$ such that

$$\forall\tau\in(0,\tau_0]\quad\left\|\bar{C}-\hat{C}\right\|\leq\min\left\{\frac{1}{2\sqrt{2}},\sin(\pi-\varphi)\right\}s.$$

But then by Statement 5.5,

$$\forall\tau\in(0,\tau_0]\ |\tilde{\alpha}-\alpha|=\left|\Psi(\bar{C})-\Psi(\hat{C})\right|\leq L_\Psi\left(\frac{s}{2}\right)\left\|C-\hat{C}\right\|,$$

where by Statement 5.7 $\left\|C-\hat{C}\right\|=O(\tau)$ as $\tau\to 0+$, which implies the validity of (5.48). $\square$

### 5.4. Approximation of the set $\boldsymbol{S}^{(n)}(\boldsymbol{A},\boldsymbol{B},\boldsymbol{\varphi})$ by a discrete set $\boldsymbol{S}^{(n)}(\boldsymbol{Q}^{(\tau)},\boldsymbol{A},\boldsymbol{B},\boldsymbol{\varphi})$

Let $A,B\in\mathbb{R}^2, A\neq B, n\geq 1$, and let $\boldsymbol{S}^{(n)}(A,B,\varphi)$ be the set of all possible sequences $\left(B^{(1)},\dots,B^{(n)}\right)$ of interior points of a polygonal line $\boldsymbol{L}$ satisfying conditions (0.3), (0.4), where $A=B^{(0)}, B=B^{(n+1)}$.
Let $n\varphi<\pi$. Then (see Theorem 1)

$$\boldsymbol{S}^{(n)}(A,B,\varphi)\subseteq[\boldsymbol{S}(A,B,n\varphi)]^n,$$

where $\boldsymbol{S}(A,B,n\varphi)$ is a bounded set; i.e., there exist $q_i,p_i\in\mathbb{R}, q_i<p_i, i=1,2$, such that

$$\boldsymbol{S}(A,B,n\varphi)\subset\boldsymbol{Q}=\{C=(c_1,c_2)\in\mathbb{R}^2|q_i\leq c_i\leq p_i, i=1,2\}, \tag{5.55}$$

where $\boldsymbol{Q}$ is a **coordinate rectangle**.

Using $\boldsymbol{Q}$, we also consider the discrete case. Let $\tau_1, \tau_2 > 0$. Consider the finite set

$$\boldsymbol{Q}^{(\tau_1,\tau_2)} = \{C = (c_1, c_2) \in \mathbb{R}^2 | c_1 = j_1\tau_1, c_2 = j_2\tau_2, q_i \le c_i \le p_i, i = 1,2,$$
$$j_1, j_2 \in \mathbb{Z}\}, \tag{5.56}$$

and for simplicity let

$$\tau_1 = \tau_2 = \tau,\ \boldsymbol{Q}^{(\tau)} = \boldsymbol{Q}^{(\tau_1,\tau_2)}. \tag{5.57}$$

Note that if $\tau \le p_i - q_i, i = 1,2,$ then

$$\forall \left(B^{(1)}, \dots, B^{(n)}\right) \in \boldsymbol{Q}^n\ \exists \left(\tilde{B}^{(1)}, \dots, \tilde{B}^{(n)}\right) \in \left[\boldsymbol{Q}^{(\tau)}\right]^n:$$
$$\left\|\tilde{B}^{(i)} - B^{(i)}\right\| \le \tau, i = 1, \dots, n, \tag{5.58}$$

and if $q_i, p_i$ are integer multiples of $\tau$, then

$$\forall \left(B^{(1)}, \dots, B^{(n)}\right) \in \boldsymbol{Q}^n\quad \exists \left(\tilde{B}^{(1)}, \dots, \tilde{B}^{(n)}\right) \in \left[\boldsymbol{Q}^{(\tau)}\right]^n:$$
$$\left\|\tilde{B}^{(i)} - B^{(i)}\right\| \le \frac{\tau}{2}, i = 1, \dots, n.$$

Let $W = \left(B^{(1)}, \dots, B^{(n)}\right) \in \boldsymbol{S}^{(n)}(A, B, \varphi) \subseteq \boldsymbol{Q}^n$. Then from (5.58) we obtain that for sufficiently small $\tau > 0$ one can find $\left(\tilde{B}^{(1)}, \dots, \tilde{B}^{(n)}\right) \in \left[\boldsymbol{Q}^{(\tau)}\right]^n$: such that

$$\left\|\tilde{B}^{(i)} - B^{(i)}\right\| \le \tau, i = 1, \dots, n.$$

Moreover, from $W = \left(B^{(1)}, \dots, B^{(n)}\right) \in \boldsymbol{S}^{(n)}(A, B, \varphi)$ it follows that

$$\alpha_i = \angle\left(B^{(i+1)} - B^{(i)}, B^{(i)} - B^{(i-1)}\right) \in [-\varphi, \varphi], i = 1, \dots, n,$$

and for some $s > 0$ (see condition 1 of Theorem 1) we have

$$s \le \left|B^{(i)} - B^{(j)}\right|, i, j = 0,1, \dots, n+1, i \ne j, \tag{5.59}$$

Two cases are possible:

1) For some $\omega > 0$ we have

$$\alpha_i = \angle(B^{(i+1)} - B^{(i)}, B^{(i)} - B^{(i-1)})) \in [-\varphi + \omega, \varphi - \omega]. \tag{5.60}$$

Then, using the fact that

$$s \le \left|B^{(i)} - B^{(i-1)}\right| \le d, i = 1, \dots, n+1, \tag{5.61}$$

where $d \le \operatorname{diam}_e \boldsymbol{S}(A, B, n\varphi) \le \sqrt{(p_1 - q_1)^2 + (p_2 - q_2)^2}$, we obtain by Statement 5.8 that there exists a number $\tau_0 > 0$ such that for all $\tau \in (0, \tau_0]$ the following holds:

$$\|\widetilde{W} - W\| \le \tau \Rightarrow |\alpha_i - \tilde{\alpha}_i| \le \omega, \tag{5.62}$$

where $\tilde{\alpha}_i = \angle(\tilde{B}^{(i+1)} - \tilde{B}^{(i)}, \tilde{B}^{(i)} - \tilde{B}^{(i-1)}), i = 1, \dots, n.$

But then for all $\tau \in (0, \tau_0]$ in $[\boldsymbol{Q}^{(\tau)}]^n$ there exists $\widetilde{W} = (\tilde{B}^{(1)}, \dots, \tilde{B}^{(n)}) \in [\boldsymbol{Q}^{(\tau)}]^n$ with $\| \widetilde{W} - W \| \le \tau$ (see (5.58)), and by virtue of (5.60), (5.62) $\widetilde{W}$ belongs to $\boldsymbol{S}^{(n)}(A, B, \varphi)$. Thus, taking into account that $\widetilde{W} = (\tilde{B}^{(1)}, \dots, \tilde{B}^{(n)}) \in [\boldsymbol{Q}^{(\tau)}]^n$, we have

$$\forall \tau \in (0, \tau_0] \quad \left[\widetilde{W} = \left(\tilde{B}^{(1)}, \dots , \tilde{B}^{(n)}\right) \in \left[\boldsymbol{Q}^{(\tau)}\right]^n, \|\widetilde{W} - W\| \le \tau\right] \Rightarrow$$

$$\Rightarrow \widetilde{W} \in \boldsymbol{S}^{(n)}(A, B, \varphi) \cap \left[\boldsymbol{Q}^{(\tau)}\right]^n,$$

i.e., in this case, for every $W = (B^{(1)}, \dots, B^{(n)}) \in \boldsymbol{S}^{(n)}(A, B, \varphi)$ and for sufficiently small $\tau > 0$, any sequence of turning points

$$\widetilde{W} = (\tilde{B}^{(1)}, \dots, \tilde{B}^{(n)}) \in \boldsymbol{S}^{(n)}(A, B, \varphi)$$

satisfying $\| \widetilde{W} - W \| \le \tau$ (such a sequence exists by (5.58)) will be found by the discrete search method announced in [2]. In this method, using Theorem 2, we sequentially iterate through the admissible turning points that belong simultaneously to the finite set $\boldsymbol{Q}^{(\tau)}$ and to the current member in the direct product on the right-hand side of condition (0.6).

2) In the second case, for some (at least one) $i \in \{1, \dots, n\}$ we have

$$\alpha_i = \angle\left(B^{(i+1)} - B^{(i)}, B^{(i)} - B^{(i-1)}\right) = \varphi.$$

In this case, as shown in Statement 5.3, there exist numbers

$$\omega = \omega(s, d, \varphi, n) > 0, \theta = \theta(s, d, \varphi, n) \ge 5,$$

such that $\forall t \in (0, \frac{s}{\theta}]$ (i.e., $t > 0$, $\frac{s}{t} \ge \theta$) one can construct (by simple algorithmically realizable formulas) a sequence

$$\widetilde{W}(t) = \left(\tilde{B}^{(1)}(t), \dots, \tilde{B}^{(n)}(t)\right) \in \boldsymbol{S}^{(n)}(A, B, \varphi)$$

of turning points with

$$\tilde{B}^{(i)}(t) \in \left[B^{(i)}, B^{(i+1)}\right], i = 1, \dots, n$$

(i.e., the points $\tilde{B}^{(1)}(t), \dots, \tilde{B}^{(n)}(t)$ belong to the polygonal line

$$\boldsymbol{L} = \left[B^{(0)}, B^{(1)}\right] \cup \left[B^{(1)}, B^{(2)}\right] \cup \dots \cup \left[B^{(n)}, B^{(n+1)}\right]$$

), such that

$$\widetilde{W}(t) = \left(\tilde{B}^{(1)}(t), \dots, \tilde{B}^{(n)}(t)\right) \in \boldsymbol{S}^{(n, s-t\kappa)}(A, B, \varphi - \omega t),$$

$$\left\|W - \widetilde{W}(t)\right\| \leq t\frac{\kappa}{2}, \qquad (5.63)$$

where $\kappa = \frac{d}{s}$.

Thus, for every $t \in (0, \frac{s}{\theta}]$ the set $\widetilde{W}(t) = (\tilde{B}^{(1)}(t), \dots, \tilde{B}^{(n)}(t))$ satisfies case 1. But then (as shown in case 1) for any fixed $t \in (0, \frac{s}{\theta}]$ there exists a number $\tau_0(t) > 0$ such that, taking (5.63) into account, we have

$$\forall \tau \in (0, \tau_0(t)]\ \exists \breve{W}(\tau) = \left(\breve{B}^{(1)}, \dots, \breve{B}^{(n)}\right) \in \boldsymbol{S}^{(n)}(A, B, \varphi) \cap \left[\boldsymbol{Q}^{(\tau)}\right]^n:$$

$$\left\| \breve{W}(\tau) - \widetilde{W}(t) \right\| \leq \tau,\ \left\|W - \breve{W}(\tau)\right\| \leq t\frac{\kappa}{2} + \tau\,. \qquad (5.64)$$

From (5.64), by the arbitrariness of $W = \left(B^{(1)}, \dots, B^{(n)}\right) \in \boldsymbol{S}^{(n)}(A, B, \varphi)$ and the arbitrariness of the choice of $t \in (0, \frac{s}{\theta}]$, it follows that

$$\lim_{\tau \to 0+} h_m\left(\boldsymbol{S}^{(n)}(A, B, \varphi), \boldsymbol{S}^{(n)}(A, B, \varphi) \cap \left[\boldsymbol{Q}^{(\tau)}\right]^n\right) = 0, \qquad (5.65)$$

i.e., in this case as well, the convergence of the discrete method (announced in [2]) to the solution set in the continuous case is proved.

Thus, we have proved

**Statement 5.9.** *Let* $\varphi \in (0, \pi), n \geq 1, n\varphi < \pi, A, B \in \mathbb{R}^2, A \neq B, \tau > 0,$ *and let* $\boldsymbol{Q}, \boldsymbol{Q}^{\tau}$ *satisfy conditions (5.55)–(5.57). Then (5.65) holds.*

Consider also a more general case. Suppose now that more general conditions hold for the set $\boldsymbol{Q}$; namely, instead of the condition $\boldsymbol{S}(A,B,n\varphi) \subset \boldsymbol{Q} \subset \mathbb{R}^2$, let $\boldsymbol{Q}$ be a convex set and assume

$$\boldsymbol{S}(A,B,n\varphi) \cap \boldsymbol{Q} \neq \emptyset, \tag{5.66}$$

and we have a set $\boldsymbol{Q}^{(\tau)}, \tau > 0$, satisfying

$$\boldsymbol{Q}^{(\tau)} \subset \boldsymbol{Q},\ h_m\big(\boldsymbol{Q}^{(\tau)}, \boldsymbol{Q}\big) \leq \tau. \tag{5.67}$$

**Statement 5.10.** *Let* $\varphi \in (0,\ \pi)$*,* $n \geq 1$*,* $n\varphi < \pi$*,* $A, B \in \mathbb{R}^2$*,* $A \neq B$*,* $\tau > 0$*,* $\boldsymbol{Q}, \boldsymbol{Q}^\tau$ *satisfy (5.66), (5.67),* $\boldsymbol{S}^{(n)}(A,B,\varphi) \cap \boldsymbol{Q}^n \neq \emptyset$ *and let* $\boldsymbol{Q}$ *be convex. Then*

$$\lim_{\tau \to 0+} h_m\Big(\boldsymbol{S}^{(n)}(A,B,\varphi) \cap \boldsymbol{Q}^n, \boldsymbol{S}^{(n)}(A,B,\varphi) \cap \big[\boldsymbol{Q}^{(\tau)}\big]^n\Big) = 0.$$

***Proof.*** One direction is obvious, since

$$\boldsymbol{S}^{(n)}(A,B,\varphi) \cap \big[\boldsymbol{Q}^{(\tau)}\big]^n \subseteq \boldsymbol{S}^{(n)}(A,B,\varphi) \cap \boldsymbol{Q}^n.$$

The opposite direction is proved similarly to the proof of Statement 3.1. Here we use the convexity of $\boldsymbol{Q}^n$ to show that $\widetilde{W}(t) \in \boldsymbol{Q}^n$. Indeed, in case 2 we use $\tilde{B}^{(i)}(t) \in \big[B^{(i)}, B^{(i+1)}\big]$, which, by convexity of $\boldsymbol{Q}$, implies $\tilde{B}^{(i)}(t) \in \boldsymbol{Q}$. (Closedness of $\boldsymbol{Q}$ is not required; we only need $h_m\big(\boldsymbol{Q}^{(\tau)}, \boldsymbol{Q}\big) \leq \tau$).

## 6. Conclusion.

Let us list the main results of the paper.

1. The problem of finite approximation of the set of piecewise-linear paths (polygonal lines) connecting two given points $A, B \in \mathbb{R}^2$ with the possibility of making $n \geq 1$ turns at some (corner) points $B^{(1)}, \dots, B^{(n)} \in \mathbb{R}^2$ (each turn angle bounded in absolute value by a given number $\varphi \in \left(0, \frac{\pi}{2}\right)$) is considered. Under the condition $n\varphi < \pi$, a set is described to which all interior vertices of such a polygonal line belong (Theorem 1). It is proved that for any point $B^{(1)}$ from this set, there exists a polygonal line satisfying the given conditions (Lemma 1). Based on these results, an explicit formula (Theorem 2) is obtained, describing the set $\mathbf{S}^{(n)}(A, B, \varphi)$ of all admissible sequences $\left(B^{(1)}, \dots, B^{(n)}\right)$ of angular points.

2. The resulting formula (Theorem 2) is used to construct a finite set of sequences $\left(B^{(1)}, \dots, B^{(n)}\right) \in \mathbf{S}^{(n)}(A, B, \varphi)$ that approximates the said set $\mathbf{S}^{(n)}(A, B, \varphi)$. This specified set is used to describe a practically implementable algorithm for constructing a finite set of admissible tuples $\left(B^{(1)}, \dots, B^{(n)}\right)$ of corner points of the polygonal line, each chosen from a discrete set $\boldsymbol{Q}^{(\tau)} \subset \boldsymbol{Q}$, where $\boldsymbol{Q}$ is a rectangle satisfying (5.55), and $\boldsymbol{Q}^{(\tau)}$ is a finite set satisfying (5.5 6)–(5.58), approximating $\boldsymbol{Q}$ with a given accuracy $\tau > 0$.

3. It is proved that the collection of corner-point tuples obtained by the specified algorithm converges (in the Hausdorff metric) to the solution set of the original (continuous) problem $\boldsymbol{S}^{(n)}(A, B, \varphi)$ as $\tau \to 0^{+}$.

4. The obtained algorithm can be used for approximate solution of optimization problems for some objective function on the set of all admissible tuples $(B^{(1)}, \dots, B^{(n)}) \in \boldsymbol{S}^{(n)}(A, B, \varphi)$, that takes into account the cost of traversing the links and the cost of turns.